\documentclass[11pt]{article}
\usepackage{amssymb,amsmath,amsthm,graphicx,color,natbib,enumerate,booktabs}
\usepackage[colorlinks=true, linkcolor=black, citecolor=black, urlcolor=black]{hyperref}
\usepackage[lmargin=1.15in,rmargin=1.15in,bmargin=1.3in]{geometry}

\usepackage{graphicx}
\usepackage{subfigure,curves}

\newcommand{\cip}{\mbox{\,$\perp\!\!\!\perp$\,}}
\newcommand{\whitebox}{\begin{flushright}\vspace{-20pt}$\Box$\end{flushright}}

\newtheorem{theorem}{Theorem}[section]
\newtheorem{lemma}[theorem]{Lemma}
\newtheorem{proposition}[theorem]{Proposition}

\theoremstyle{definition}		

\newtheorem{definition}[theorem]{Definition}
\newtheorem{example}[theorem]{Example}

\definecolor{orange}{rgb}{1,0.5,0}
\definecolor{pink}{rgb}{1,0.75,0.79}
\definecolor{purple}{rgb}{0.65, 0.12, 0.82}
\definecolor{darkgreen}{rgb}{0,.8,0}

\title{Lattices of Graphical Gaussian Models with Symmetries} 

\author{Helene Gehrmann\thanks{\emph{Address for correspondence}. Helene Gehrmann, Department of Statistics, 1 South Parks Road, Oxford OX1 3TG, United Kingdom, 
\texttt{gehrmann@stats.ox.ac.uk}.}\\
Department of Statistics, University of Oxford
}

\begin{document}

\hyphenation{Lau-ritzen}
\pretolerance=0
\setlength{\emergencystretch}{8em}
\hyphenpenalty=1000 

\maketitle

\begin{abstract} In order to make graphical Gaussian models a viable modelling tool when the number of variables outgrows the number of
observations,\citet{lauritzen_sym} introduced model classes which place equality restrictions on concentrations or partial correlations. The models can be represented by vertex and edge coloured graphs. The need for model selection methods makes it imperative to understand the structure of model classes. We identify four model classes that form complete lattices of models with respect to model inclusion, which qualifies them for an Edwards--Havr{\'a}nek model selection procedure \citep{edwards_lattice}. Two classes turn out most suitable for a corresponding model search. We obtain an explicit search algorithm for one of them and provide a model search example for the
other.\\

\noindent \emph{Keywords:} Conditional independence; Covariance selection; Invariance; Model selection; Patterned covariance matrices; Permutation symmetry
\end{abstract}

\section{Introduction}

Graphical models are probabilistic models which use graphs to
represent dependencies between random variables. This article is
concerned with models represented by undirected graphs, in which
each variable corresponds to a vertex and a pair of vertices is
joined by an edge unless the corresponding variables are
conditionally independent given the remaining variables. In addition
to providing a concise form of visualisation of the conditional
independence structure of a model, the graphical representation can
be exploited to make statistical inference computations more
efficient \citep{lauritzen}.

Motivated by the growing need for parsimonious models in modern day
applications, in particular when the number of variables outgrows
the number of observations, in recent years graphical models with
additional equality constraints on model parameters are becoming of
increasing interest, in discrete models \citep{gottard,ramirez} as
well as in multivariate Gaussian models, which are the central
object of interest in this article. First
studies \citep{lauritzen_sym,uhler} show that equality constraints
reduce the minimal number of observations required to ensure
estimability of the model parameters with probability one, which
makes graphical Gaussian models with equality constraints a
promising model class.

Symmetry constraints, induced by distribution invariance under a
permutation group applied to the variable labels, are a special
instance of equality constraints and have a long history for the
Gaussian distribution before the advent of graphical
models \citep{wilks, votaw, olkin, olkin72, andersson, jensen}. First
studies of models combining symmetry constraints with conditional
independence relations are given in \citet{hylleberg, andersen,
madsen}. The models we study have been introduced
in \citet{lauritzen_sym} and contain the models in \citet{hylleberg} as a special case. The types of restrictions are:
equality between specified elements of the concentration matrix
(RCON) and equality between specified partial correlations (RCOR).
The models can be represented by vertex and edge coloured graphs,
where parameters associated with equally coloured vertices or edges
are restricted to being identical.



In order for RCON and RCOR models to become widely applicable in
practice model selection methods need to be developed, which
motivates the study of model structures. This is the main objective
of this article. As both model types RCON and RCOR models form
complete lattices, both qualify for the Edwards--Havr{\'a}nek model
selection procedure for lattices \citep{edwards_lattice}. However due
to the large number of models it is more feasible to search, at
least initially, in suitable subsets of the model space.
Particularly favourable are subsets of models in which equality
constraints are placed in a pattern which makes them more readily
interpretable.

Four model classes with desirable statistical properties which
express themselves in regularity of graph colouring have been
previously identified in the literature: The most restrictive is
given by \emph{graphical symmetry models} studied
in \citet{hylleberg}, also appearing in \citet{lauritzen_sym} under
the name \emph{RCOP models}. The corresponding graph colourings are
given by vertex and edge orbits of a permutation group acting on the
variable labels and we therefore term them
\emph{permutation-generated}. Colourings representing models which
place the same equality constraints on the concentrations and
partial correlations were termed \emph{edge regular}
in \citet{gehrmann}. Two further model types, ensuring estimability
of a non-zero model mean subject to equality constraints, were
introduced in \citet{gehrmann}, the colourings representing them
termed \emph{vertex regular} and \emph{regular} respectively.

The main results presented in this article are that each of the
model classes forms a complete lattice of models and the
identification of their meet and join operations. The former is
established by showing that each model class is stable under model
intersection, which gives the shared meet operation, and by
demonstrating that whenever a model does not fall inside a given
class there is a unique smallest larger model, or supremum, which
does, giving the distinct join operations. The found lattice
structure qualifies each model class for an Edwards--Havr{\'a}nek
model search, giving a first model selection procedure for RCON and
RCOR models.

We focus on models represented by edge regular and
permutation-generated colourings as their structure is generally
more tractable and their constraints more readily interpretable. We
present an Edwards--Havr{\'a}nek model selection algorithm for
models with edge regular colourings and illustrate it by means of an
example with five variables, with a very encouraging performance. We
further provide an example of a model search within models with
permutation-generated colourings with four variables, commonly known
as \emph{Fret's heads} \citep{frets,mardia}.

\section{Preliminaries and Notation}

\subsection{Notation}

Let $G=(V,E)$ be an undirected uncoloured graph with vertex set $V$
and edge set $E$. For a $|V| \times |V|$ matrix
$A=(a_{\alpha\beta})$, $A(G)$ shall denote the matrix defined by
$A(G)_{\alpha\beta} = 0$ whenever there is no edge between $\alpha$
and $\beta$ in $G$ for $\alpha\not=\beta$, and $A(G)_{\alpha\beta} =
a_{\alpha\beta}$ otherwise. For a set of matrices $M$ we let $M^{+}$
denote the set of positive-definite matrices inside $M$.
$\mathcal{S}$ shall denote the set of symmetric matrices, so that
$\mathcal{S}^{+}$ denotes the set of symmetric positive definite
matrices and $\mathcal{S}^{+}(G)$ the set of symmetric positive
definite matrices indexed by $V$ whose $\alpha\beta$-entry is zero
for $\alpha\beta \not\in E$ for $\alpha\not=\beta$. We indicate that
a matrix is symmetric by only writing its elements on the diagonal
and above. Asterisks as matrix entries indicate that the
corresponding entry is unconstrained apart from any restrictions
stated explicitly.

For a discrete set $D$ we let $S(D)$ denote the set symmetric group
acting on $D$, consisting of all permutations of the elements in
$D$. For $F \subseteq S(D)$, $\langle F \rangle$ denotes the group
generated by $F$, containing all permutation which can be expressed
as products of elements in $F$ and their inverses. Permutations are
written in cycle notation, meaning that $\sigma=(i_{1}i_{2}\ldots
i_{r})$ maps $i_{j}$ to $i_{j+1}$ for $1 \leq j < r$ and $i_{r}$ to
$i_{1}$.

For a graph $G=(V,E)$, $\textup{Aut}(G)$ denotes the automorphism
group of $G$, containing all permutations in $S(V)$ which leave $G$
invariant. For a partition $P$ of a set $S$ and $a,b \in S$ we write
$a \equiv b \ (P)$ to denote that $a$ and $b$ lie in the same set in
$P$. For $n \in \mathbb{N}$, we denote sets of the form $\{1, 2,
\ldots, n\}$ by $[n]$.

\subsection{Graphical Gaussian Models} \label{ugm}

A \emph{graphical Gaussian model} is concerned with the distribution
of a multivariate random vector \linebreak $Y=(Y_{\alpha})_{\alpha
\in V}$. Let $G=(V,E)$ be an undirected graph with vertex set $V$
and edge set $E$. Then the graphical Gaussian model represented by
$G$ is given by assuming that $Y$ follows a multivariate Gaussian
$\mathcal{N}_{|V|}(\mu, \Sigma)$ distribution with concentration
matrix $K=\Sigma^{-1} \in \mathcal{S}^{+}(G)$.

The entries in $K=(k_{\alpha\beta})_{\alpha,\beta \in V}$ have a
simple interpretation. The diagonal elements $k_{\alpha\alpha}$ are
reciprocals of the conditional variances given the remaining
variables
\begin{equation*} \label{eq_k_interpretation}
k_{\alpha\alpha} = \text{Var}(Y_{\alpha}|Y_{V\setminus
\{\alpha\}})^{-1}
\end{equation*}
for $\alpha \in V$. The scaled elements of the concentration matrix
\begin{equation} \label{eq:c}
c_{\alpha\beta} =
\frac{k_{\alpha\beta}}{\sqrt{k_{\alpha\alpha}k_{\beta\beta}}}
\end{equation}
for $\alpha,\beta \in V$ are the negative partial correlation
coefficients
\begin{equation} \label{partial_cor}
\rho_{\alpha\beta | V \setminus \{\alpha,\beta\}} =
\frac{\text{Cov}(Y_{\alpha},Y_{\beta} | Y_{V \setminus
\{\alpha,\beta\}})}{\text{Var}(Y_{\alpha} | Y_{V \setminus
\{\alpha\}})^{1/2}\text{Var}(Y_{\beta} | Y_{V \setminus
\{\beta\}})^{1/2}} = - c_{\alpha\beta}
\end{equation}
for $\alpha,\beta \in V$. It follows that
\begin{equation} \label{eq_cip}
\alpha\beta \not\in E \ \ \ \Longleftrightarrow \ \ \
k_{\alpha\beta} = 0 \ \ \ \Longleftrightarrow \ \ \ Y_{\alpha} \cip
Y_{\beta} \mid Y_{V \setminus \{\alpha,\beta\}}
\end{equation}
see e.g.\ Chapter 5 in \citet{lauritzen} for further details.

\subsection{Graph Colouring}

For general graph terminology we refer to \citet{bollobas_book}.
Following \citet{lauritzen_sym}, for $G=(V,E)$ an undirected graph, a
\emph{vertex colouring} of $G$ is a partition $\mathcal{V}=\{V_{1},
\ldots, V_{k}\}$ of $V$, where we refer to $V_{1}, \ldots, V_{k}$ as
the \emph{vertex colour classes}. Similarly, an \emph{edge
colouring} of $G$ is a partition $\mathcal{E}=\{E_{1}, \ldots,
E_{l}\}$ of $E$ into $l$ \emph{edge colour classes} $E_{1}, \ldots,
E_{l}$. A colour class with a single element is \emph{atomic} and a
colour class which is not atomic is \emph{composite}. We let
$\mathcal{G}=(\mathcal{V},\mathcal{E})$ denote the \emph{coloured
graph} with vertex colouring $\mathcal{V}$ and edge colouring
$\mathcal{E}$ and let $(\mathcal{V},\mathcal{E})$ denote its
\emph{graph colouring}. For $\mathcal{V}$ and $\mathcal{E}$ as above
and $u \in \mathcal{V}$, we let $T^{u}$ denote the $|V| \times |V|$
diagonal matrix with $T^{u}_{\alpha\alpha} = 1$ if and only if
$\alpha \in u$ and zero otherwise. Similarly for $u \in
\mathcal{E}$, we let $T^{u}$ be the symmetric $|V| \times |V|$
matrix with $T^{u}_{\alpha\beta} = 1$ if and only if $\alpha\beta
\in u$ and zero otherwise.

In our display of vertex and edge coloured graphs, we indicate the
colour class of a vertex by the number of asterisks we place next to
it. Similarly we indicate the colour class of an edge by dashes.
Vertices and edges which are displayed in black correspond to atomic
colour classes.

\subsection{Lattices} \label{lattice_theory}

A \emph{binary relation} $\rho$ on a set $A$ is a subset of $A
\times A$ with two elements $a , b \in A$ being \emph{in relation}
with respect to $\rho$ if and only if $(a,b) \in \rho$, which we
denote by $a \ \rho \  b$. If $\rho$ is \emph{reflexive} [$a \ \rho
\ a\ \forall a \in A$], \emph{antisymmetric} [$a \ \rho \ b$ and $b
\ \rho \ a \ \Rightarrow \ a = b\ \forall a,b \in A$] and
\emph{transitive} [$a \ \rho \ b$ and $b \ \rho \ c \ \Rightarrow \
a \ \rho \ c\ \forall a,b,c \in A $], it is a \emph{partial ordering
relation} and $A$ a \emph{partially ordered set} or \emph{poset}. We
denote a poset $A$ with partial ordering relation $\rho$ by $\langle
A; \rho\rangle$, abbreviated by simply $A$ if the binary relation is
clear.

For $H \subseteq A$ and $a \in A$, $a$ is an \emph{upper bound} of
$H$ if $h \leq a$ for all $h \in H$. $a$ is the \emph{least upper
bound} or \emph{supremum} of $H$ if every upper bound $b$ of $H$
satisfies $a \leq b$ and we then write $a = \sup H$. \emph{Lower
bound} and \emph{greatest lower bound} or \emph{infimum}, denoted
$\inf H$, are defined similarly. $\sup \emptyset$ is the smallest
element in $A$, called $zero$, if it exists, and $\inf \emptyset$ is
the largest element in $A$, called $unit$, if it exists.

A poset $\langle L; \leq \rangle$ is a \emph{lattice} if $\inf H$
and $\sup H$ exist for any finite nonempty subset $H$ of $L$. It is
called \emph{complete} if $\inf H$ and $\sup H$ also exist for $H =
\emptyset$. A poset can be shown to be a complete lattice with the
following result.

\begin{lemma} \label{lemma_complete_lattice}
If $\langle P; \leq \rangle$ is a poset in which $\inf H$ exists for
all $H \subseteq P$, then $\langle P; \leq \rangle$ is a complete
lattice.
\end{lemma}

For a lattice $L$ and $a,b \in L$, we write $a \wedge b$ for
$\inf\{a,b\}$ and $a \vee b$ for $\sup \{a,b\}$, and refer to
$\wedge$ as the \emph{meet operation} and to $\vee$ as the
\emph{join operation}. $L$ is \emph{distributive} if for all $a,b,c
\in L$,
\begin{equation} \label{eq_distr}
a \vee (b \wedge c) = (a \vee b) \wedge (a \vee c)
\end{equation}

The structure of a lattice $\langle L; \leq \rangle$ may be
visualised by a \emph{Hasse diagram}, in which each element pair $a,
b \in L$ is joined by an edge whenever $a \leq b$ and there is no $x
\in L \setminus \{a, b\}$ such that $a \leq x \leq b$.

We denote most partial orderings by $\leq$. Which partial ordering
the symbol refers to will be determined by the context. For an
overview of lattice theory see \citet{graetzer}.

\section{Model Types: RCON and RCOR Models}

\subsection{RCON Models: Equality Restrictions on Concentrations}

RCON models are graphical Gaussian models which place equality
constraints on the entries of the concentration matrix $K =
\Sigma^{-1}$. For a model whose conditional independence structure
is represented by graph $G=(V,E)$, the restrictions can be
represented by a graph colouring $(\mathcal{V},\mathcal{E})$, with
the vertex colouring $\mathcal{V}$ representing constraints on the
entries on the diagonal of $K$ and the edge colouring $\mathcal{E}$
representing constraints in the off-diagonal entries. Whenever two
vertices $\alpha, \beta \in V$ belong to the same vertex colour
class, the corresponding two diagonal entries $k_{\alpha\alpha}$ and
$k_{\beta\beta}$ are restricted to being identical. Similarly, two
edges $\alpha\beta, \gamma\delta \in E$ of the same colour represent
the constraint $k_{\alpha\beta} = k_{\gamma\delta}$.

We denote the set of positive definite matrices which satisfy such
constraints for a graph colouring $(\mathcal{V},\mathcal{E})$ by
$\mathcal{S}^{+}(\mathcal{V},\mathcal{E})$. Put formally, the
distribution of a random vector $Y \in \mathbb{R}^{V}$ is said to
lie in the RCON model represented by the coloured graph
$\mathcal{G}=(\mathcal{V},\mathcal{E})$ if
\begin{equation} \label{rcon_par}
Y \sim \mathcal{N}_{V}(0,\Sigma), \quad K = \Sigma^{-1} \in
\mathcal{S}^{+}(\mathcal{V},\mathcal{E}) = \left\{ \sum_{u \in
\mathcal{V} \cup \mathcal{E}} \lambda_{u}T^{u}, \ \lambda \in
\mathbb{R}^{\mathcal{V} \cup \mathcal{E}} \right\}^{+}
\end{equation}

Since the constraints are linear in $K$, by standard exponential
family theory \citep{brown}, just as unconstrained graphical Gaussian
models, RCON models are regular exponential families. Thus the
maximum likelihood estimate of $\lambda$ is uniquely determined,
provided it exists. For the corresponding computation algorithm
see \citet{lauritzen_sym} and \citet{grc}. Note that RCON models are
instances of models considered in \citet{anderson70}.

\begin{example} \label{rcon_example}
The data consist of the examination marks of 88 students in the
mathematical subjects Algebra, Analysis, Mechanics, Statistics and
Vectors \citep{mardia}. \citet{whittaker, edwards} previously
demonstrated an excellent fit to the unconstrained model represented
by the graph shown in Figure
\ref{fig:mathmark_uncol}. \citet{lauritzen_sym} show the data to also
support the RCON model represented by the graph in Figure
\ref{fig:mathmark_col}. The model specifications are
\begin{equation*}
Y \sim \mathcal{N}_{5}(0,\Sigma), \ \ \ \ \ \Sigma^{-1} \in M
\end{equation*}
with $M$ given below the corresponding graphs. If the subjects are
indexed in alphabetical order, the graph {\spaceskip=1.6pt colouring
representing the constraints of the RCON model is given by
$(\mathcal{V},\mathcal{E})$ with $
\mathcal{V}=\{\{1\},\{2,5\},\{3,4\}\}$} and
$\mathcal{E}=\{\{12\},\{13,14,15,24,35\}\}$. Note that the number of
model parameters has been reduced from 11 to 5.

\newpage
\begin{figure}[ht!]
  \centering
  \subfigure[Inferred conditional independence structure of Mathematics marks data.]{\label{fig:mathmark_uncol}
\begin{picture}(180,190) (10,-10)
        \begin{normalsize}
        \thicklines
        \put(40,100){\line(0,1){60}}
        \put(160,100){\line(0,1){60}}
        \put(40,100){\line(2,1){60}}
        \put(40,160){\line(2,-1){60}}
        \put(160,100){\line(-2,1){60}}
        \put(160,160){\line(-2,-1){60}}

        \put(40,100){\circle*{6}}
        \put(40,160){\circle*{6}}
        \put(160,100){\circle*{6}}
        \put(100,130){\circle*{6}}
        \put(160,160){\circle*{6}}

        \put(20,80){\footnotesize Mechanics}
        \put(140,80){\footnotesize Statistics}
        \put(20,170){\footnotesize Vectors}
        \put(140,170){\footnotesize Analysis}
        \put(85,143){\footnotesize Algebra}

        \put(25,30){\footnotesize $M= \left\{\left(
        \begin{array}{ccccc}
        * & * & * & * & *\\
          & * & 0 & * & 0\\
          &   & * & 0 & * \\
          &   &   & * & 0 \\
          &   &   &   & *
        \end{array}\right)\right\}^{+}$}
        \end{normalsize}
\end{picture}
  }
\hspace{20pt}
  \subfigure[RCON model supported by Mathematics marks data.]{\label{fig:mathmark_col}
\begin{picture}(180,190) (5,-10)
\begin{normalsize}
\thicklines
        \put(40,100){\color{darkgreen}{\line(0,1){60}}}
        \put(160,100){\color{darkgreen}{\line(0,1){60}}}
        \put(40,100){\color{darkgreen}{\line(2,1){60}}}
        \put(40,160){\color{darkgreen}{\line(2,-1){60}}}
        \put(160,100){\color{darkgreen}{\line(-2,1){60}}}
        \put(160,160){\line(-2,-1){60}}

        \put(40,100){\color{blue}{\circle*{6}}}
        \put(40,160){\color{red}{\circle*{6}}}
        \put(160,100){\color{blue}{\circle*{6}}}
        \put(100,130){\circle*{6}}
        \put(160,160){\color{red}{\circle*{6}}}

                \put(37,130){\line(1,0){6}}
                \put(157,130){\line(1,0){6}}
                \put(70,112){\line(0,1){6}}
                \put(70,142){\line(0,1){6}}
                \put(130,112){\line(0,1){6}}

                \put(28,98){$*$}
                \put(167,98){$*$}
                \put(22,158){$**$}
                \put(167,158){$**$}

        \put(20,80){\footnotesize Mechanics}
        \put(140,80){\footnotesize Statistics}
        \put(20,170){\footnotesize Vectors}
        \put(140,170){\footnotesize Analysis}
        \put(85,143){\footnotesize Algebra}

                \put(15,30){\footnotesize $M= \left\{\left(
        \begin{array}{ccccc}
        * & * & \color{darkgreen}{\lambda_{13}} & \color{darkgreen}{\lambda_{13}} & \color{darkgreen}{\lambda_{13}}\\
          & \color{red}{\lambda_{22}} & 0 & \color{darkgreen}{\lambda_{13}} & 0\\
          &   & \color{blue}{\lambda_{33}} & 0 & \color{darkgreen}{\lambda_{13}} \\
          &   &   & \color{blue}{\lambda_{33}} & 0 \\
          &   &   &   & \color{red}{\lambda_{22}}
        \end{array}\right)\right\}^{+}$}

        \end{normalsize}
\end{picture}
  }
  \caption{Mathematics marks example.}
  \label{fig:mathmark_ex}
\end{figure}
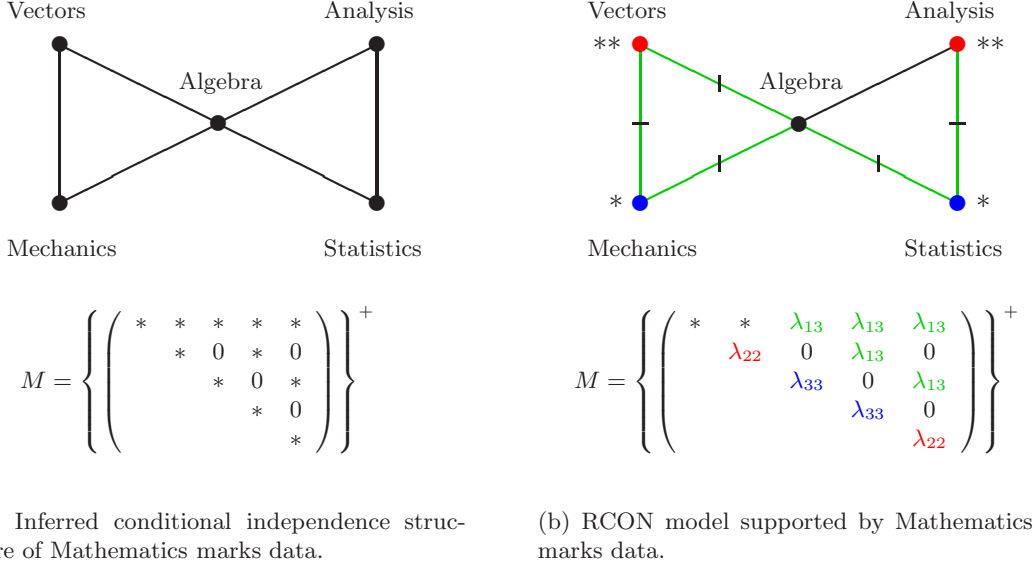

\end{example}

\subsection{RCOR Models: Equality Restrictions on Partial Correlations}

RCOR models place symmetry restrictions on the diagonal elements of
the concentration matrix $K=\Sigma^{-1}$ and on the
partial correlations as given in equation (\ref{partial_cor}). Just
as for RCON models, for a model with graph $G=(V,E)$, the
constraints can be represented by a graph colouring
$(\mathcal{V},\mathcal{E})$: Vertices of the same colour represent
restrictions on the diagonal entries of $K$ (exactly as in RCON
models), and whenever two edges $\alpha\beta, \gamma\delta \in E$
belong to the same edge colour class in $\mathcal{E}$, the
corresponding partial correlations $\rho_{\alpha\beta | V \setminus
\{\alpha,\beta\}}$ and $\rho_{\gamma\delta | V \setminus
\{\gamma,\delta\}}$, defined in equation (\ref{partial_cor}), are
restricted to being identical.

We denote the set of positive definite matrices which satisfy such
restrictions for a graph colouring $(\mathcal{V},\mathcal{E})$ by
$\mathcal{R}^{+}(\mathcal{V},\mathcal{E})$. If we let the $|V|
\times |V|$ matrices $A=(a_{\alpha\beta})$ and $C=(c_{\alpha\beta})$
be given by $a_{\alpha\beta} = \sqrt{k_{\alpha\beta}}$ for $\alpha =
\beta$ and zero otherwise, and let $c_{\alpha\beta}$ as in
equation (\ref{eq:c}) for $\alpha\not=\beta$ and $c_{\alpha\beta} =
1$ otherwise, then $K=ACA$ and the distribution of a random vector
$Y \in \mathbb{R}^{V}$ lies in the RCOR model represented by the
coloured graph $\mathcal{G}=(\mathcal{V},\mathcal{E})$ if
\begin{eqnarray} \label{rcor_par}
&&Y \sim \mathcal{N}_{V}(0,\Sigma), \quad K = \Sigma^{-1} \in \mathcal{R}^{+}(\mathcal{V},\mathcal{E}) = \left\{ ACA \mid A = \sum_{u \in \mathcal{V}} \eta_{u} T^{u}, \ \eta \in \mathbb{R}^{\mathcal{V}}_{+}, \right. \nonumber\\
&& \left. \hspace{250pt} C = I + \sum_{u \in \mathcal{E}} \tau_{u} T^{u},\ \tau \in (-1,1)^{\mathcal{E}}\right\}^{+}
\end{eqnarray}

Thus the constraints of RCOR models define a differentiable manifold
in $\mathcal{S}^{+}$ which makes them curved exponential
families \citep{brown}. Therefore, the maximum likelihood estimates
of $\eta$ and $\tau$, if they exist, may not be unique. For a
discussion and computation algorithm we refer
to \citet{lauritzen_sym} and \citet{grc}.

RCOR models are scale invariant if variables inside the same vertex
colour class are manipulated in the same way \citep{lauritzen_sym}.
Thus they are particularly suitable for variables measured on
different scales.

We highlight that both RCON and RCOR models generally do not place
the same equality restrictions on $\Sigma$ as they do on
$\Sigma^{-1}$ and on partial correlations.

\begin{example} \label{rcor_example}
The data is concerned with anxiety and anger in a trait and state
version of 684 students \citep{wermuth} and strongly support the
conditional independence model displayed in Figure
\ref{fig:anxiety_uncol}. As shown in \citet{lauritzen_sym}, they also
support the RCOR model represented by the coloured graph in Figure
\ref{fig:anxiety_col}, parametrised by 6 parameters rather than 8.
The variable names are combinations between $T$ or $S$, for ``trait''
and ``state'', and $X$ or $N$, standing for ``anxiety'' and ``anger''. The
model specifications are
\begin{equation*}
Y \sim \mathcal{N}_{4}(0,\Sigma), \ \ \ \ \ \Sigma^{-1} \in M
\end{equation*}
with $M$ given below the graphs. The variables are indexed
anti-clockwise starting from $TX$.

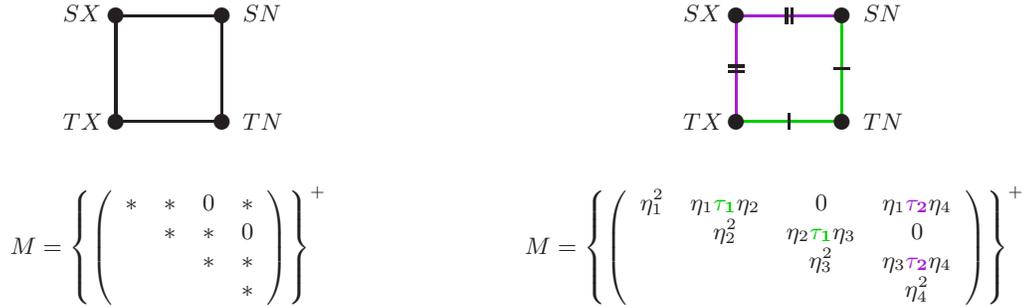
\begin{figure}[ht!]
  \centering
  \subfigure[Conditional independence structure supported by personality characteristics data.]{\label{fig:anxiety_uncol}
\begin{picture}(180,135) (20,-10)
        \begin{normalsize}
        \thicklines
        \put(60,67){\footnotesize $TX$}
        \put(60,108){\footnotesize $SX$}
        \put(128,67){\footnotesize $TN$}
        \put(128,108){\footnotesize $SN$}
        \put(80,70){\line(1,0){40}}
        \put(80,70){\line(0,1){40}}
        \put(120,110){\line(-1,0){40}}
        \put(120,110){\line(0,-1){40}}
        \put(80,70){\circle*{6}}
        \put(80,110){\circle*{6}}
        \put(120,70){\circle*{6}}
        \put(120,110){\circle*{6}}

        \put(40,20){\footnotesize $M= \left\{\left(
        \begin{array}{cccc}
        * & * & 0 & * \\
          & * & * & 0 \\
          &   & * & * \\
          &   &   & *
        \end{array}\right)\right\}^{+}$}

        \end{normalsize}
\end{picture}
  }
\hspace{30pt}
  \subfigure[RCOR model supported by personality characteristics data.]{\label{fig:anxiety_col}
\begin{picture}(180,115) (20,-10)
        \begin{normalsize}
        \thicklines

        \put(70,67){\footnotesize $TX$}
        \put(70,108){\footnotesize $SX$}
        \put(138,67){\footnotesize $TN$}
        \put(138,108){\footnotesize $SN$}
        \put(90,70){\color{darkgreen}{\line(1,0){40}}}
        \put(90,70){\color{purple}{\line(0,1){40}}}
        \put(130,110){\color{purple}{\line(-1,0){40}}}
        \put(130,110){\color{darkgreen}{\line(0,-1){40}}}
        \put(90,70){\circle*{6}}
        \put(90,110){\circle*{6}}
        \put(130,70){\circle*{6}}
        \put(130,110){\circle*{6}}

        \put(127,90){\line(1,0){6}}
        \put(87,89){\line(1,0){6}}
        \put(87,91){\line(1,0){6}}
        \put(110,67){\line(0,1){6}}
        \put(109,107){\line(0,1){6}}
        \put(111,107){\line(0,1){6}}

        \put(10,20){\footnotesize $M= \left\{\left(
        \begin{array}{cccc}
        \eta_{1}^{2} & \eta_{1}\mathbf{\color{darkgreen}{\tau_{1}}}\eta_{2} & 0 & \eta_{1}\mathbf{\color{purple}{\tau_{2}}}\eta_{4} \\
          & \eta_{2}^{2} & \eta_{2}\mathbf{\color{darkgreen}{\tau_{1}}}\eta_{3} & 0 \\
          &   & \eta_{3}^{2} & \eta_{3}\mathbf{\color{purple}{\tau_{2}}}\eta_{4} \\
          &   &   & \eta_{4}^{2}
        \end{array}\right)\right\}^{+}$}

        \end{normalsize}
\end{picture}
  }
  \vspace{-5pt}
  \caption{Personality characteristics example.}
\end{figure}
\end{example}

\subsection{Number of RCON and RCOR Models}

Let $\mathcal{S}^{+}_{V}$ and $\mathcal{R}^{+}_{V}$ denote the sets
of RCON and RCOR models with variable set $V$ and let
$\mathcal{C}_{V}$ be the set of vertex and edge coloured graphs with
vertex set $V$. Further, let $M_{V}$ be the set of unconstrained
graphical Gaussian models with variable set $V$ and $U_{V}$ the set
of undirected graphs with vertex set $V$.

As, by equations (\ref{rcon_par}) and (\ref{rcor_par}), for both
model types there is one model parameter for each vertex colour
class in $\mathcal{V}$ and one for each edge colour class in
$\mathcal{E}$ in the coloured dependence graph $\mathcal{G} =
(\mathcal{V},\mathcal{E})$, there are as many RCON and RCOR models
with variables $V$ as there are coloured graphs with with vertex set
$V$, \textit{i.e.\/}, $|\mathcal{S}^{+}_{V}| = |\mathcal{R}^{+}_{V}|
= |\mathcal{C}_{V}|$. Given that the number of graph colourings of a
particular graph $G=(V,E)$ is given by the product $|P(V)||P(E)|$ of
the number of partitions of $V$ multiplied by the number of
partitions of $E$, we obtain
\begin{equation} \label{eq_cv}
|\mathcal{C}_{V}| = \sum_{G=(V,E) \in U_{V}} |P(V)||P(E)| = |P(V)|
\sum_{G=(V,E) \in U_{V}}|P(E)|
\end{equation}

For a discrete set $S$ of size $d$, $|P(S)|$ is given the $d^{th}$
\emph{Bell number} $B_{d}$ \citep{bell, pitman}, which satisfy the
recursive relationship $B_{d+1} = \sum_{k=0}^{d}{d \choose k}B_{k}$,
with $B_{0}=1$. Hence,
\begin{equation*}
|\mathcal{S}^{+}_{V}| = |\mathcal{R}^{+}_{V}| = B_{|V|}\sum_{G=(V,E)
\in U_{V}} B_{|E|} = B_{|V|} \sum_{k=0}^{{|V| \choose 2}} {{|V|
\choose 2} \choose k} B_{k} = B_{|V|} B_{{|V| \choose 2}+1}
\end{equation*}
For each $d$, $B_{d}$ can be evaluated as the least integer greater
than the sum of the first $2d$ terms in \emph{Dobi\'{n}ski's
formula} \citep{dobinski,comlet}
\begin{equation} \label{bell}
B_{d} = e^{-1}\sum_{k=0}^{\infty} \frac{k^{d}}{k!} = e^{-1}
\left(\frac{0^{d}}{0!} + \frac{1^{d}}{1!} + \frac{2^{d}}{2!} +
\ldots \right)
\end{equation}
so that clearly $|\mathcal{S}^{+}_{V}| = |\mathcal{R}^{+}_{V}|$ grow
super-exponentially in $|V|$. For illustration, observe that while
$|M_{[4]}| =64$ and $|M_{[5]}|= 1,024$, $|\mathcal{S}^{+}_{[4]}| =
|\mathcal{R}^{+}_{[4]}| = 13,155$ and $|\mathcal{S}^{+}_{[5]}| =
|\mathcal{R}^{+}_{[5]}| = 35,285,640$.

\subsection{Structure of the Sets of RCON and RCOR Models}

It is a well-known fact that $M_{V}$ is a complete distributive
lattice with respect to model inclusion, with partial ordering
induced by the partial ordering on $U_{V}$ given by edge set
inclusion: for $G_{1}=(V,E_{1}), G_{2}=(V,E_{2}) \in
U_{V}$, $G_{1} \leq G_{2}$ whenever $E_{1} \subseteq E_{2}$, with
$G_{1} \land G_{2} = (V, E_{1} \cap E_{2})$ and $G_{1}
\lor G_{2} = (V, E_{1} \cup E_{2})$. For $M_{1}, M_{2} \in M_{V}$
represented by $G_{1}, G_{2}$ as above, $M_{1} \subseteq M_{2}$ if
and only if $G_{1} \leq G_{2}$. The zero in $\langle U_{V}; \leq
\rangle$ is the empty graph, and the unit the complete graph, in
which every edge is present.

RCON and RCOR models are specified by partitions of $V$ and $E$. For
any finite discrete set $S$, the set $P(S)$ of partitions of $S$
forms a complete non-distributive lattice, with $P_{1} \leq P_{2}$
for $P_{1}, P_{2} \in P(S)$ whenever $P_{1}$ is \emph{finer} than
$P_{2}$, or, put differently, whenever $P_{2}$ is \emph{coarser}
than $P_{1}$, \textit{i.e.\/}, if every set in $P_{2}$ can be
expressed as a union of sets in $P_{1}$. This allows the
identification of a partial ordering $\preceq$ on $\mathcal{C}_{V}$
which corresponds to model inclusion in $\mathcal{S}^{+}_{V}$ and
$\mathcal{R}^{+}_{V}$: For
$\mathcal{G}=(\mathcal{V}_{\mathcal{G}},\mathcal{E}_{\mathcal{G}}),
\mathcal{H}=(\mathcal{V}_{\mathcal{H}},\mathcal{E}_{\mathcal{H}})
\in \mathcal{C}_{V}$ with underlying uncoloured graphs $G$ and $H$,
$\mathcal{G} \preceq \mathcal{H}$ whenever\\

   \vspace{-5pt}
\noindent (i) $\ G \leq H$;  \quad (ii) $\ \mathcal{V}_{\mathcal{G}}
\geq \mathcal{V}_{\mathcal{H}}$;  \quad (iii)$\ $ every colour
class in $\mathcal{E}_{\mathcal{G}}$ is a union of colour classes in
$\mathcal{E}_{\mathcal{H}}$. \\ \vspace{-5pt}

Put in words, if we let $\mathcal{M}_{\mathcal{G}},
\mathcal{M}_{\mathcal{H}}$ denote two RCON or RCOR models (both of
the same type) represented by $\mathcal{G}, \mathcal{H} \in
\mathcal{C}_{V}$, then $\mathcal{M}_{\mathcal{G}} \subseteq
\mathcal{M}_{\mathcal{H}}$ if $\mathcal{H}$ can be obtained from
$\mathcal{G}$ by splitting of colour classes and adding new edge
colour classes, or equivalently if $\mathcal{G}$ can be obtained
from $\mathcal{H}$ by merging colour classes and dropping edge
colour classes.

\begin{figure}[ht!]
  \centering
\begin{picture}(390,45)(0,20)
\begin{normalsize}
\thicklines
        \put(6,10){1}
        \put(6,50){4}
        \put(66,10){2}
        \put(66,50){3}
        \put(20,10){\color{darkgreen}{\line(1,0){40}}}
        \put(60,50){\color{darkgreen}{\line(0,-1){40}}}
        \put(20,10){\color{blue}{\circle*{6}}}
        \put(60,10){\color{blue}{\circle*{6}}}
        \put(60,50){\circle*{6}}
        \put(20,50){\circle*{6}}

        \put(17,-2){$*$}
        \put(57,-2){$*$}

        \put(40,7){\line(0,1){6}}
        \put(57,30){\line(1,0){6}}

                \put(35,60){\footnotesize $\mathcal{G}_{1}$}

                \put(85,30){$\preceq$}

                \put(106,10){1}
        \put(106,50){4}
        \put(166,10){2}
        \put(166,50){3}
        \put(120,10){\line(1,0){40}}
        \put(120,10){\color{purple}{\line(0,1){40}}}
        \put(120,10){\color{purple}{\line(1,1){40}}}
        \put(160,50){\color{purple}{\line(-1,0){40}}}
        \put(160,50){\line(0,-1){40}}
        \put(120,10){\color{blue}{\circle*{6}}}
        \put(120,50){\circle*{6}}
        \put(160,10){\color{blue}{\circle*{6}}}
        \put(160,50){\circle*{6}}

        \put(157,-2){$*$}
        \put(117,-2){$*$}

                \put(139,47){\line(0,1){6}}
                \put(141,47){\line(0,1){6}}
                \put(117,29){\line(1,0){6}}
                \put(117,31){\line(1,0){6}}
                \put(139,27){\line(0,1){6}}
                \put(141,27){\line(0,1){6}}

                \put(135,60){\footnotesize $\mathcal{G}_{2}$}

        \put(226,10){1}
        \put(226,50){4}
        \put(286,10){2}
        \put(286,50){3}
        \put(240,10){\color{darkgreen}{\line(1,0){40}}}
        \put(280,50){\color{darkgreen}{\line(0,-1){40}}}
        \put(240,10){\color{blue}{\circle*{6}}}
        \put(280,10){\color{blue}{\circle*{6}}}
        \put(280,50){\circle*{6}}
        \put(240,50){\circle*{6}}

        \put(237,-2){$*$}
        \put(277,-2){$*$}

        \put(260,7){\line(0,1){6}}
        \put(277,30){\line(1,0){6}}

                \put(255,60){\footnotesize $\mathcal{G}_{1}$}

                \put(305,30){$\not\preceq$}

                \put(326,10){1}
        \put(326,50){4}
        \put(386,10){2}
        \put(386,50){3}
        \put(340,10){\color{darkgreen}{\line(1,0){40}}}
        \put(340,10){\color{purple}{\line(0,1){40}}}
        \put(340,10){\color{darkgreen}{\line(1,1){40}}}
        \put(380,50){\color{purple}{\line(-1,0){40}}}
        \put(380,50){\color{darkgreen}{\line(0,-1){40}}}
        \put(340,10){\color{blue}{\circle*{6}}}
        \put(340,50){\color{red}{\circle*{6}}}
        \put(380,10){\color{blue}{\circle*{6}}}
        \put(380,50){\color{red}{\circle*{6}}}

        \put(377,-2){$*$}
        \put(337,-2){$*$}
                \put(375,57){$**$}
        \put(335,57){$**$}

                \put(359,47){\line(0,1){6}}
                \put(361,47){\line(0,1){6}}
                \put(337,29){\line(1,0){6}}
                \put(337,31){\line(1,0){6}}
                \put(377,30){\line(1,0){6}}
                \put(360,27){\line(0,1){6}}
                \put(360,7){\line(0,1){6}}

                \put(355,60){\footnotesize $\mathcal{G}_{3}$}
        \end{normalsize}
\end{picture}\\
\vspace{20pt}
  \caption{Partial ordering in $\mathcal{C}_{[4]}$.}
  \label{fig:part_ord_ex}
\end{figure}
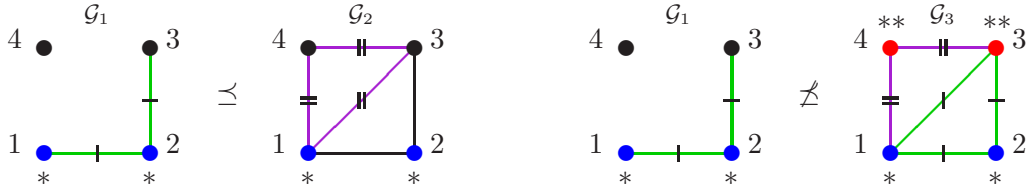

For example, for the graphs
$\mathcal{G}_{i}=(\mathcal{V}_{i},\mathcal{E}_{i})$ for $i = 1,2,3$
in Figure \ref{fig:part_ord_ex}, $\mathcal{G}_{1} \preceq
\mathcal{G}_{2}$ as conditions (i)--(iii) above are satisfied
whereas $\mathcal{G}_{1} \not\preceq \mathcal{G}_{3}$ because (ii)
and (iii) are violated. Thus the corresponding RCON or RCOR models
$\mathcal{M}_{1},\mathcal{M}_{2},\mathcal{M}_{3}$ (all of the same
type) satisfy $\mathcal{M}_{1} \subseteq \mathcal{M}_{2}$ and
$\mathcal{M}_{1} \not\subseteq \mathcal{M}_{3}$.

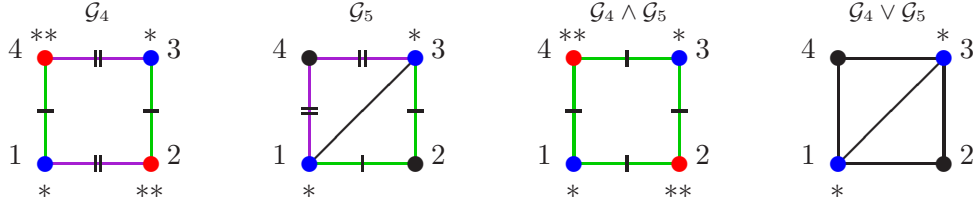
\begin{figure}[ht!]
  \centering
\begin{picture}(360,90)(0,-20)
\begin{normalsize}
\thicklines
        \put(6,0){1}
        \put(6,40){4}
        \put(66,0){2}
        \put(66,40){3}
        \put(20,0){\color{purple}{\line(1,0){40}}}
        \put(20,0){\color{darkgreen}{\line(0,1){40}}}
        \put(60,40){\color{purple}{\line(-1,0){40}}}
        \put(60,40){\color{darkgreen}{\line(0,-1){40}}}
        \put(20,0){\color{blue}{\circle*{6}}}
        \put(20,40){\color{red}{\circle*{6}}}
        \put(60,0){\color{red}{\circle*{6}}}
        \put(60,40){\color{blue}{\circle*{6}}}
        \put(17,-13){$*$}
        \put(54,-13){$**$}
        \put(14,46){$**$}
        \put(57,46){$*$}
                \put(39,-3){\line(0,1){6}}
                \put(41,-3){\line(0,1){6}}
                \put(39,37){\line(0,1){6}}
                \put(41,37){\line(0,1){6}}
                \put(17,20){\line(1,0){6}}
                \put(57,20){\line(1,0){6}}

                \put(35,55){\footnotesize $\mathcal{G}_{4}$}

                \put(106,0){1}
        \put(106,40){4}
        \put(166,0){2}
        \put(166,40){3}
        \put(120,0){\color{darkgreen}{\line(1,0){40}}}
        \put(120,0){\color{purple}{\line(0,1){40}}}
        \put(120,0){\line(1,1){40}}
        \put(160,40){\color{purple}{\line(-1,0){40}}}
        \put(160,40){\color{darkgreen}{\line(0,-1){40}}}
        \put(120,0){\color{blue}{\circle*{6}}}
        \put(120,40){\circle*{6}}
        \put(160,0){\circle*{6}}
        \put(160,40){\color{blue}{\circle*{6}}}
        \put(117,-13){$*$}
        \put(157,46){$*$}
                \put(140,-3){\line(0,1){6}}
                \put(139,37){\line(0,1){6}}
                \put(141,37){\line(0,1){6}}
                \put(117,19){\line(1,0){6}}
                \put(117,21){\line(1,0){6}}
                \put(157,20){\line(1,0){6}}

                \put(135,55){\footnotesize $\mathcal{G}_{5}$}

                \put(206,0){1}
        \put(206,40){4}
        \put(266,0){2}
        \put(266,40){3}
        \put(220,0){\color{darkgreen}{\line(1,0){40}}}
        \put(220,0){\color{darkgreen}{\line(0,1){40}}}
        \put(260,40){\color{darkgreen}{\line(-1,0){40}}}
        \put(260,40){\color{darkgreen}{\line(0,-1){40}}}
        \put(220,0){\color{blue}{\circle*{6}}}
        \put(220,40){\color{red}{\circle*{6}}}
        \put(260,0){\color{red}{\circle*{6}}}
        \put(260,40){\color{blue}{\circle*{6}}}
        \put(217,-13){$*$}
        \put(254,-13){$**$}
        \put(214,46){$**$}
        \put(257,46){$*$}
                \put(240,-3){\line(0,1){6}}
                \put(240,37){\line(0,1){6}}
                \put(217,20){\line(1,0){6}}
                \put(257,20){\line(1,0){6}}

                \put(226,55){\footnotesize $\mathcal{G}_{4} \wedge \mathcal{G}_{5}$}

                \put(306,0){1}
        \put(306,40){4}
        \put(366,0){2}
        \put(366,40){3}
        \put(320,0){\line(1,0){40}}
        \put(320,0){\line(0,1){40}}
        \put(320,0){\line(1,1){40}}
        \put(360,40){\line(-1,0){40}}
        \put(360,40){\line(0,-1){40}}
        \put(320,0){\color{blue}{\circle*{6}}}
        \put(320,40){\circle*{6}}
        \put(360,0){\circle*{6}}
        \put(360,40){\color{blue}{\circle*{6}}}
        \put(317,-13){$*$}
        \put(357,46){$*$}

                \put(324,55){\footnotesize $\mathcal{G}_{4} \vee \mathcal{G}_{5}$}

        \end{normalsize}
\end{picture}
\vspace{-10pt}
  \caption{Meet and join operations in $\mathcal{C}_{[4]}$.}
  \label{fig:join_col_ex}

\end{figure}

It is then straight forward to show that $\langle \mathcal{C}_{V};
\preceq \rangle$ is a complete lattice with meet and join operations
\begin{equation} \label{meet_join_c}
\mathcal{G} \wedge \mathcal{H} = (\mathcal{V}_{\mathcal{G}} \vee
\mathcal{V}_{\mathcal{H}}, \mathcal{E}^{*}_{\mathcal{G}} \vee
\mathcal{E}^{*}_{\mathcal{H}}) \ \ \ \ \ \  \textup{and} \ \ \ \ \ \
\mathcal{G} \vee \mathcal{H} = (\mathcal{V}_{\mathcal{G}} \wedge
\mathcal{V}_{\mathcal{H}}, \mathcal{E}^{**}_{\mathcal{G}} \wedge
\mathcal{E}^{**}_{\mathcal{H}})
\end{equation}
where $\mathcal{E}^{*}_{\mathcal{G}} \subseteq
\mathcal{E}_{\mathcal{G}}$ and $\mathcal{E}^{*}_{\mathcal{H}}
\subseteq \mathcal{E}_{\mathcal{H}}$ are maximal with the property
that they are partitions of the same set of edges inside
$E_{\mathcal{G}} \cap E_{\mathcal{H}}$,
$\mathcal{E}^{**}_{\mathcal{G}} = \mathcal{E}_{\mathcal{G}} \cup
\{\{E_{\mathcal{H}} \setminus E_{\mathcal{G}}\}\}$ and
$\mathcal{E}^{**}_{\mathcal{H}} = \mathcal{E}_{\mathcal{H}} \cup
\{\{E_{\mathcal{G}} \setminus E_{\mathcal{H}}\}\}$. The graphs in
Figure \ref{fig:join_col_ex} illustrate the operations. The zero in
$\langle \mathcal{C}_{V}; \preceq \rangle$ is given by the empty
graph in which all vertices are of the same colour and the unit is
the complete graph with atomic colour classes.


\begin{proposition} \label{prop_rcon_rcor_c}
Let
$\mathcal{G}=(\mathcal{V}_{\mathcal{G}},\mathcal{E}_{\mathcal{G}}),
\mathcal{H}=(\mathcal{V}_{\mathcal{H}},\mathcal{E}_{\mathcal{H}})
\in \mathcal{C}_{V}$ and let
$\mathcal{S}^{+}(\mathcal{V}_{\mathcal{G}},\mathcal{E}_{\mathcal{G}}),
\mathcal{S}^{+}(\mathcal{V}_{\mathcal{H}},\mathcal{E}_{\mathcal{H}})
\in \mathcal{S}^{+}_{V}$ and
$\mathcal{R}^{+}(\mathcal{V}_{\mathcal{G}},\mathcal{E}_{\mathcal{G}}),
\mathcal{R}^{+}(\mathcal{V}_{\mathcal{H}},\mathcal{E}_{\mathcal{H}})
\in \mathcal{R}^{+}_{V}$ be the RCON and RCOR models represented by
$\mathcal{G}$ and $\mathcal{H}$. Then
\begin{equation*}
\mathcal{S}^{+}(\mathcal{V}_{\mathcal{G}},\mathcal{E}_{\mathcal{G}})
\subseteq
\mathcal{S}^{+}(\mathcal{V}_{\mathcal{H}},\mathcal{E}_{\mathcal{H}})
\ \ \ \Longleftrightarrow \ \ \ \mathcal{G} \preceq \mathcal{H} \ \
\ \Longleftrightarrow \ \ \
\mathcal{R}^{+}(\mathcal{V}_{\mathcal{G}},\mathcal{E}_{\mathcal{G}})
\subseteq
\mathcal{R}^{+}(\mathcal{V}_{\mathcal{H}},\mathcal{E}_{\mathcal{H}})
\end{equation*}
and $\mathcal{S}^{+}_{V}$ and $\mathcal{R}^{+}_{V}$ are complete
non-distributive lattices with meet and join operations induced by
the meet and join operations in $\langle \mathcal{C}_{V} \preceq
\rangle$, given in equation (\ref{meet_join_c}).
\end{proposition}

\section{Model Classes within RCON and RCOR Models} \label{sym_model_classes}

The motivation to study model classes strictly
within the sets of RCON and RCOR models is three-fold: firstly,
having demonstrated that the number of RCON and RCOR models grows
dramatically with the number of variables, especially for model
selection, smaller model (search) spaces are desirable. Secondly,
generic equality constraints of RCON and RCOR models are generally
not readily interpretable and, lastly, do not guarantee the
corresponding model to have any particular statistical properties.

Four model classes within the sets of RCON and RCOR models which are
characterised by desirable statistical properties expressing
themselves in regularity of colouring have been previously
identified in the literature. This section is devoted to their
definition and first properties. Three of the four colouring
regularities were termed \emph{edge regularity},
with the corresponding models appearing
in \citet{lauritzen_sym}, \emph{vertex regularity} and
\emph{regularity} in \citet{gehrmann}. We term colourings of the
fourth type \emph{permutation-generated}. The corresponding models
are referred to as \emph{graphical symmetry models}
in \citet{hylleberg} and as \emph{RCOP
models} in \citet{lauritzen_sym}.

\subsection{Models Represented by Edge Regular Colourings}

RCON and RCOR models place restrictions on different parameter sets,
which translates into different model properties. While
the restrictions in RCON models ensure the models to be regular
exponential families, RCOR models are scale invariant within vertex
colour classes. Thus if a graph colouring
$(\mathcal{V},\mathcal{E})$ yields the same model restrictions
representing the constraints of an RCON model as it does
representing those of an RCOR model, it represents a model with both
of the above desirable properties. Such models can be identified by
their graph colouring.

\begin{definition} \label{def_edge_reg}
For $\mathcal{G} = (\mathcal{V},\mathcal{E}) \in \mathcal{C}_{V}$ we
say that $(\mathcal{V},\mathcal{E})$ is \emph{edge regular} if any
pair of edges in the same edge colour class in $\mathcal{E}$
connects the same vertex colour classes in $\mathcal{V}$.
\end{definition}

It then holds:

\begin{proposition}[\cite{lauritzen_sym}] \label{prop_edge_reg}
The RCON and RCOR models determined by $(\mathcal{V},\mathcal{E})$
yield identical restrictions
\begin{equation*}
\mathcal{S}^{+}(\mathcal{V},\mathcal{E}) =
\mathcal{R}^{+}(\mathcal{V},\mathcal{E})
\end{equation*}
\noindent if and only if $(\mathcal{V},\mathcal{E})$ is edge
regular.
\end{proposition}

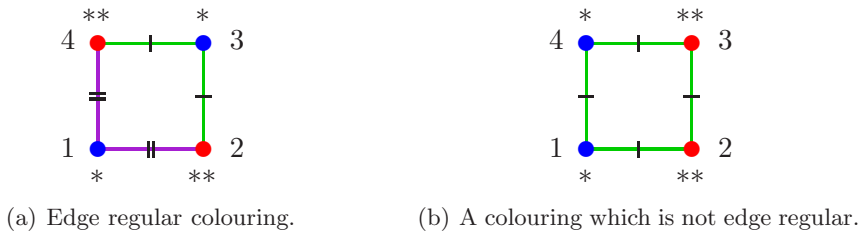
\begin{figure}[ht!]
  \centering
  \subfigure[Edge regular colouring.]{\label{fig:edge_reg}
  \begin{picture}(140,70)(-30,-15) 
\begin{normalsize}
\thicklines
        \put(6,-3){1}
        \put(6,38){4}
        \put(70,-3){2}
        \put(70,38){3}
        \put(20,0){\color{purple}{\line(1,0){40}}}
        \put(20,0){\color{purple}{\line(0,1){40}}}
        \put(60,40){\color{darkgreen}{\line(-1,0){40}}}
        \put(60,40){\color{darkgreen}{\line(0,-1){40}}}
        \put(20,0){\color{blue}{\circle*{6}}}
        \put(20,40){\color{red}{\circle*{6}}}
        \put(60,0){\color{red}{\circle*{6}}}
        \put(60,40){\color{blue}{\circle*{6}}}
        \put(17,-13){$*$}
        \put(54,-13){$**$}
        \put(14,48){$**$}
        \put(57,48){$*$}
                \put(39,-3){\line(0,1){6}}
                \put(41,-3){\line(0,1){6}}
                \put(40,37){\line(0,1){6}}
                \put(17,21){\line(1,0){6}}
                \put(17,19){\line(1,0){6}}
                \put(57,20){\line(1,0){6}}

        \end{normalsize}
        \end{picture}
  } \hspace{20pt}
  \subfigure[A colouring which is not edge regular.]{\label{fig:not_edge_reg}
\begin{picture}(160,70)(-40,-15) 
\begin{normalsize}
\thicklines
        \put(6,-3){1}
        \put(6,38){4}
        \put(70,-3){2}
        \put(70,38){3}
        \put(20,0){\color{darkgreen}{\line(1,0){40}}}
        \put(20,0){\color{darkgreen}{\line(0,1){40}}}
        \put(60,40){\color{darkgreen}{\line(-1,0){40}}}
        \put(60,40){\color{darkgreen}{\line(0,-1){40}}}
        \put(20,0){\color{blue}{\circle*{6}}}
        \put(20,40){\color{blue}{\circle*{6}}}
        \put(60,0){\color{red}{\circle*{6}}}
        \put(60,40){\color{red}{\circle*{6}}}
        \put(17,-13){$*$}
        \put(54,-13){$**$}
        \put(17,48){$*$}
        \put(54,48){$**$}
                \put(40,-3){\line(0,1){6}}
                \put(40,37){\line(0,1){6}}
                \put(17,20){\line(1,0){6}}
                \put(57,20){\line(1,0){6}}
        \end{normalsize}
\end{picture}
  }
  \caption{An edge regular colouring and one which is not edge regular.}
  \label{fig:edge_reg_ex}
\end{figure}

We provide a simple example for illustration. While the colouring in
Figure \ref{fig:edge_reg} is edge regular (both green edges (single
dash) connect a blue vertex (single asterisk) to a red one (two
asterisks), and the same is true for the purple edges (two dashes)),
the colouring in Figure \ref{fig:not_edge_reg} is not, as the green
edges appear between different pairs of vertex colours.

\subsection{Models Represented by Vertex Regular Colourings}

Vertex regular colourings are of relevance to the estimation of a
non-zero mean vector $\mu$ in a $\mathcal{N}_{|V|}(\mu, \Sigma)$
distribution if $\mu$ is subject to equality constraints and
$\Sigma^{-1}$ is restricted to lie inside
$\mathcal{S}^{+}(\mathcal{V},\mathcal{E})$ or inside
$\mathcal{R}^{+}(\mathcal{V},\mathcal{E})$ for some coloured graph
$\mathcal{G}=(\mathcal{V},\mathcal{E})$.

\begin{proposition}[\cite{gehrmann}] \label{prop_vertex_reg}
Let $\mathcal{G} = (\mathcal{V},\mathcal{E}) \in \mathcal{C}_{V}$
and let $\mathcal{M}$ be a partition of $V$. For $\alpha \in V$ let
$v_{\alpha}$ denote the set in $\mathcal{M}$ which contains $\alpha$
and let
\begin{equation*}
\Omega=\Omega(\mathcal{M}) = \{(x_{\alpha})_{\alpha \in V} \in
\mathbb{R}^{V}: x_{\alpha} = x_{\beta}\  \mbox{whenever} \ \alpha
\equiv \beta \ (\mathcal{M}) \}
\end{equation*}
Further let $(Y^{i})_{1 \leq i \leq n}$ be a sample of independent
identically distributed observations $Y^{i} \sim
\mathcal{N}_{|V|}(\mu,\Sigma)$ with $\mu$ restricted to lie inside
$\Omega$. Then the following are equivalent
\begin{enumerate}[(i)]
\itemsep=0pt
\parskip=0pt
\item the likelihood function based on $(y^{i})_{1 \leq i \leq n}$ is
maximised in $\mu$ by the least-squares estimator $\mu^{*}$ for all
$\Sigma$  with $\Sigma^{-1} \in
\mathcal{S}^+(\mathcal{V},\mathcal{E})$ or with $\Sigma^{-1} \in
\mathcal{R}^+(\mathcal{V},\mathcal{E})$ where
\begin{equation}
\mu^*_\alpha=   \frac{\sum_{i=1}^{n}\sum_{\beta \in
v_{\alpha}}y^{i}_{\beta}}{|v_{\alpha}|n} \label{mle2}
\end{equation}
\item $\mathcal{M}$ is finer than $\mathcal{V}$ and $(\mathcal{M},\mathcal{E})$ is vertex regular.
\end{enumerate}
\end{proposition}

For the definition of a vertex regular colouring we require the
concept of an \emph{equitable partition}, first defined
in~\cite{sachs}. For an undirected graph $G=(V,E)$, a vertex
colouring $\mathcal{V}$ of $V$ is called \emph{equitable with
respect to $G$} if for all $v, w \in \mathcal{V}$ and all $\alpha,
\beta \in v$, we have $|\textup{ne}(\alpha) \cap w| =
|\textup{ne}(\beta) \cap w|$. \emph{Vertex regular} graph colourings
are the analogue to equitable partitions for vertex and edge
coloured graphs.

\begin{definition} \label{def_vertex_reg}
For $\mathcal{G} = (\mathcal{V},\mathcal{E}) \in \mathcal{C}_{V}$
let the subgraph induced by the edge colour class $u \in
\mathcal{E}$ be denoted by $G^{u}=(V,u)$. We say that
$(\mathcal{V},\mathcal{E})$ is \emph{vertex regular} if
$\mathcal{V}$ is equitable with respect to $G^{u}$ for all $u \in
\mathcal{E}$.
\end{definition}

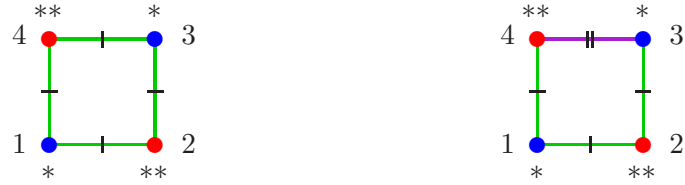
\begin{figure}[ht!]
  \centering
  \subfigure[Vertex regular colouring.]{\label{fig:vertex_reg}
  \begin{picture}(140,70)(-30,-15) 
\begin{normalsize}
\thicklines
        \put(6,-3){1}
        \put(6,38){4}
        \put(70,-3){2}
        \put(70,38){3}
        \put(20,0){\color{darkgreen}{\line(1,0){40}}}
        \put(20,0){\color{darkgreen}{\line(0,1){40}}}
        \put(60,40){\color{darkgreen}{\line(-1,0){40}}}
        \put(60,40){\color{darkgreen}{\line(0,-1){40}}}
        \put(20,0){\color{blue}{\circle*{6}}}
        \put(20,40){\color{red}{\circle*{6}}}
        \put(60,0){\color{red}{\circle*{6}}}
        \put(60,40){\color{blue}{\circle*{6}}}
        \put(17,-13){$*$}
        \put(54,-13){$**$}
        \put(14,48){$**$}
        \put(57,48){$*$}
                \put(40,-3){\line(0,1){6}}
                \put(40,37){\line(0,1){6}}
                \put(17,20){\line(1,0){6}}
                \put(57,20){\line(1,0){6}}

        \end{normalsize}
        \end{picture}
  } \hspace{20pt}
    \subfigure[A colouring which is not vertex regular.]{\label{fig:not_vertex_reg}
\begin{picture}(165,70)(-40,-15) 
\begin{normalsize}
\thicklines
        \put(6,-3){1}
        \put(6,38){4}
        \put(70,-3){2}
        \put(70,38){3}
        \put(20,0){\color{darkgreen}{\line(1,0){40}}}
        \put(20,0){\color{darkgreen}{\line(0,1){40}}}
        \put(60,40){\color{purple}{\line(-1,0){40}}}
        \put(60,40){\color{darkgreen}{\line(0,-1){40}}}
        \put(20,0){\color{blue}{\circle*{6}}}
        \put(20,40){\color{red}{\circle*{6}}}
        \put(60,0){\color{red}{\circle*{6}}}
        \put(60,40){\color{blue}{\circle*{6}}}
        \put(17,-13){$*$}
        \put(54,-13){$**$}
        \put(14,48){$**$}
        \put(57,48){$*$}
                \put(40,-3){\line(0,1){6}}
                \put(41,37){\line(0,1){6}}
                \put(39,37){\line(0,1){6}}
                \put(17,20){\line(1,0){6}}
                \put(57,20){\line(1,0){6}}

        \end{normalsize}
        \end{picture}
  }  
  \caption{A vertex regular colouring and one which is not vertex regular.}
  \label{fig:vert_reg_ex}
\end{figure}

While the colouring in Figure \ref{fig:vertex_reg} is vertex
regular, the colouring in Figure \ref{fig:not_vertex_reg} is not.
The former has only one edge colour class, so that it is vertex
regular if and only if its vertex colouring is equitable with
respect to $G$, which it is. The colouring on the right cannot be
vertex regular as while vertex 4 is incident to a purple edge (two
dashes), vertex 2 isn't even though they are of the same colour.

\subsection{Models Represented by Regular Colourings}

RCON or RCOR models with restrictions represented by colourings
which are both edge regular and vertex regular combine the
properties of both model classes. It can be shown that the
colourings of such models are precisely those which in the
terminology of~\cite{siemons} are \emph{regular}:


\begin{definition}[\cite{siemons}] \label{def_reg_col}
For $\mathcal{G} = (\mathcal{V},\mathcal{E}) \in \mathcal{C}_{V}$,
$(\mathcal{V},\mathcal{E})$ is \emph{regular} if
\begin{enumerate}[(i)]
\item every pair of equally coloured edges in $\mathcal{E}$ connects
the same vertex colour classes in $\mathcal{V}$;
\item every pair of equally coloured vertices in $\mathcal{V}$ has the
same degree in every edge colour class in $\mathcal{E}$.
\end{enumerate}
\end{definition}

By the above, the colourings shown in Figure \ref{fig:not_edge_reg}
and Figure \ref{fig:not_vertex_reg} cannot be regular. While the
colouring given in Figure \ref{fig:vertex_reg} is regular, the
colouring in Figure \ref{fig:edge_reg} is not.

\subsection{Models Represented by Permutation-Generated Colourings}

\emph{Permutation-generated} colourings are a special instance of
regular colourings (for a proof see \citet{gehrmann}), and thus by
definition also of edge regular and vertex regular colourings. They
represent models in which equality constraints on the parameters are
induced by permutation symmetry and allow a particularly simple
maximisation of the likelihood function. In brief, maximum
likelihood estimates can be obtained by standard methods for
unconstrained models after taking averages within colour
classes~\cite{lauritzen_sym}. Further, models represented by
permutation-generated colourings form the only model class discussed
here which restricts $\Sigma^{-1}$ and $\Sigma$ in the same fashion.

The corresponding models are defined through distribution invariance
under a permutation group $\Gamma$ acting on the variable labels
$V$. If $Y \sim \mathcal{N}_{|V|}(0,\Sigma)$, then permutations
acting on $V$ simultaneously permute rows and columns of $\Sigma$ so
that the distribution of $Y$ is invariant under $\Gamma \subseteq
S(V)$ if and only if
\begin{equation} \label{eq_sigma}
P(\sigma)\Sigma P(\sigma)^{T} = \Sigma \ \ \ \ \Longleftrightarrow \
\ \ \ P(\sigma)\Sigma = \Sigma P(\sigma) \ \ \ \ \Longleftrightarrow
\ \ \ \  P(\sigma)\Sigma^{-1} = \Sigma^{-1} P(\sigma)
\end{equation}
for all $\sigma \in \Gamma$, where for $\alpha, \beta \in V$,
$P(\sigma)_{\alpha\beta} = 1$ if and only if $\sigma$ maps $\beta$
to $\alpha$ and zero otherwise. A necessary condition for
equation~(\ref{eq_sigma}) to hold is that the zero entries in
$\Sigma^{-1}$ are preserved for all $\Sigma$ in the model and all
$\sigma \in \Gamma$. Thus if the distribution of $Y$ is assumed to
lie in the graphical Gaussian model represented by graph $G=(V,E)$,
by equation~(\ref{eq_cip}), we in particular require that $\Gamma
\subseteq \textup{Aut}(G)$.

Therefore, in the notation in~\cite{lauritzen_sym}, a graphical
Gaussian model with conditional independence structure represented
by graph $G$ which is permutation invariant under group $\Gamma
\subseteq \textup{Aut}(G)$ is given by assuming
\begin{equation*}
\Sigma^{-1} \in \mathcal{S}^{+}(G) \cap \mathcal{S}^{+}(\Gamma)
\end{equation*}
where $\mathcal{S}^{+}(\Gamma)$ is the set of
positive definite matrices $\Sigma$ satisfying the equivalent conditions in equation (\ref{eq_sigma}).

By definition, permutation invariant models place constraints on all
model parameters and thus in particular on concentrations and
partial correlations, which they restrict in the same fashion. Thus
symmetry constraints in permutation invariant models can be
represented by a vertex and edge colouring
$(\mathcal{V},\mathcal{E})$ of $G$ given by the orbits of $\Gamma$
in $V$ and $E$ respectively, \textit{i.e.\/},~by giving two vertices
$\alpha, \beta \in V$ the same colour whenever there exists $\sigma
\in \Gamma$ mapping $\alpha$ to $\beta$, and similarly for the
edges. We term such colourings \emph{permutation-generated},
formally defined below.

\begin{definition} \label{def_perm_gen}
For $\mathcal{G} = (\mathcal{V},\mathcal{E}) \in \mathcal{C}_{V}$
with underlying uncoloured graph $G=(V,E)$ we say that
$(\mathcal{V},\mathcal{E})$ is \emph{permutation-generated} if there
exists a group $\Gamma \subseteq \textup{Aut}(G)$ acting on $V$ such
that $\mathcal{V}$ and $\mathcal{E}$ are given by the orbits of
$\Gamma$ in $V$ and $E$ respectively.
\end{definition}

The following example illustrates that in addition to the
aforementioned desirable statistical properties, models with
permutation-generated restrictions allow a very intuitive
interpretation.

\begin{example} \label{rcop_example}
The data, commonly referred to as \emph{Fret's heads}, is concerned
with the head dimensions of 25 pairs of first and second
sons \citep{frets,mardia}. Previous analyses \citep{whittaker} support
a model represented by the graph in Figure \ref{fig:four_cycle},
where $L_{i}$ and $B_{i}$ denote the head length and head breadth of
son $i$ for $i=1,2$.~\cite{lauritzen_sym} showed the model generated
by $\Gamma=\langle(B_{1}B_{2})(L_{1}L_{2}) \rangle$, corresponding
to permuting the two sons, represented by the first graph in Figure
\ref{fig:rcop_four_cycle} to be an excellent fit.

Another model with constraints generated by permutation symmetry
which fits the data very well is the complete symmetry model,
generated by $\Gamma= S(V)$, which is represented by the second
coloured graph in Figure \ref{fig:rcop_four_cycle}. Interestingly,
it is further favourable over the former with regards to parameter
estimation. While the graph on the left in Figure
\ref{fig:rcop_four_cycle} is non-decomposable and symmetry arguments
combined with results in~\cite{buhl} give that at least 2
observations are required for almost sure existence of
$\hat{\Sigma}$, see also \citet{uhler}, the complete symmetry model
only requires one observation for $\hat{\Sigma}$ to exist almost
surely.

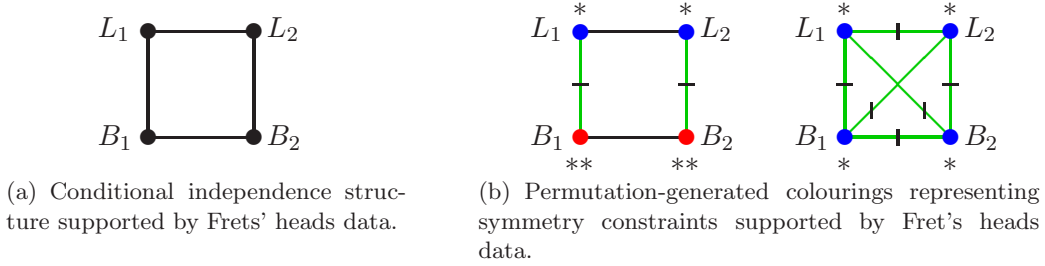
\begin{figure}[ht!]
  \centering
  \subfigure[Conditional independence structure supported by Frets' heads data.]{\label{fig:four_cycle}
  \begin{picture}(144,70)(20,-10) 
\begin{normalsize}
\thicklines
        \put(51,-3){$B_{1}$}
        \put(51,38){$L_{1}$}
        \put(115,-3){$B_{2}$}
        \put(115,38){$L_{2}$}
        \put(70,0){\line(1,0){40}}
        \put(70,0){\line(0,1){40}}
        \put(110,40){\line(-1,0){40}}
        \put(110,40){\line(0,-1){40}}
        \put(70,0){\circle*{6}}
        \put(70,40){\circle*{6}}
        \put(110,0){\circle*{6}}
        \put(110,40){\circle*{6}}
        \end{normalsize}
        \end{picture}
  } \hspace{20pt}
  \subfigure[Permutation-generated colourings representing symmetry constraints
supported by Fret's heads data.]{\label{fig:rcop_four_cycle}
\begin{picture}(205,70)(35,-10) 
\begin{normalsize}
\thicklines
        \put(51,-3){$B_{1}$}
        \put(51,38){$L_{1}$}
        \put(115,-3){$B_{2}$}
        \put(115,38){$L_{2}$}
        \put(70,0){\line(1,0){40}}
        \put(70,0){\color{darkgreen}{\line(0,1){40}}}
        \put(110,40){\line(-1,0){40}}
        \put(110,40){\color{darkgreen}{\line(0,-1){40}}}
        \put(70,0){\color{red}{\circle*{6}}}
        \put(70,40){\color{blue}{\circle*{6}}}
        \put(110,0){\color{red}{\circle*{6}}}
        \put(110,40){\color{blue}{\circle*{6}}}
        \put(64,-13){$**$}
        \put(104,-13){$**$}
        \put(67,46){$*$}
        \put(107,46){$*$}
                \put(67,20){\line(1,0){6}}
                \put(107,20){\line(1,0){6}}

        \put(151,-3){$B_{1}$}
        \put(151,38){$L_{1}$}
        \put(215,-3){$B_{2}$}
        \put(215,38){$L_{2}$}
        \put(170,0){\color{darkgreen}{\line(1,0){40}}}
        \put(170,0){\color{darkgreen}{\line(0,1){40}}}
        \put(170,0){\color{darkgreen}{\line(1,1){40}}}
        \put(210,40){\color{darkgreen}{\line(-1,0){40}}}
        \put(210,40){\color{darkgreen}{\line(0,-1){40}}}
        \put(210,0){\color{darkgreen}{\line(-1,1){40}}}
        \put(170,0){\color{blue}{\circle*{6}}}
        \put(170,40){\color{blue}{\circle*{6}}}
        \put(210,0){\color{blue}{\circle*{6}}}
        \put(210,40){\color{blue}{\circle*{6}}}
        \put(167,-13){$*$}
        \put(207,-13){$*$}
        \put(167,46){$*$}
        \put(207,46){$*$}
                \put(167,20){\line(1,0){6}}
                \put(207,20){\line(1,0){6}}
                \put(190,-3){\line(0,1){6}}
                \put(190,37){\line(0,1){6}}
                \put(200,7){\line(0,1){6}}
                \put(180,7){\line(0,1){6}}

        \end{normalsize}
\end{picture}
  }
  \caption{Frets' heads example.}
  \label{fig:rcop_ex}
\end{figure}
\end{example}

\subsection{Relations Between Model Classes}

\begin{figure}[ht!]
  \centering
\begin{picture}(150,125)(-10,-20)
\put(40,60){\oval(150,75)} \put(100,45){\oval(150,75)}
\put(70,50){\circle{45}} \put(-15,55){$P$} \put(95,62){$R$}
\put(66,45){$\Pi$} \put(145,35){$B$}
\end{picture}
\vspace{-30pt}
  \caption{Structural relations between colouring classes.}
  \label{fig:relations}
  \vspace{5pt}
\end{figure}
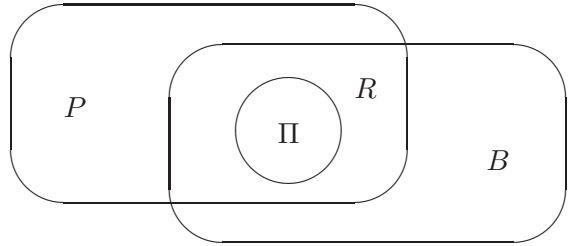

Let $B$, $P$, $R$ and $\Pi$ denote the sets of edge regular, vertex
regular, regular and permutation-generated colourings respectively.
The structural relations between colouring classes are summarised in
the diagram displayed in Figure \ref{fig:relations}. In fact, we
already saw examples of colourings in three of the four disjoint
sets in the diagram. The colouring displayed in Figure
\ref{fig:not_edge_reg} lies in $P \setminus R$, Figure
\ref{fig:rcop_four_cycle} shows a colouring in $\Pi$ and the
colouring in Figure \ref{fig:not_vertex_reg} lies in $B \setminus
R$. (A graph colouring in $R \setminus \Pi$ is given by
$(\mathcal{V},\mathcal{E})$ with $V=[11]$,
$\mathcal{V}=\{\{1,2,3\},\{4,5,6,7,8,9\},\{10,11\}\}$ and
$\mathcal{E}=\{\{14,15,26,27,38,39\},\{(4,10),(5,10),\linebreak(6,10),(7,11),(8,11),(9,11)\}\}$
where $(i,j)$ denotes an edge between vertices $i$ and $j$.)

By Proposition \ref{prop_edge_reg}, a graph colouring yields the
same model restrictions representing an RCON model as it does
representing an RCOR model if and only if it lies in $B$. Therefore,
the model type only needs to be specified whenever a graph colouring
lies in $P \setminus B$. Put formally:

\begin{proposition} \label{prop_mod_classes}
For $\mathcal{G}=(\mathcal{V},\mathcal{E}) \in \mathcal{C}_{V}$,
$\mathcal{S}^{+}(\mathcal{V},\mathcal{E}) =
\mathcal{R}^{+}(\mathcal{V},\mathcal{E})$ for
$(\mathcal{V},\mathcal{E}) \in B$ and
$\mathcal{S}^{+}(\mathcal{V},\mathcal{E}) \not=
\mathcal{R}^{+}(\mathcal{V},\mathcal{E})$ for
$(\mathcal{V},\mathcal{E}) \in P \setminus B$. Thus if for $X \in
\{B,P,R,\Pi\}$ we let $\mathcal{S}^{+}_{X}$ denote the set of RCON
models represented by graphs in $X$ and similarly for
$\mathcal{R}^{+}_{X}$ and RCOR models, then
\begin{equation*}
\mathcal{S}^{+}_{B} = \mathcal{R}^{+}_{B}, \ \ \ \ \ \ \ \ \
\mathcal{S}^{+}_{R} = \mathcal{R}^{+}_{R}, \ \ \ \ \ \ \ \ \
\mathcal{S}^{+}_{\Pi} = \mathcal{R}^{+}_{\Pi}, \ \ \ \ \ \ \ \ \
\mathcal{S}^{+}_{P} \not= \mathcal{R}^{+}_{P}
\end{equation*}
giving rise to five model classes lying strictly within
$\mathcal{S}^{+}_{V}$ and $\mathcal{R}^{+}_{V}$ with
\begin{equation*}
\mathcal{S}^{+}_{\Pi} \subset \mathcal{S}^{+}_{R} \subset
\mathcal{S}^{+}_{B}, \ \ \ \ \ \ \ \ \ \mathcal{S}^{+}_{B} \cap
\mathcal{S}^{+}_{P} = \mathcal{S}^{+}_{B} \cap \mathcal{R}^{+}_{P} =
\mathcal{S}^{+}_{R}.
\end{equation*}
\end{proposition}

Let $B_{V}$, $P_{V}$, $R_{V}$ and $\Pi_{V}$ denote the sets of graph
colourings inside $B$, $P$, $R$ and $\Pi$ with vertex set $V$. By
Proposition \ref{prop_mod_classes}, there are five corresponding
model classes: $\mathcal{S}^{+}_{B_{V}}$, $\mathcal{S}^{+}_{P_{V}}$,
$\mathcal{R}^{+}_{P_{V}}$, $\mathcal{S}^{+}_{R_{V}}$ and
$\mathcal{S}^{+}_{\Pi_{V}}$. For illustration we give the
corresponding model class sizes for $V=[4]$ together with
$|M_{[4]}|$ in Table \ref{tab_sizes}.

\begin{table}[ht!]
\begin{center}
\begin{tabular}{lcccccc}
\toprule
    \textbf{Model class}                & $\mathcal{S}^{+}_{[4]}, \mathcal{R}^{+}_{[4]}$ & $\mathcal{S}^{+}_{B_{[4]}}$ & $\mathcal{S}^{+}_{P_{[4]}}, \mathcal{R}^{+}_{P_{[4]}}$ & $\mathcal{S}^{+}_{R_{[4]}}$ & $\mathcal{S}^{+}_{\Pi_{[4]}}$ & $M_{[4]}$ \\[.5ex]
\midrule
Size                 &  13,155 &  3065 &  1380 & 251 & 251 & 64\\
\bottomrule
\end{tabular}
\end{center}
\vspace{-10pt} \caption{Model set sizes for $V=[4]$.}
\label{tab_sizes}
\end{table}

The relative sizes in Table \ref{tab_sizes} are representative of
the general case: $\mathcal{S}^{+}_{B_{V}}$ is the largest model
class, followed by $\mathcal{S}^{+}_{P_{V}}$,
$\mathcal{R}^{+}_{P_{V}}$ and $\mathcal{S}^{+}_{R_{V}}$.
$|\mathcal{S}^{+}_{P_{V}}| = |\mathcal{R}^{+}_{P_{V}}|$ will
generally be considerably smaller than $\mathcal{S}^{+}_{B_{V}}$ as
the defining conditions of $P_{V}$ are far more restrictive than
those for $B_{V}$. $|\mathcal{S}^{+}_{\Pi_{V}}|$ is the smallest
class of the four for all $V$, however may equal
$|\mathcal{S}^{+}_{R_{V}}|$ for some $V$, as for example for
$V=[4]$.

\section{Structures of Model Classes}

Below we show that each model class defined above forms a complete
non-distributive lattice, starting with $\mathcal{S}^{+}_{B_{V}}$
and $\mathcal{S}^{+}_{\Pi_{V}}$ as their structure turns out most
tractable. For brevity, we only outline results for the remaining
model classes.

\subsection{Models Represented by Edge Regular Colourings}

\begin{proposition} \label{prop_b}
$B_{V}$ is stable under the meet operation $\land$ in $\langle
\mathcal{C}_{V}; \preceq \rangle$ given in
equation~(\ref{meet_join_c}).
\end{proposition}

\noindent \emph{Proof:} Let $\mathcal{G} =
(\mathcal{V}_{\mathcal{G}},\mathcal{E}_{\mathcal{G}}), \mathcal{H} =
(\mathcal{V}_{\mathcal{H}},\mathcal{E}_{\mathcal{H}}) \in B_{V}$.
$\mathcal{G} \land \mathcal{H} = (\mathcal{V}_{\mathcal{G} \land
\mathcal{H}}, \mathcal{E}_{\mathcal{G} \land \mathcal{H}})$ is
obtained from $\mathcal{G}$ and $\mathcal{H}$ by dropping of edge
colour classes and merging of colour classes. The only operation
potentially leading to $\mathcal{G} \land \mathcal{H}$ lying outside
of $B_{V}$ is merging of edge colour classes.

So let $\alpha\beta$ and $\gamma\delta$ be two edges in $\mathcal{G}
\land \mathcal{H}$ of equal colour. Then there exists a sequence
$\alpha_{0}\beta_{0}, \ldots, \alpha_{k}\beta_{k}$ in
$E_{\mathcal{G} \wedge \mathcal{H}}$ such that $\alpha_{0}\beta_{0}
= \alpha\beta$, $\alpha_{k}\beta_{k} = \gamma\delta$, and
$\alpha_{i-1}\beta_{i-1} \equiv \alpha_{i}\beta_{i} \
(\mathcal{E}_{\mathcal{G}})$ or $\alpha_{i-1}\beta_{i-1} \equiv
\alpha_{i}\beta_{i} \ (\mathcal{E}_{\mathcal{H}})$ for $1 \leq i
\leq k$. As both $\mathcal{G}$ and $\mathcal{H}$ have edge regular
colourings, $\alpha_{i-1}\beta_{i-1}$ and $\alpha_{i}\beta_{i}$
connect the same vertex colour classes in the graph in which they
are of equal colour, which we denote by
$\{\alpha_{i-1},\beta_{i-1}\} \equiv \{\alpha_{i},\beta_{i}\} \
(\mathcal{V}_{\mathcal{X}})$ with $\mathcal{X} \in
\{\mathcal{G},\mathcal{H}\}$. This gives that $\{\alpha, \beta\} =
\{\alpha_{0}, \beta_{0}\} \equiv \{\alpha_{k}, \beta_{k}\} =
\{\gamma, \delta\} \ (\mathcal{V}_{\mathcal{G}} \vee
\mathcal{V}_{\mathcal{H}})$. As $\mathcal{V}_{\mathcal{G}} \vee
\mathcal{V}_{\mathcal{H}} = \mathcal{V}_{\mathcal{G} \wedge
\mathcal{H}}$, $\alpha\beta$ and $\gamma\delta$ connect the same
vertex colour classes in $\mathcal{G} \land \mathcal{H}$. \whitebox

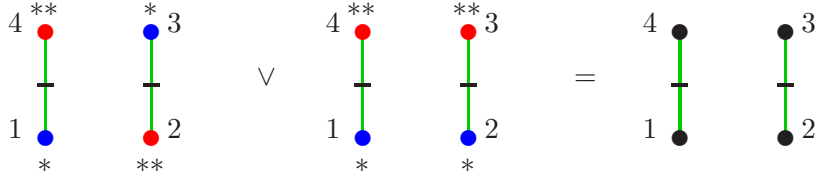
\begin{figure}[ht!]
  \centering
\begin{picture}(320,70)(10,-10)
\begin{normalsize}
\thicklines
         \put(6,0){1}
        \put(6,40){4}
        \put(66,0){2}
        \put(66,40){3}
        \put(20,0){\color{darkgreen}{\line(0,1){40}}}
        \put(60,40){\color{darkgreen}{\line(0,-1){40}}}
        \put(20,0){\color{blue}{\circle*{6}}}
        \put(20,40){\color{red}{\circle*{6}}}
        \put(60,0){\color{red}{\circle*{6}}}
        \put(60,40){\color{blue}{\circle*{6}}}
        \put(17,-13){$*$}
        \put(54,-13){$**$}
        \put(14,46){$**$}
        \put(57,46){$*$}
                \put(17,20){\line(1,0){6}}
                \put(57,20){\line(1,0){6}}

                \put(100,20){$\lor$}

                \put(126,0){1}
        \put(126,40){4}
        \put(186,0){2}
        \put(186,40){3}
        \put(140,0){\color{darkgreen}{\line(0,1){40}}}
        \put(180,40){\color{darkgreen}{\line(0,-1){40}}}
        \put(140,0){\color{blue}{\circle*{6}}}
        \put(140,40){\color{red}{\circle*{6}}}
        \put(180,0){\color{blue}{\circle*{6}}}
        \put(180,40){\color{red}{\circle*{6}}}
        \put(137,-13){$*$}
        \put(177,-13){$*$}
        \put(134,46){$**$}
        \put(174,46){$**$}
                \put(137,20){\line(1,0){6}}
                \put(177,20){\line(1,0){6}}

                \put(220,20){=}

                \put(246,0){1}
        \put(246,40){4}
        \put(306,0){2}
        \put(306,40){3}
        \put(260,0){\color{darkgreen}{\line(0,1){40}}}
        \put(300,40){\color{darkgreen}{\line(0,-1){40}}}
        \put(260,0){\circle*{6}}
        \put(260,40){\circle*{6}}
        \put(300,0){\circle*{6}}
        \put(300,40){\circle*{6}}
                \put(297,20){\line(1,0){6}}
                \put(257,20){\line(1,0){6}}

        \end{normalsize}
\end{picture}
  \caption{$B_{[4]}$ is not stable under $\lor$.}
  \label{fig:join_edge_reg}
\end{figure}

Graphs $\mathcal{G}_{4}$, $\mathcal{G}_{5}$ and $\mathcal{G}_{4}
\land \mathcal{G}_{5}$ in Figure \ref{fig:join_col_ex} illustrate
the stability of $B_{V}$ under $\wedge$. That $B_{V}$ is generally
not stable under the join operation $\lor$ in
equation (\ref{meet_join_c}) is established by the example in Figure
\ref{fig:join_edge_reg}.

\begin{proposition} \label{prop_barbell_meet}
Let $\mathcal{G} = (\mathcal{V},\mathcal{E}) \in \mathcal{C}_{V}$
and let $\mathcal{E}_{B}$ be the partition of $E$ which puts
$\alpha\beta, \gamma\delta \in E$ in the same set whenever they
connect the same vertex color classes. Then $\mathcal{G}$ has a
supremum in $B_{V}$, given by
\begin{equation*}
s_{B}(\mathcal{G}) = (\mathcal{V},\mathcal{E} \wedge
\mathcal{E}_{B})
\end{equation*}
\end{proposition}

\noindent \emph{Proof:} The claim is trivially true for $\mathcal{G}
\in B_{V}$. If $\mathcal{G} \in \mathcal{C}_{V}\setminus B_{V}$, an
edge regular colouring cannot be achieved through splitting vertex
colour classes or adding edge colour classes. The only effective
manipulation is therefore the splitting of edge colour classes. The
coarsest partition which is finer than $\mathcal{E}$ and gives an
edge regular colouring is $\mathcal{E} \wedge \mathcal{E}_{B}$, as
it splits edge colour classes only if they connect different vertex
colour classes.
\whitebox

We deduce:

\begin{theorem} \label{thm_b}
$\mathcal{S}^{+}_{B_{V}}$ is a complete non-distributive lattice
with respect to model inclusion. The meet operation is induced by
the meet operation in $\langle \mathcal{C}_{V}; \leq \rangle$ given
in equation~(\ref{meet_join_c}). The join of two models represented
by graphs $\mathcal{G},\mathcal{H} \in B_{V}$ is represented by the
graph $s_{B}(\mathcal{G} \lor \mathcal{H})$.
\end{theorem}

\noindent \emph{Proof}: Proposition \ref{prop_b} implies that $\inf
H$ exists for all finite $H \subseteq B_{V}$, which by Lemma
\ref{lemma_complete_lattice} gives that $B_{V}$ is a complete
lattice, with the same meet operation as $\langle \mathcal{C}_{V};
\preceq \rangle$. As the zero and unit in $\langle \mathcal{C}_{V};
\preceq \rangle$ have edge regular colourings, they are also the
zero and unit in $B_{V}$.

By definition, for $\mathcal{G}, \mathcal{H} \in B_{V}$, the join
$\mathcal{K}$ of $\mathcal{G}$ and $\mathcal{H}$ in $B_{V}$ is the
smallest graph with respect to partial ordering $\preceq$ which
satisfies $\mathcal{K} \succeq \mathcal{G}$, $\mathcal{K} \succeq
\mathcal{H}$ and $\mathcal{K} \in B_{V}$. The supremum $\mathcal{G}
\vee \mathcal{H}$ of $\mathcal{G}$ and $\mathcal{H}$ in $\langle
\mathcal{C}_{V}; \preceq \rangle$ is the smallest graph satisfying
the first two relations so that $\mathcal{K}$ is in fact the
smallest graph satisfying $\mathcal{K} \succeq (\mathcal{G} \lor
\mathcal{H})$ and $\mathcal{K} \in B_{V}$. Thus $\mathcal{G}
\lor_{B} \mathcal{H}$ is given by the supremum of $\mathcal{G} \lor
\mathcal{H}$ in $B_{V}$, which by Proposition
\ref{prop_barbell_meet} equals $s_{B}(\mathcal{G} \lor
\mathcal{H})$.

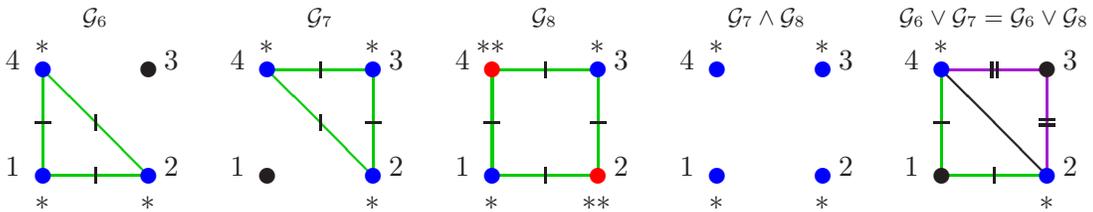
\begin{figure}[ht!]
  \centering
  \vspace{90pt}
\begin{picture}(400,95)(10,-30)
\begin{normalsize}
\thicklines
        \put(6,90){1}
        \put(6,130){4}
        \put(66,90){2}
        \put(66,130){3}
        \put(20,90){\color{darkgreen}{\line(1,0){40}}}
        \put(20,90){\color{darkgreen}{\line(0,1){40}}}
        \put(60,90){\color{darkgreen}{\line(-1,1){40}}}
        \put(20,90){\color{blue}{\circle*{6}}}
        \put(20,130){\color{blue}{\circle*{6}}}
        \put(60,90){\color{blue}{\circle*{6}}}
        \put(60,130){{\circle*{6}}}
        \put(17,77){$*$}
        \put(57,77){$*$}
        \put(17,136){$*$}
                \put(40,87){\line(0,1){6}}
                \put(40,107){\line(0,1){6}}
                \put(17,110){\line(1,0){6}}

                \put(35,147){\footnotesize $\mathcal{G}_{6}$}

        \put(91,90){1}
        \put(91,130){4}
        \put(151,90){2}
        \put(151,130){3}
        \put(105,130){\color{darkgreen}{\line(1,0){40}}}
        \put(145,90){\color{darkgreen}{\line(0,1){40}}}
        \put(145,90){\color{darkgreen}{\line(-1,1){40}}}
        \put(105,90){{\circle*{6}}}
        \put(105,130){\color{blue}{\circle*{6}}}
        \put(145,90){\color{blue}{\circle*{6}}}
        \put(145,130){\color{blue}{\circle*{6}}}
        \put(142,136){$*$}
        \put(142,77){$*$}
        \put(102,136){$*$}
                \put(125,127){\line(0,1){6}}
                \put(125,107){\line(0,1){6}}
                \put(142,110){\line(1,0){6}}

                \put(120,147){\footnotesize $\mathcal{G}_{7}$}

                \put(176,90){1}
        \put(176,130){4}
        \put(236,90){2}
        \put(236,130){3}
        \put(190,90){\color{darkgreen}{\line(1,0){40}}}
        \put(190,90){\color{darkgreen}{\line(0,1){40}}}
        \put(230,130){\color{darkgreen}{\line(-1,0){40}}}
        \put(230,130){\color{darkgreen}{\line(0,-1){40}}}
        \put(190,90){\color{blue}{\circle*{6}}}
        \put(190,130){\color{red}{\circle*{6}}}
        \put(230,90){\color{red}{\circle*{6}}}
        \put(230,130){\color{blue}{\circle*{6}}}
        \put(187,77){$*$}
        \put(224,77){$**$}
        \put(184,136){$**$}
        \put(227,136){$*$}
                \put(210,87){\line(0,1){6}}
                \put(210,127){\line(0,1){6}}
                \put(187,110){\line(1,0){6}}
                \put(227,110){\line(1,0){6}}

                \put(205,147){\footnotesize $\mathcal{G}_{8}$}

        \put(261,90){1}
        \put(261,130){4}
        \put(321,90){2}
        \put(321,130){3}
        \put(275,90){\color{blue}{\circle*{6}}}
        \put(275,130){\color{blue}{\circle*{6}}}
        \put(315,90){\color{blue}{\circle*{6}}}
        \put(315,130){\color{blue}{\circle*{6}}}
        \put(272,77){$*$}
        \put(312,77){$*$}
        \put(272,136){$*$}
        \put(312,136){$*$}

                \put(279,147){\footnotesize $\mathcal{G}_{7} \wedge \mathcal{G}_{8}$}

                \put(346,90){1}
        \put(346,130){4}
        \put(406,90){2}
        \put(406,130){3}
        \put(360,130){{\line(1,-1){40}}}
        \put(360,130){\color{darkgreen}{\line(0,-1){40}}}
        \put(360,90){\color{darkgreen}{\line(1,0){40}}}
        \put(400,130){\color{purple}{\line(-1,0){40}}}
        \put(400,130){\color{purple}{\line(0,-1){40}}}
        \put(360,90){{\circle*{6}}}
        \put(360,130){\color{blue}{\circle*{6}}}
        \put(400,90){\color{blue}{\circle*{6}}}
        \put(400,130){{\circle*{6}}}
        \put(397,77){$*$}
        \put(357,136){$*$}
                \put(379,127){\line(0,1){6}}
                \put(381,127){\line(0,1){6}}
                \put(380,87){\line(0,1){6}}
                \put(397,111){\line(1,0){6}}
                \put(397,109){\line(1,0){6}}
                \put(357,110){\line(1,0){6}}

                \put(344,147){\footnotesize $\mathcal{G}_{6} \vee \mathcal{G}_{7}=\mathcal{G}_{6} \vee \mathcal{G}_{8}$}

        \end{normalsize}
\end{picture}\\
\vspace{-110pt} \caption{Non-distributivity of $B_{V}$.}
\label{fig_non_dist}
\end{figure}

Non-distributivity of $B_{V}$ is established by observing that
equation~(\ref{eq_distr}) is violated for $a = \mathcal{G}_{6}$, $b
= \mathcal{G}_{7}$ and $c = \mathcal{G}_{8}$ displayed in Figure
\ref{fig_non_dist}, as $\mathcal{G}_{6} = \mathcal{G}_{6} \vee
(\mathcal{G}_{7} \wedge \mathcal{G}_{8}) \not = (\mathcal{G}_{6}
\vee \mathcal{G}_{7}) \wedge (\mathcal{G}_{6} \vee \mathcal{G}_{8})
= \mathcal{G}_{6} \vee \mathcal{G}_{7}$.

The results on the structure of $B_{V}$ naturally translate to the
set of models $\mathcal{S}^{+}_{B_{V}}$, proving the claim.
\whitebox

\subsection{Models Represented by Permutation-Generated Colourings}

Let $\Gamma_{V}$ denote the set of permutation groups acting on $V$.
Then $\langle \Gamma_{V}; \subseteq \rangle$ is a complete
lattice~\cite{schmidt} with meet and join operations given by
$\Gamma_{1} \wedge \Gamma_{2} = \Gamma_{1} \cap \Gamma_{2}$ and
$\Gamma_{1} \vee \Gamma_{2} = \langle \Gamma_{1} \cup \Gamma_{2}
\rangle$. We obtain:

\begin{proposition} \label{prop_pi}
$\Pi_{V}$ is stable under the meet operation $\land$ in $\langle
\mathcal{C}_{V}; \preceq \rangle$ given in
equation~(\ref{meet_join_c}). If\linebreak
$\mathcal{G}=(\mathcal{V}_{\mathcal{G}}$,
$\mathcal{E}_{\mathcal{G}}),\mathcal{H}=(\mathcal{V}_{\mathcal{H}},\mathcal{E}_{\mathcal{H}})
\in \Pi_{V}$ are generated by $\Gamma_{\mathcal{G}},
\Gamma_{\mathcal{H}} \in \Gamma_{V}$, then the colouring of
$\mathcal{G} \land \mathcal{H}$ is generated by
$\Gamma_{\mathcal{G}} \lor \Gamma_{\mathcal{H}}$.
\end{proposition}

\noindent \emph{Proof:} Let $\mathcal{G},\mathcal{H} \in \Pi_{V}$
and $\Gamma_{\mathcal{G}}, \Gamma_{\mathcal{H}} \in \Gamma_{V}$ be
as in the claim. Then $\mathcal{V}_{\mathcal{G}}$ and
$\mathcal{E}_{\mathcal{G}}$ are unions of orbits of
$\Gamma_{\mathcal{G}}$ in $V$ and $E_{\mathcal{G}}$, and similarly
for $\mathcal{V}_{\mathcal{H}}$,  $\mathcal{E}_{\mathcal{H}}$ and
$\Gamma_{\mathcal{H}}$. By definition of the meet operation in
$\langle \mathcal{C}_{V}; \preceq \rangle$, each vertex colour class
in $\mathcal{G} \land \mathcal{H} = (\mathcal{V}_{\mathcal{G} \wedge
\mathcal{H}}$, $\mathcal{E}_{\mathcal{G} \wedge \mathcal{H}})$ can
be expressed as a union of vertex colour classes in
$\mathcal{V}_{\mathcal{G}}$, and as a union of vertex colour classes
in $\mathcal{V}_{\mathcal{H}}$, and similarly for the edges. Thus
the colouring of $\mathcal{G} \land \mathcal{H}$ is invariant under
the action of both groups $\Gamma_{\mathcal{G}}$ and
$\Gamma_{\mathcal{H}}$, and therefore also under
$\Gamma_{\mathcal{G}} \lor \Gamma_{\mathcal{H}}$.


To show that the colouring of $\mathcal{G} \land \mathcal{H}$ is
generated by $\Gamma_{\mathcal{G}} \vee \Gamma_{\mathcal{H}}$, we
need to show that whenever two vertices or edges in $\mathcal{G}
\land \mathcal{H}$ are of the same colour, then there exists $\sigma
\in \Gamma_{\mathcal{G}} \vee \Gamma_{\mathcal{H}}$ which maps one
of them to the other. We present the argument for the edges only, as
it is can be trivially transferred to the vertices. So let
$\alpha\beta$ and $\gamma\delta$ be two edges in $\mathcal{G} \land
\mathcal{H}$ of equal colour. Then there exists a sequence
$\alpha_{0}\beta_{0}, \ldots, \alpha_{k}\beta_{k}$ in
$E_{\mathcal{G} \wedge \mathcal{H}}$ such that $\alpha_{0}\beta_{0}
= \alpha\beta$, $\alpha_{k}\beta_{k} = \gamma\delta$, and
$\alpha_{i-1}\beta_{i-1} \equiv \alpha_{i}\beta_{i} \
(\mathcal{E}_{\mathcal{G}})$ or $\alpha_{i-1}\beta_{i-1} \equiv
\alpha_{i}\beta_{i} \ (\mathcal{E}_{\mathcal{H}})$ for $1 \leq i
\leq k$. There must therefore exist $\sigma_{i} \in
\Gamma_{\mathcal{G}} \cup \Gamma_{\mathcal{H}}$ such that
$\alpha_{i-1}\beta_{i-1}$ is mapped to $\alpha_{i}\beta_{i}$ by
$\sigma_{i}$ for $1 \leq i \leq k$, giving that the product
$\sigma_{k}\ldots \sigma_{1} \in \Gamma_{\mathcal{G}} \vee
\Gamma_{\mathcal{H}}$ maps $\alpha\beta$ to $\gamma\delta$.
\whitebox

Graphs $\mathcal{G}_{4}$, $\mathcal{G}_{5}$ and $\mathcal{G}_{4}
\land \mathcal{G}_{5}$ in Figure \ref{fig:join_col_ex} illustrate
the above result, as each of the graphs lies in $\Pi_{V}$, with
generating groups $\Gamma_{4}=\langle(13)(24)\rangle$, $\Gamma_{5}=
\langle(13)\rangle$ and $\Gamma_{4 \land 5} = \Gamma_{4} \vee
\Gamma_{5}= \langle(13),(24)\rangle$ respectively. Observing that
the join $\mathcal{G}_{4} \lor \mathcal{G}_{5}$, also displayed in
Figure \ref{fig:join_col_ex}, does not lie in $\Pi_{V}$ establishes
that $\Pi_{V}$ is generally not stable under the join operation in
$\langle \mathcal{C}_{V}; \preceq \rangle$. However the following
holds:

\begin{proposition} \label{prop_pi_sup}
For $\mathcal{G}=(\mathcal{V},\mathcal{E}) \in \mathcal{C}_{V}$, let
$\textup{Aut}(\mathcal{V},\mathcal{E}) \leq S(V)$ denote the largest
group leaving $(\mathcal{V},\mathcal{E})$ invariant and let
$(\mathcal{V}_{\textup{Aut}}, \mathcal{E}_{\textup{Aut}})$ denote
the graph colouring of $G=(V,E)$ given by the orbits of
$\textup{Aut}(\mathcal{V},\mathcal{E})$ in $V$ and $E$ respectively.
Then $\mathcal{G}$ has a supremum in $\Pi_{V}$ given by
\vspace{-3pt}
\begin{equation*}
s_{\Pi}(\mathcal{G}) = (\mathcal{V}_{\textup{Aut}},
\mathcal{E}_{\textup{Aut}})
\end{equation*}
\end{proposition}

\noindent \emph{Proof:} The claim is trivially true if $\mathcal{G}
\in \Pi_{V}$. So suppose $\mathcal{G} \in \mathcal{C}_{V} \setminus
\Pi_{V}$. $\mathcal{G}$ is modified to a larger graph by adding edge
colour classes and splitting colour classes. As the former will not
enforce a permutation-generated colouring, to prove the claim we
need to that $(\mathcal{V}_{\textup{Aut}},
\mathcal{E}_{\textup{Aut}})$ is the coarsest refinement of
$(\mathcal{V},\mathcal{E})$ which lies in $\Pi_{V}$. This is clearly
the case as $\textup{Aut}(\mathcal{V},\mathcal{E})$ is the largest
group which leaves $\mathcal{V}$ and $\mathcal{E}$ invariant. \whitebox

\begin{theorem} \label{thm_pi}
$\mathcal{S}^{+}_{\Pi_{V}}$ is a complete non-distributive lattice
with respect to model inclusion. The meet operation is induced by
the meet operation in $\langle \mathcal{C}_{V}; \leq \rangle$ given
in equation~(\ref{meet_join_c}). The join of two models represented
by graphs $\mathcal{G},\mathcal{H} \in \Pi_{V}$ is represented by
the graph $s_{\Pi}(\mathcal{G} \lor \mathcal{H})$.
\end{theorem}

\noindent \emph{Proof}: The proof is analogous to the proof of
Theorem \ref{thm_b}. In brief, Proposition \ref{prop_pi} and Lemma
\ref{lemma_complete_lattice} give that $\Pi_{V}$ is a complete
lattice with meet operation as claimed. As the zero and unit in
$\langle \mathcal{C}_{V}; \preceq \rangle$ have
permutation-generated colourings, with $\Gamma_{0}=S(V)$ and
$\Gamma_{1}=\langle Id \rangle$, they are the zero and unit in
$\Pi_{V}$.

By Proposition \ref{prop_pi_sup}, the join of two graphs
$\mathcal{G}, \mathcal{H} \in \Pi_{V}$ is given by
$s_{\Pi}(\mathcal{G} \vee \mathcal{H})$. The graphs displayed in
Figure \ref{fig_non_dist} establish non-distributivity of $\Pi_{V}$
as each of them has a permutation-generated colouring, with
$\Gamma_{6}=\langle(124)\rangle$, $\Gamma_{7}=\langle(234)\rangle$,
$\Gamma_{8}=\langle(13),(24)\rangle$, $\Gamma_{7 \wedge 8}=S([4])$
and $\Gamma_{6 \vee 7} = \Gamma_{6 \vee 8}  =\langle(24)\rangle$
respectively. The above directly translate to
$\mathcal{S}^{+}_{\Pi_{V}}$, proving the claim. \whitebox

\subsection{Models Represented by Regular and Vertex Regular Colourings}

The structures of $R_{V}$ and $P_{V}$ turn out to be closely related
and it is for that reason that we treat them together. We abstain
from giving explicit proofs of intermediate results for brevity,
however give all facts the reader will require to construct them. We
shall employ the notion of a factor graph:

\begin{definition}[\cite{frey}] \label{def_factor_graph}
Let $f(Y)$ be a function in $Y$ which factorises as $f(Y) = \Pi_{i
\in I} f_{A_{i}}(Y_{A_{i}})$ where $A_{i} \subseteq V$ and
$f_{A_{i}}$ cannot be factorised further for $i \in I$. Then the
\emph{factor graph} of $f$ is the graph $G_{F} = (V \cup F,E_{F})$
with $F=\{f_{A_{i}}\}_{i \in I}$ being the set of factor vertices
and $E_{F} = \{\alpha f_{A_{i}} \mid \alpha \in V, f_{A_{i}} \in F
\text{ with } \alpha \in A_{i}\}$.
\end{definition}

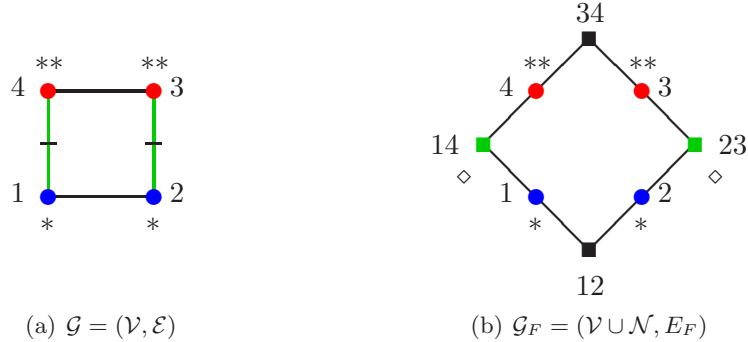
\begin{figure}[ht!]
  \centering
  \subfigure[$\mathcal{G}=(\mathcal{V},\mathcal{E})$]{\label{fig:col_graph}
\begin{picture}(190,105)(65,-28)
\begin{normalsize}
\thicklines
          \put(126,8){$1$}
        \put(126,48){$4$}
        \put(186,8){$2$}
        \put(186,48){$3$}

        \put(137,-3){$*$}
        \put(177,-3){$*$}
        \put(135,58){$**$}
        \put(175,58){$**$}

        \put(140,10){\line(1,0){40}}
        \put(140,10){\color{darkgreen}{\line(0,1){40}}}
        \put(140,50){\line(1,0){40}}
        \put(180,50){\color{darkgreen}{\line(0,-1){40}}}

        \put(140,10){\color{blue}{\circle*{6}}}
        \put(140,50){\color{red}{\circle*{6}}}
        \put(180,10){\color{blue}{\circle*{6}}}
        \put(180,50){\color{red}{\circle*{6}}}

                \put(177,30){\line(1,0){6}}
                \put(137,30){\line(1,0){6}}
        \end{normalsize}
\end{picture}
\vspace{20pt}
  }
\hspace{-20pt}
  \subfigure[$\mathcal{G}_{F}=(\mathcal{V} \cup \mathcal{N},E_{F})$]{\label{fig:col_factor_graph}
\begin{picture}(190,105)(45,-28)
\begin{normalsize}
\thicklines
        \put(106,8){$1$}
        \put(106,48){$4$}
        \put(166,8){$2$}
        \put(166,48){$3$}
        \put(135,-27){$12$}
        \put(189,26){$23$}
        \put(135,76){$34$}
        \put(80,26){$14$}

        \put(100,30){\line(1,1){40}}
        \put(100,30){\line(1,-1){40}}
        \put(180,30){\line(-1,-1){40}}
        \put(180,30){\line(-1,1){40}}

        \put(120,10){\color{blue}{\circle*{6}}}
        \put(120,50){\color{red}{\circle*{6}}}
        \put(160,10){\color{blue}{\circle*{6}}}
        \put(160,50){\color{red}{\circle*{6}}}

        \put(97,27){\scriptsize{\color{darkgreen}{$\blacksquare$}}}
        \put(177,27){\scriptsize{\color{darkgreen}{$\blacksquare$}}}
        \put(137,-13){\scriptsize{{$\blacksquare$}}}
        \put(137,67){\scriptsize$\blacksquare$}

        \put(117,-3){$*$}
        \put(157,-3){$*$}
        \put(115,58){$**$}
        \put(155,58){$**$}

        \put(90,15){$\diamond$}
        \put(185,15){$\diamond$}


        \end{normalsize}
\end{picture}
 }
  \caption{A coloured graph and the corresponding coloured factor graph.}
  \label{fig:factor_graph}
\end{figure}

For $Y = (Y_{\alpha})_{\alpha \in V}$ assumed to follow a
$\mathcal{N}_{|V|}(\mu,\Sigma)$ distribution, the density $f(y)$
factorises as
\begin{equation*}
f(y) \propto \prod_{\alpha \in V}
\exp\{-k_{\alpha\alpha}(y_{\alpha}-\mu_{\alpha})^{2}/2\} \ \cdot
\prod_{\substack{\alpha,\beta \in V, \\\alpha \not= \beta}}
\exp\{-k_{\alpha\beta}(y_{\alpha}-\mu_{\alpha})(y_{\beta}-\mu_{\beta})\}
\end{equation*}
giving that for the Gaussian distribution, either $A_{i}=\{\alpha\}$
or $A_{i}=\{\alpha, \beta\}$ for $\alpha, \beta \in V$, with a
factor being present if and only if the corresponding entry in
$\Sigma^{-1}$ is non-zero. Thus if the distribution of $Y$ is
assumed lie in the graphical Gaussian model represented by graph
$G=(V,E)$, by equation~(\ref{eq_cip}), each factor corresponds to a
vertex in $V$ or edge in $E$. The vertices in $V$ can clearly be
identified with their factors so that the factor graph of a
graphical Gaussian model with graph $G=(V,E)$ equals
\begin{equation*}
G_{F}=(V \cup F, E_{F}) \ \ \ \text{with} \ \ \  F=\{e \mid e \in
E\}  \ \ \ \text{and} \ \ \ E_{F} = \{\alpha e \mid \alpha \in V, e
\in E \textup{ is incident with $\alpha$ in } G\}.
\end{equation*}
This can be extended to the notion of a coloured factor graph, with
an example given in Figure \ref{fig:factor_graph}.

\begin{definition}
For $\mathcal{G}=(\mathcal{V},\mathcal{E}) \in \mathcal{C}_{V}$
representing a graphical Gaussian model with equality constraints,
let $N$ be a set of nodes with each node representing an edge in $E$
and let $\mathcal{N}$ be the colouring of $N$ in which nodes receive
the same colour if and only if the corresponding edges are equally
coloured in $\mathcal{E}$. The \emph{coloured factor graph} of the
model is defined to be the vertex and node coloured graph
$\mathcal{G}_{F} = (\mathcal{V} \cup \mathcal{N}, E_{F})$ with
$E_{F} = \{\alpha n \mid \alpha \in V, n \in N \mbox{\textup{ and
$n$ represents an edge incident with }} \alpha \}$. The set of
coloured factor graphs with vertex set $V$ is denoted
$\mathcal{F}_{V}$.
\end{definition}

We give four intermediate results whose proofs we omit for brevity.

\begin{lemma} \label{lemma_iso}
$\mathcal{F}_{V}$ and $\mathcal{C}_{V}$ are isomorphic. We denote
the isomorphism by $\phi_{V}: \mathcal{C}_{V} \rightarrow
\mathcal{F}_{V}$.
\end{lemma}

\begin{lemma} \label{lemma_join}
If $\mathcal{G}=(\mathcal{V},\mathcal{E}) \in \mathcal{C}_{V}$ and
$\phi_{V}(\mathcal{G})=\mathcal{G}_{F}=(\mathcal{V} \cup
\mathcal{N},E_{F})$ is the corresponding factor graph, then
$(\mathcal{V},\mathcal{E}) \in R_{V}$ if and only if $(\mathcal{V}
\cup \mathcal{N})$ is equitable with respect to $G_{F} = (V \cup N,
E_{F})$.
\end{lemma}

\begin{lemma}[\cite{mckay}] \label{lemma_mckay1}
If $P_{1}$ and $P_{2}$ are two partitions of $V$ both equitable with
respect to the same graph $G=(V,E)$, then so is their join $P_{1}
\lor P_{2}$.
\end{lemma}

\begin{lemma}[\cite{mckay}] \label{lemma_mckay2}
For $G=(V,E)$ and a partition $P$ of $V$, (up to the order of cells)
there exists a unique coarsest partition that is finer than $P$ and
equitable with respect to $G$, to be denoted by $r_{G}(P)$.
\end{lemma}

Combined, Lemmas \ref{lemma_iso}, \ref{lemma_join} and
\ref{lemma_mckay1} can be used to prove:

\begin{proposition} \label{prop_r}
$R_{V}$ is stable under the meet operation $\land$ in $\langle
\mathcal{C}_{V}; \preceq \rangle$ given in
equation~(\ref{meet_join_c}).
\end{proposition}

Note that the graphs in Figure \ref{fig:join_col_ex} illustrate the
stability of $R_{V}$ under $\wedge$ while showing that $R_{V}$ is
generally not stable under $\vee$ in $\langle \mathcal{C}_{V};
\preceq \rangle$. Further, Lemma \ref{lemma_mckay2} implies:

\begin{proposition} \label{prop_sup_r}
Let $\mathcal{G} =(\mathcal{V},\mathcal{E})\in \mathcal{C}_{V}$ and
let $\phi_{V}(\mathcal{G}) = \mathcal{G}_{F} = (\mathcal{V} \cup
\mathcal{N}, E_{F})$ be the corresponding coloured factor graph.
Then $\mathcal{G}$ has a supremum in $R_{V}$ given by
\begin{equation*}
s_{R}(\mathcal{G}) = \phi_{V}^{-1}((r_{G_{F}}(\mathcal{V} \cup
\mathcal{N}), E_{F}))
\end{equation*}
\end{proposition}

We conclude:

\begin{theorem} \label{thm_r}
$\mathcal{S}^{+}_{R_{V}}$ is a complete non-distributive lattice
with respect to model inclusion. The meet operation is induced by
the meet operation in $\langle \mathcal{C}_{V}; \leq \rangle$ given
in equation~(\ref{meet_join_c}). The join of two models represented
by graphs $\mathcal{G},\mathcal{H} \in R_{V}$ is represented by the
graph $s_{R}(\mathcal{G} \lor \mathcal{H})$.
\end{theorem}

\noindent \emph{Proof}: The proof is analogous to the proofs of
Theorems \ref{thm_b} and \ref{thm_pi}.
\whitebox

Lemma \ref{lemma_mckay2} further implies:

\begin{proposition} \label{prop_p}
$P_{V}$ is stable under the meet operation $\land$ in $\langle
\mathcal{C}_{V}; \preceq \rangle$ given in
equation~(\ref{meet_join_c}).
\end{proposition}

The graphs in Figure \ref{fig:join_col_ex} also illustrate the
stability of $P_{V}$ under $\wedge$ while establishing that $P_{V}$
is generally not stable under $\vee$ in $\langle \mathcal{C}_{V};
\preceq \rangle$. Further:

\begin{lemma} \label{lem_vr}
If $\mathcal{G} =(\mathcal{V},\mathcal{E})\in \mathcal{C}_{V}$ and
$s_{R}(\mathcal{G}) = (\mathcal{V}_{R},\mathcal{E}_{R})$, then for
all $\mathcal{G}'=(\mathcal{V}',\mathcal{E})\in \mathcal{C}_{V}$
with $\mathcal{V}_{R} \leq \mathcal{V}' \leq \mathcal{V}$,
$s_{R}(\mathcal{G}') = (\mathcal{V}_{R},\mathcal{E}_{R})$.
\end{lemma}

Lemma \ref{lem_vr} can be shown to imply:

\begin{proposition} \label{prop_sup_p}
For $\mathcal{G} =(\mathcal{V},\mathcal{E})\in \mathcal{C}_{V}$ let
$s_{R}(\mathcal{G}) = (\mathcal{V}_{R},\mathcal{E}_{R})$. Then
$\mathcal{G}$ has a supremum in $P_{V}$ given by
\begin{equation*}
s_{P}(\mathcal{G}) = (\mathcal{V}_{R},\mathcal{E})
\end{equation*}
\end{proposition}

We conclude:

\begin{theorem}
$\mathcal{S}^{+}_{P_{V}}$ and $\mathcal{R}^{+}_{P_{V}}$ are complete
lattices with respect to model inclusion. Their meet operation is
induced by the meet operation in $\langle \mathcal{C}_{V}; \leq
\rangle$ given in equation~(\ref{meet_join_c}). The join of two
models represented by $\mathcal{G},\mathcal{H} \in P_{V}$ is
represented by the graph $s_{P}(\mathcal{G} \lor \mathcal{H})$.
\end{theorem}

\noindent \emph{Proof}: In complete analogy to the proofs of
Theorems \ref{thm_b} and \ref{thm_pi}. \whitebox


\section{Model Selection}

One way to develop model selection procedures for the model classes
considered in this article is by adapting existing model search
algorithms for unconstrained graphical Gaussian models. Having shown
each of the model classes to be complete lattices, just as the set
of standard graphical models, it is natural to consider methods
which exploit this structural property. Prominent methods among them
are stepwise procedures \citep{whittaker, edwards}, the
Edwards--Havr{\'a}nek model selection
procedure \citep{edwards_lattice}, and, more recently, neighbourhood
selection with the lasso \citep{meinshausen2006}, stability
selection \citep{meinshausen2010}, and the SINful
approach \citep{drton2008}.

A crucial difference between the search spaces of unconstrained
graphical Gaussian models and the models studied here is that while
the former constitute a distributive lattice, the latter are all non-distributive. This directly disqualifies
neighbourhood selection with the lasso, stability selection and the
SINful approach, as they all require distributivity. Put explicitly,
while the just mentioned methods are algorithms for determining for
each edge whether it is to be present in graph of the accepted
model(s) or not, for the model classes considered here not only the
edge set needs to be determined, but also partitions of the vertices
and present edges into sets corresponding to equal model parameters.
This turns model selection into a principally different problem.

For the rest of the article we focus on the Edwards--Havr{\'a}nek
model selection procedure and develop a corresponding algorithm for
the lattice of models $\mathcal{S}^{+}_{B_{V}}$ represented by edge
regular colourings. We illustrate the algorithm with a brief summary
of its, very encouraging, performance for the data set described in
Example \ref{rcon_example}. We further give a summary of an
Edwards--Havr{\'a}nek model search within the lattice of models
$\mathcal{S}^{+}_{\Pi_{V}}$ with permutation-generated colourings
for the Fret's heads data described in Example \ref{rcop_example}.
All formal results are given without proof, however they can be
obtained by considering the partial ordering of $\mathcal{C}_{V}$.

\subsection{The Edwards--Havr{\'a}nek Model Selection Procedure}

The Edwards--Havr{\'a}nek model selection procedure operates on
model search spaces which are lattices and is closely related to the
all possible models approach but considerably faster. It is based on
the following two principles: (i) if a model is accepted then all
models that include it are (weakly) accepted, and (ii) if a model is
rejected then all of its submodels are considered to be (weakly)
rejected.

The procedure starts by initially testing a set of models and
assigns the accepted models to a set $\mathcal{A}$ and the rejected
models to set $\mathcal{R}$. By assumption, all models larger than
$\mathcal{A}$ are (weakly) accepted and the ones smaller than
$\mathcal{R}$ are (weakly) rejected, so that only $\min
\mathcal{A}$, the smallest models in $\mathcal{A}$, and $\max
\mathcal{R}$, the largest models in $\mathcal{R}$, are of interest.
The procedure repeatedly updates $\min \mathcal{A}$ and $\max
\mathcal{R}$ and terminates once the set to be updated remains
unchanged, when it returns $\min \mathcal{A}$. The method to
determine whether a model is to be rejected can be any suitable
statistical test in accordance to the principle of
\emph{coherence} \citep{gabriel}, stating that a test should not
accept a model while rejecting a larger one.

Following \citet{edwards_lattice}, for a set of models $S$ let
$D_{a}(S)$ denote the set of models in the search lattice $L$ say
which are smallest with the property that they are not contained in
any model in $S$,
\begin{equation*}
D_{a}(S) = \min \{d \in L \mid d \not\subseteq s \textup{  for all }
s\in S\}
\end{equation*}
and let $D_{r}(S)$ be the set of largest models that do not contain
any model in $S$,
\begin{equation*}
D_{r}(S) = \max \{d \in L \mid s \not\subseteq d \textup{  for all }
s\in S\}
\end{equation*}
$D_{a}(S)$ is referred to as the \emph{acceptance dual} of $S$, and
$D_{r}(S)$ as the \emph{rejection dual} of $S$. The procedure may
then be summarised as:
\begin{itemize}
\item[1.] Test an initial set of models and assign the accepted
models to $\mathcal{A}$ and the rejected models to $\mathcal{R}$.
\item[2.] Choose between 3 and 4.
\item[3.] Test the models in $D_{r}(\mathcal{A})\setminus \mathcal{R}$.
If all are rejected, stop; otherwise, update $\mathcal{A}$ and $\mathcal{R}$ and go to 2.
\item[4.] Test the models in $D_{a}(\mathcal{R})\setminus \mathcal{A}$.
If all are accepted, stop; otherwise, update $\mathcal{A}$ and $\mathcal{R}$ and go to 2.
\end{itemize}


Acceptance and rejection duals of sets of models can be computed in
a recursive manner by using the following two relations. If $S$ and
$T$ are two sets of models, then
\begin{eqnarray*}
D_{a}(S \cup T) &=& \min\{s \vee t \mid s \in D_{a}(S), t \in D_{a}(T)\} \label{eq_dual1}\\
D_{r}(S \cup T) &=& \max\{s \wedge t \mid s \in D_{r}(S), t \in
D_{r}(T)\} \label{eq_dual2}
\end{eqnarray*}
Thus describing the duals of a single model is enough.

\subsection{Models Represented by Edge Regular Colourings}

\begin{proposition}
Let $\mathcal{S}^{+}(\mathcal{V}, \mathcal{E}) \in
\mathcal{S}^{+}_{B_{V}}$ be a model represented by graph
$\mathcal{G}=(\mathcal{V}, \mathcal{E})$ with edge regular colouring
and underlying uncoloured graph $G=(V,E)$. Then the acceptance dual
$D_{a}(\mathcal{S}^{+}(\mathcal{V}, \mathcal{E}))$ of
$\mathcal{S}^{+}(\mathcal{V}, \mathcal{E})$ in
$\mathcal{S}^{+}_{B_{V}}$ contains all models represented by
coloured graphs $\mathcal{G}_{a}=(\mathcal{V}_{a}, \mathcal{E}_{a})$
satisfying
\begin{itemize}
\item[(1i)] $\mathcal{V}_{a} = \{V_{1}, V_{2}\}$ such that $\mathcal{V}_{a} \not\geq \mathcal{V}$ and $\mathcal{E}_{a} = \emptyset$
\item[(1ii)] $\mathcal{V}_{a} = \{V\}$ and $\mathcal{E}_{a} = \{E_{a}\}$ with $\mathcal{E}_{a} \not= \emptyset$ and $\mathcal{E}_{a} \not\geq \mathcal{E}$.
\end{itemize}
\end{proposition}

\begin{figure}[ht!]
  \centering
  \subfigure[A graph of type (1i).]{\label{fig_adual_i}
  \begin{picture}(144,50)(25,-5) 
\begin{normalsize}
\thicklines
        \put(66,2){1}
        \put(66,43){4}
        \put(130,2){2}
        \put(130,43){3}
        \put(80,5){\color{blue}{\circle*{6}}}
        \put(80,45){\color{blue}{\circle*{6}}}
        \put(120,5){\color{red}{\circle*{6}}}
        \put(120,45){\color{red}{\circle*{6}}}
        \put(77,-8){$*$}
        \put(114,-8){$**$}
        \put(77,51){$*$}
        \put(114,51){$**$}
        \end{normalsize}
        \end{picture}
  }
  \subfigure[A graph of type (1ii).]{\label{fig_adual_ii}
\hspace{-10pt}
\begin{picture}(144,50)(25,-5) 
\begin{normalsize}
\thicklines
        \put(66,2){1}
        \put(66,43){4}
        \put(130,2){2}
        \put(130,43){3}
        \put(80,5){\color{darkgreen}{\line(1,0){40}}}
        \put(80,5){\color{darkgreen}{\line(0,1){40}}}
        \put(120,45){\color{darkgreen}{\line(-1,-1){40}}}
        \put(120,45){\color{darkgreen}{\line(0,-1){40}}}
        \put(80,5){\color{blue}{\circle*{6}}}
        \put(80,45){\color{blue}{\circle*{6}}}
        \put(120,5){\color{blue}{\circle*{6}}}
        \put(120,45){\color{blue}{\circle*{6}}}
        \put(77,-8){$*$}
        \put(117,-8){$*$}
        \put(77,51){$*$}
        \put(117,51){$*$}
                \put(100,2){\line(0,1){6}}
                \put(100,22){\line(0,1){6}}
                \put(77,25){\line(1,0){6}}
                \put(117,25){\line(1,0){6}}
        \end{normalsize}
\end{picture}
  }
  \caption{Acceptance dual corresponding to the graph in Figure \ref{fig:edge_reg}.}
  \label{fig_adual}
\end{figure}
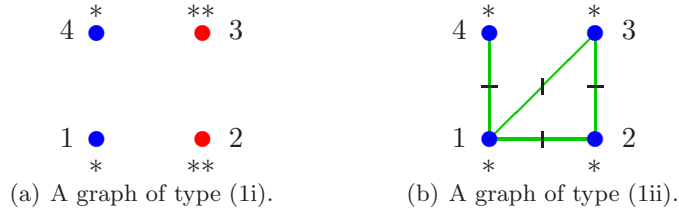

Put into words, models in the acceptance dual of
$\mathcal{S}^{+}(\mathcal{V}, \mathcal{E}) \in
\mathcal{S}^{+}_{B_{V}}$ are either represented by the empty graph
with two vertex colour classes which are not unions of colour
classes in $\mathcal{V}$, or they are represented by graphs in which
all vertices are of the same colour, as are the edges, and the edge
set is not a union of colour classes in $\mathcal{E}$. For example,
the coloured graphs in Figure \ref{fig_adual} display models which
lie in the acceptance dual of the model represented by the edge
regular colouring in Figure \ref{fig:edge_reg}.

By definition, acceptance duals are used to test whether there exist
models immediately larger than $\max \mathcal{R}$ which can be
rejected. Effectively, graphs of type (1i) refine the colouring of
the maximally rejected models, while graphs of type (1ii) add edge
colour classes, which both give larger models.

\begin{proposition}
Let $\mathcal{S}^{+}(\mathcal{V}, \mathcal{E}) \in
\mathcal{S}^{+}_{B_{V}}$ be a model represented by graph
$\mathcal{G}=(\mathcal{V}, \mathcal{E})$ with edge regular colouring
and underlying uncoloured graph $G=(V,E)$, and for a discrete set
$A$, let $\textup{atom}(A)$ denote the partition of $A$ into atomic
sets. Then the rejection dual $D_{r}(\mathcal{S}^{+}(\mathcal{V},
\mathcal{E}))$ of $\mathcal{S}^{+}(\mathcal{V}, \mathcal{E})$ in
$\mathcal{S}^{+}_{B_{V}}$ contains all models represented by
coloured graphs $\mathcal{G}_{r}=(\mathcal{V}_{r}, \mathcal{E}_{r})$
satisfying
\begin{itemize}
\item[(2i)] $\mathcal{V}_{r} = \{\{\alpha, \beta\}\} \cup \textup{atom}(V \setminus \{\alpha, \beta\})$ such that $\mathcal{V}_{r} \not\leq \mathcal{V}$ and $\mathcal{E}_{r} = \{\alpha\beta \mid \alpha, \beta \in V\}$, or
\item[(2ii)] $\mathcal{V}_{r} = \textup{atom}(V)$ and $\mathcal{E}_{r} = \textup{atom}(\{\alpha\beta \mid \alpha, \beta \in V\} \setminus \{e\})$ with $e \in E$, or
\item[(2iii)] $\mathcal{V}_{r} = \{\{\alpha, \beta\},\{\gamma, \delta\}\} \cup \textup{atom}(V \setminus \{\alpha, \beta, \gamma, \delta\})$ with $\alpha, \beta$ and $\gamma, \delta$ being of the same colour in $\mathcal{V}$ and $\mathcal{E}_{r} = \{\{\alpha\gamma, \beta\delta\}\} \cup \textup{atom}(\{\alpha\beta \mid \alpha, \beta \in V\} \setminus \{\alpha\gamma, \beta\delta\})$, where we may have $\alpha = \beta$ or $\gamma=\delta$ but not both, such that $(\mathcal{V},\mathcal{E}) \not\preceq (\mathcal{V}_{r},\mathcal{E}_{r})$.
\end{itemize}
\end{proposition}

\begin{figure}[ht!]
  \centering
  \subfigure[A graph of type (2i)]{\label{fig_rdual_i}
  \begin{picture}(144,60)(25,-5) 
\begin{normalsize}
\thicklines
        \put(66,2){1}
        \put(66,43){4}
        \put(127,2){2}
        \put(127,43){3}
        \put(80,5){{\line(1,0){40}}}
        \put(80,5){{\line(0,1){40}}}
        \put(80,5){{\line(1,1){40}}}
        \put(80,45){{\line(1,-1){40}}}
        \put(120,45){{\line(-1,0){40}}}
        \put(120,45){{\line(0,-1){40}}}
        \put(80,5){\color{blue}{\circle*{6}}}
        \put(80,45){\color{blue}{\circle*{6}}}
        \put(120,5){{\circle*{6}}}
        \put(120,45){{\circle*{6}}}
        \put(77,-8){$*$}
        \put(77,51){$*$}
        \end{normalsize}
        \end{picture}
} \hspace{-35pt}
  \subfigure[A graph of type (2ii).]{\label{fig_rdual_ii}
\begin{picture}(144,60)(25,-5) 
\begin{normalsize}
\thicklines
        \put(66,2){1}
        \put(66,43){4}
        \put(127,2){2}
        \put(127,43){3}
        \put(80,5){{\line(1,0){40}}}
        \put(80,5){{\line(0,1){40}}}
        \put(120,45){{\line(-1,-1){40}}}
        \put(80,45){{\line(1,-1){40}}}
        \put(120,45){{\line(0,-1){40}}}
        \put(80,5){{\circle*{6}}}
        \put(80,45){{\circle*{6}}}
        \put(120,5){{\circle*{6}}}
        \put(120,45){{\circle*{6}}}
     \end{normalsize}
\end{picture}
  } \hspace{-35pt}
  \subfigure[A graph of type (2iii).]{\label{fig_rdual_iii}
\begin{picture}(144,60)(25,-5) 
\begin{normalsize}
\thicklines
        \put(66,2){1}
        \put(66,43){4}
        \put(127,2){2}
        \put(127,43){3}
        \put(80,5){\color{darkgreen}{\line(1,0){40}}}
        \put(80,5){{\line(0,1){40}}}
        \put(120,45){{\line(-1,-1){40}}}
        \put(80,45){{\line(1,-1){40}}}
        \put(120,45){\color{darkgreen}{\line(0,-1){40}}}
        \put(120,45){{\line(-1,0){40}}}
        \put(80,5){\color{blue}{\circle*{6}}}
        \put(80,45){{\circle*{6}}}
        \put(120,5){{\circle*{6}}}
        \put(120,45){\color{blue}{\circle*{6}}}
        \put(77,-8){$*$}
        \put(117,51){$*$}
        \put(100,2){{\line(0,1){6}}}
        \put(117,25){{\line(1,0){6}}}
                \end{normalsize}
\end{picture}
  }
  \caption{Rejection dual corresponding to the graph in Figure \ref{fig:edge_reg}.}
  \label{fig_r_dual}
\end{figure}
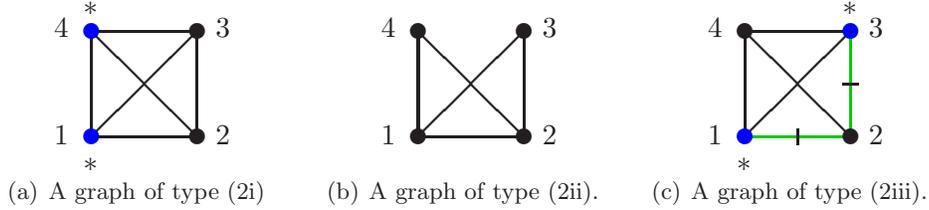

Graphs representing models in the rejection dual of a model
represented by $\mathcal{G}=(\mathcal{V},\mathcal{E})$ almost
represent the unrestricted saturated model except for a minor
modification: In graphs of type (2i) two vertices which are not of
the same colour in $\mathcal{V}$ form the only composite colour
class in $\mathcal{V}_{r}$, while graphs of~type (2ii) are missing
an edge present in $\mathcal{G}=(\mathcal{V},\mathcal{E})$. Graphs
of type (2iii) have a pair of equally coloured edges which are not
of the same colour in $\mathcal{E}$, and give the end vertices of
the edges an appropriate colouring for the graph colouring to be
edge regular. Examples of coloured graphs representing models in the
rejection dual of the model represented by the graph displayed in
Figure \ref{fig:edge_reg} are given in Figure \ref{fig_r_dual}.

Rejection duals contain models which lie immediately below $\min
\mathcal{A}$. Graphs of type (2i) merge vertex colour classes,
graphs of type (2ii) cause edge colour classes to be dropped and
graphs of type (2iii) merge edge colour classes, as well as vertex
colour classes to ensure the resulting model to be edge regular. All
operations give graphs which represent smaller models.

It can be shown that while
$|D_{a}(\mathcal{S}^{+}(\mathcal{V},\mathcal{E}))|$ grows
super-exponentially in the number of variables $|V|$ (at rate
$O(2^{|V|^{2}/2})$), the size of the rejection dual
$|D_{r}(\mathcal{S}^{+}(\mathcal{V},\mathcal{E}))|$ grows
polynomially in $|V|$ (at rate $O(|V|^{4})$), so that from a
computational point of view working with rejection duals only is
much more efficient. The following algorithm is therefore most
efficient for models with edge regular colourings.

\begin{itemize}
\item[1.] Test an initial set of models and assign each to $\mathcal{R}$ if it is rejected and to $\mathcal{A}$
otherwise.
\item[2.] Test the models in $D_{r}(\mathcal{A}) \setminus \mathcal{R}$. If all are rejected, stop. Otherwise update $\mathcal{A}$ and $\mathcal{R}$ and repeat.
\end{itemize}

\begin{figure}[ht!]
  \centering
  \subfigure{
\begin{picture}(180,190) (10,-10)
        \begin{normalsize}
        \thicklines

        \put(95,170){$\mathcal{G}_{1}$}
        \put(40,100){\color{darkgreen}{\line(0,1){60}}}
        \put(160,100){\line(0,1){60}}
        \put(40,100){\color{purple}{\line(2,1){60}}}
        \put(40,160){\line(2,-1){60}}
        \put(160,100){\color{purple}{\line(-2,1){60}}}
        \put(160,160){\line(-2,-1){60}}
        {\color{darkgreen}{\curve(40,100, 100,115, 160,160)}}

        \put(28,98){$*$}
                \put(167,98){$*$}
                \put(22,158){$**$}
                \put(167,158){$**$}

        \put(40,100){\color{blue}{\circle*{6}}}
        \put(40,160){\color{red}{\circle*{6}}}
        \put(160,100){\color{blue}{\circle*{6}}}
        \put(100,130){\circle*{6}}
        \put(160,160){\color{red}{\circle*{6}}}

        \put(20,80){\footnotesize Mechanics}
        \put(140,80){\footnotesize Statistics}
        \put(20,170){\footnotesize Vectors}
        \put(140,170){\footnotesize Analysis}
        \put(85,143){\footnotesize Algebra}

        \put(37,130){\line(1,0){6}}
                \put(71,112){\line(0,1){6}}
                \put(69,112){\line(0,1){6}}
                \put(129,112){\line(0,1){6}}
                \put(131,112){\line(0,1){6}}
                \put(100,112){\line(0,1){6}}

        \put(75,80){\footnotesize BIC 2601.617}

        \put(95,50){$\mathcal{G}_{3}$}
        \put(40,-20){{\line(0,1){60}}}
        \put(160,-20){\color{darkgreen}{\line(0,1){60}}}
        \put(40,-20){\color{purple}{\line(2,1){60}}}
        \put(40,40){\line(2,-1){60}}
        \put(160,-20){\color{purple}{\line(-2,1){60}}}
        \put(160,40){\line(-2,-1){60}}

        {\color{darkgreen}{\curve(40,-20, 100,-5, 160,40)}}

        \put(40,-20){\color{blue}{\circle*{6}}}
        \put(40,40){\color{red}{\circle*{6}}}
        \put(160,-20){\color{blue}{\circle*{6}}}
        \put(100,10){\circle*{6}}
        \put(160,40){\color{red}{\circle*{6}}}

        \put(28,-22){$*$}
                \put(167,-22){$*$}
                \put(22,38){$**$}
                \put(167,38){$**$}

        \put(20,-40){\footnotesize Mechanics}
        \put(140,-40){\footnotesize Statistics}
        \put(20,50){\footnotesize Vectors}
        \put(140,50){\footnotesize Analysis}
        \put(85,23){\footnotesize Algebra}

        \put(157,10){\line(1,0){6}}
                \put(71,-8){\line(0,1){6}}
                \put(69,-8){\line(0,1){6}}
                \put(129,-8){\line(0,1){6}}
                \put(131,-8){\line(0,1){6}}
                \put(100,-8){\line(0,1){6}}

        \put(75,-40){\footnotesize BIC 2603.376}

                \end{normalsize}
\end{picture}
  }
\hspace{6pt}
  \subfigure{
\begin{picture}(180,190) (10,-10)
\begin{normalsize}
\thicklines
        \put(95,170){$\mathcal{G}_{2}$}
        \put(40,100){\color{darkgreen}{\line(0,1){60}}}
        \put(160,100){\color{darkgreen}{\line(0,1){60}}}
        \put(40,100){{\line(2,1){60}}}
        \put(40,160){{\line(2,-1){60}}}
        \put(160,100){{\line(-2,1){60}}}
        \put(160,160){\line(-2,-1){60}}
        {\color{darkgreen}{\curve(40,100, 100,115, 160,160)}}

        \put(40,100){\color{blue}{\circle*{6}}}
        \put(40,160){\color{red}{\circle*{6}}}
        \put(160,100){\color{blue}{\circle*{6}}}
        \put(100,130){\circle*{6}}
        \put(160,160){\color{red}{\circle*{6}}}

                \put(37,130){\line(1,0){6}}
                \put(157,130){\line(1,0){6}}
                \put(100,112){\line(0,1){6}}

        \put(75,80){\footnotesize BIC 2600.017}

              \put(20,80){\footnotesize Mechanics}
        \put(140,80){\footnotesize Statistics}
        \put(20,170){\footnotesize Vectors}
        \put(140,170){\footnotesize Analysis}
        \put(85,143){\footnotesize Algebra}

                \put(28,98){$*$}
                \put(167,98){$*$}
                \put(22,158){$**$}
                \put(167,158){$**$}

        \put(95,50){$\mathcal{G}_{4}$}
        \put(40,-20){\color{darkgreen}{\line(0,1){60}}}
        \put(160,-20){\color{darkgreen}{\line(0,1){60}}}
        \put(40,-20){\color{purple}{\line(2,1){60}}}
        \put(40,40){\line(2,-1){60}}
        \put(160,-20){\color{purple}{\line(-2,1){60}}}
        \put(160,40){\line(-2,-1){60}}

        \put(28,-22){$*$}
                \put(167,-22){$*$}
                \put(22,38){$**$}
                \put(167,38){$**$}

        \put(40,-20){\color{blue}{\circle*{6}}}
        \put(40,40){\color{red}{\circle*{6}}}
        \put(160,-20){\color{blue}{\circle*{6}}}
        \put(100,10){\circle*{6}}
        \put(160,40){\color{red}{\circle*{6}}}

        \put(20,-40){\footnotesize Mechanics}
        \put(140,-40){\footnotesize Statistics}
        \put(20,50){\footnotesize Vectors}
        \put(140,50){\footnotesize Analysis}
        \put(85,23){\footnotesize Algebra}

        \put(37,10){\line(1,0){6}}
        \put(157,10){\line(1,0){6}}
                \put(71,-8){\line(0,1){6}}
                \put(69,-8){\line(0,1){6}}
                \put(129,-8){\line(0,1){6}}
                \put(131,-8){\line(0,1){6}}

        \put(75,-40){\footnotesize BIC 2591.468}

        \end{normalsize}
\end{picture}
  }
\vspace{25pt}
  \caption{Graphs of minimally accepted models in $\mathcal{S}^{+}_{B_{[5]}}$.}
  \label{fig_edge_reg_search}
\end{figure}
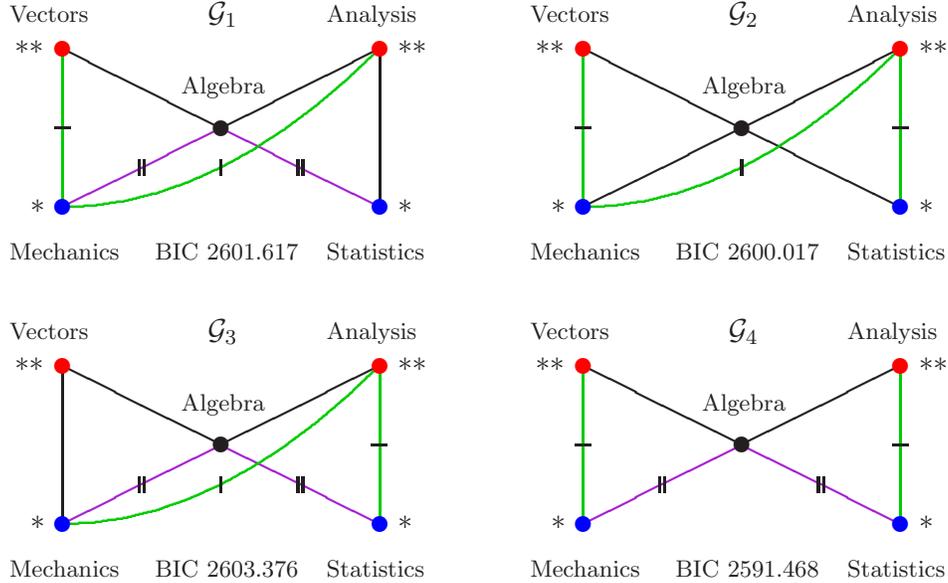

We executed the above algorithm for the Mathematics marks data set
described in Example \ref{rcon_example} with the saturated
uncoloured model as our initial set of accepted models $\mathcal{A}$
and let $\mathcal{R} = \emptyset$ initially. Models were tested for
acceptance by performing a likelihood ratio test relative to the
saturated unconstrained model at significance level $5\%$ using
functionality implemented in the \texttt{R} package
\texttt{gRc}~\cite{grc}. The algorithm fitted 232 models, out of a
total of $1.3 \cdot 10^{6}$, in 8 stages before arriving at 4
minimally accepted models whose graphs are displayed in Figure
\ref{fig_edge_reg_search} together with their BIC values. (The 232
models are distributed among the stages as follows. 1: 20 (6
accepted), 2: 21 (19 accepted), 3: 41 (40 accepted), 4: 56 (56
accepted), 5: 55 (55 accepted), 6: 29 (29 accepted), 7: 9 (9
accepted) and 8: 1 (1 accepted).)

The uncoloured graphs underlying the graphs in Figure
\ref{fig_edge_reg_search} contain the graph in Figure
\ref{fig:mathmark_uncol} as a subgraph. $\mathcal{G}_{1}$,
$\mathcal{G}_{2}$ and $\mathcal{G}_{3}$ only differ by one edge and
$\mathcal{G}_{4}$ has exactly the same edge set. Thus the
conditional independence structures largely agree. The minimally
accepted model with the lowest BIC value is represented by
$\mathcal{G}_{4}$. It is different from but not dissimilar to the
RCON model fitted in \citet{lauritzen_sym}, whose graph is displayed
in Figure \ref{fig:mathmark_col}, which has a slightly lower BIC
value of 2587.404 but no specific properties and it is chosen in an
ad hoc manner. Note that the model fitted in \citet{lauritzen_sym} is
not edge regular so that it could not have been considered by the
algorithm.

This example suggests that an Edwards--Havr{\'a}nek model selection
procedure for models with edge regular colourings may be feasible in
general.

\subsection{Models Represented by Permutation-Generated Colourings}

The class of permutation-generated colourings $\Pi_{V}$ is more
complex in its structure than $B_{V}$ and therefore the duals
$D_{a}(\mathcal{S}^{+}(\mathcal{V},\mathcal{E}))$ and
$D_{r}(\mathcal{S}^{+}(\mathcal{V},\mathcal{E}))$ of a model
$\mathcal{S}^{+}(\mathcal{V},\mathcal{E})$ cannot be given in a
purely combinatorial form in the graph colouring. A sound
understanding of the relationship between a group $\Gamma$ and its
orbits in $V$ and $V \times V$ is required in order to design a
general algorithm, in principle, applicable to any variable set $V$.
For illustrative purposes of the underlying principles we provide a
brief summary of an Edwards--Havr{\'a}nek model search in
$\mathcal{S}^{+}_{\Pi_{V}}$ for Fret's heads data described in
Example \ref{rcop_example}.

\begin{figure}[ht!]
  \centering
  \setlength{\unitlength}{.97pt}
\begin{picture}(450,180)(30,-30)
\begin{normalsize}
\thicklines
         \put(6,90){1}
        \put(6,130){4}
        \put(66,90){2}
        \put(66,130){3}
        \put(20,90){{\line(0,1){40}}}
        \put(60,130){{\line(0,-1){40}}}
        \put(20,90){{\line(1,0){40}}}
        \put(60,130){{\line(-1,0){40}}}
        \put(20,90){{\line(1,1){40}}}
        \put(60,90){{\line(-1,1){40}}}
        \put(20,90){{\circle*{6}}}
        \put(20,130){{\circle*{6}}}
        \put(60,90){{\circle*{6}}}
        \put(60,130){{\circle*{6}}}

        \put(37,140){\footnotesize $\mathcal{G}_{1}$}
        \put(18,68){\footnotesize $\Gamma_{1}=\langle Id \rangle$}
        \put(37,50){\footnotesize $\mathcal{G}_{7}$}

                \put(91,90){1}
        \put(91,130){4}
        \put(151,90){2}
        \put(151,130){3}
        \put(105,90){\color{darkgreen}{\line(0,1){40}}}
        \put(145,130){\color{purple}{\line(0,-1){40}}}
        \put(105,90){{\line(1,0){40}}}
        \put(145,130){{\line(-1,0){40}}}
        \put(105,90){\color{purple}{\line(1,1){40}}}
        \put(145,90){\color{darkgreen}{\line(-1,1){40}}}
        \put(105,90){\color{blue}{\circle*{6}}}
        \put(105,130){{\circle*{6}}}
        \put(145,90){\color{blue}{\circle*{6}}}
        \put(145,130){{\circle*{6}}}
        \put(102,79){$*$}
        \put(142,79){$*$}
                \put(102,110){\line(1,0){6}}
                \put(142,111){\line(1,0){6}}
                \put(142,109){\line(1,0){6}}
                \put(114,97){\line(0,1){6}}
                \put(116,97){\line(0,1){6}}
                \put(135,97){\line(0,1){6}}

                \put(122,140){\footnotesize $\mathcal{G}_{2}$}
        \put(100,68){\footnotesize $\Gamma_{2}=\langle(12)\rangle$}
        \put(122,50){\footnotesize $\mathcal{G}_{8}$}

                \put(176,90){1}
        \put(176,130){4}
        \put(236,90){2}
        \put(236,130){3}
        \put(190,90){\color{darkgreen}{\line(0,1){40}}}
        \put(230,130){\color{purple}{\line(0,-1){40}}}
        \put(190,90){\color{purple}{\line(1,0){40}}}
        \put(230,130){\color{darkgreen}{\line(-1,0){40}}}
        \put(190,90){{\line(1,1){40}}}
        \put(230,90){{\line(-1,1){40}}}
        \put(190,90){\color{blue}{\circle*{6}}}
        \put(190,130){{\circle*{6}}}
        \put(230,90){{\circle*{6}}}
        \put(230,130){\color{blue}{\circle*{6}}}
        \put(187,79){$*$}
        \put(227,136){$*$}
                \put(187,110){\line(1,0){6}}
                \put(227,111){\line(1,0){6}}
                \put(227,109){\line(1,0){6}}
            \put(210,127){\line(0,1){6}}
            \put(211,87){\line(0,1){6}}
            \put(209,87){\line(0,1){6}}

            \put(207,140){\footnotesize $\mathcal{G}_{3}$}
        \put(185,68){\footnotesize $\Gamma_{3}=\langle(13)\rangle$}
        \put(207,50){\footnotesize $\mathcal{G}_{9}$}

                \put(261,90){1}
        \put(261,130){4}
        \put(321,90){2}
        \put(321,130){3}
        \put(275,90){{\line(0,1){40}}}
        \put(315,130){{\line(0,-1){40}}}
        \put(275,90){\color{purple}{\line(1,0){40}}}
        \put(315,130){\color{darkgreen}{\line(-1,0){40}}}
        \put(275,90){\color{darkgreen}{\line(1,1){40}}}
        \put(315,90){\color{purple}{\line(-1,1){40}}}
        \put(275,90){\color{blue}{\circle*{6}}}
        \put(275,130){\color{blue}{\circle*{6}}}
        \put(315,90){{\circle*{6}}}
        \put(315,130){{\circle*{6}}}
        \put(272,79){$*$}
        \put(272,136){$*$}
                \put(285,97){\line(0,1){6}}
                \put(306,97){\line(0,1){6}}
                \put(304,97){\line(0,1){6}}
            \put(295,127){\line(0,1){6}}
            \put(296,87){\line(0,1){6}}
            \put(294,87){\line(0,1){6}}

            \put(292,140){\footnotesize $\mathcal{G}_{4}$}
        \put(270,68){\footnotesize $\Gamma_{4}=\langle(14)\rangle$}
        \put(290,50){\footnotesize $\mathcal{G}_{10}$}

                \put(346,90){1}
        \put(346,130){4}
        \put(406,90){2}
        \put(406,130){3}
        \put(360,90){{\line(0,1){40}}}
        \put(400,130){{\line(0,-1){40}}}
        \put(360,90){\color{purple}{\line(1,0){40}}}
        \put(400,130){\color{darkgreen}{\line(-1,0){40}}}
        \put(360,90){\color{purple}{\line(1,1){40}}}
        \put(400,90){\color{darkgreen}{\line(-1,1){40}}}
        \put(360,90){{\circle*{6}}}
        \put(360,130){{\circle*{6}}}
        \put(400,90){\color{blue}{\circle*{6}}}
        \put(400,130){\color{blue}{\circle*{6}}}
        \put(397,79){$*$}
        \put(397,136){$*$}
                \put(369,97){\line(0,1){6}}
                \put(371,97){\line(0,1){6}}
                \put(390,97){\line(0,1){6}}
            \put(380,127){\line(0,1){6}}
            \put(379,87){\line(0,1){6}}
            \put(381,87){\line(0,1){6}}

            \put(377,140){\footnotesize $\mathcal{G}_{5}$}
        \put(355,68){\footnotesize $\Gamma_{5}=\langle(23)\rangle$}
        \put(375,50){\footnotesize $\mathcal{G}_{11}$}

                \put(431,90){1}
        \put(431,130){4}
        \put(491,90){2}
        \put(491,130){3}
        \put(445,90){\color{purple}{\line(0,1){40}}}
        \put(485,130){\color{darkgreen}{\line(0,-1){40}}}
        \put(445,90){\color{purple}{\line(1,0){40}}}
        \put(485,130){\color{darkgreen}{\line(-1,0){40}}}
        \put(445,90){{\line(1,1){40}}}
        \put(485,90){{\line(-1,1){40}}}
        \put(445,90){{\circle*{6}}}
        \put(445,130){\color{blue}{\circle*{6}}}
        \put(485,90){\color{blue}{\circle*{6}}}
        \put(485,130){{\circle*{6}}}
        \put(482,79){$*$}
        \put(442,136){$*$}
                \put(442,109){\line(1,0){6}}
                \put(442,111){\line(1,0){6}}
                \put(482,110){\line(1,0){6}}
            \put(465,127){\line(0,1){6}}
            \put(466,87){\line(0,1){6}}
            \put(464,87){\line(0,1){6}}

        \put(462,140){\footnotesize $\mathcal{G}_{6}$}
        \put(440,68){\footnotesize $\Gamma_{6}=\langle(24)\rangle$}
        \put(460,50){\footnotesize $\mathcal{G}_{12}$}

        \put(6,0){1}
        \put(6,40){4}
        \put(66,0){2}
        \put(66,40){3}
        \put(20,0){\color{darkgreen}{\line(0,1){40}}}
        \put(60,40){\color{purple}{\line(0,-1){40}}}
        \put(20,0){{\line(1,0){40}}}
        \put(60,40){{\line(-1,0){40}}}
        \put(20,0){\color{darkgreen}{\line(1,1){40}}}
        \put(60,0){\color{purple}{\line(-1,1){40}}}
        \put(20,0){{\circle*{6}}}
        \put(20,40){\color{blue}{\circle*{6}}}
        \put(60,0){{\circle*{6}}}
        \put(60,40){\color{blue}{\circle*{6}}}
        \put(17,46){$*$}
        \put(57,46){$*$}
                \put(17,20){\line(1,0){6}}
                \put(57,21){\line(1,0){6}}
                \put(57,19){\line(1,0){6}}
                \put(30,7){\line(0,1){6}}
                \put(51,7){\line(0,1){6}}
                \put(49,7){\line(0,1){6}}

        \put(15,-22){\footnotesize $\Gamma_{7}=\langle (34) \rangle$}
        \put(35,-40){\footnotesize $\mathcal{G}_{13}$}

                \put(91,0){1}
        \put(91,40){4}
        \put(151,0){2}
        \put(151,40){3}
        \put(105,0){\color{darkgreen}{\line(0,1){40}}}
        \put(145,40){\color{purple}{\line(0,-1){40}}}
        \put(105,0){\color{purple}{\line(1,0){40}}}
        \put(145,40){\color{darkgreen}{\line(-1,0){40}}}
        \put(105,0){\color{purple}{\line(1,1){40}}}
        \put(145,0){\color{darkgreen}{\line(-1,1){40}}}
        \put(105,0){\color{blue}{\circle*{6}}}
        \put(105,40){{\circle*{6}}}
        \put(145,0){\color{blue}{\circle*{6}}}
        \put(145,40){\color{blue}{\circle*{6}}}
        \put(102,-11){$*$}
        \put(142,-11){$*$}
        \put(142,46){$*$}
                \put(102,20){\line(1,0){6}}
                \put(142,21){\line(1,0){6}}
                \put(142,19){\line(1,0){6}}
                \put(114,7){\line(0,1){6}}
                \put(116,7){\line(0,1){6}}
                \put(135,7){\line(0,1){6}}
                \put(124,-3){\line(0,1){6}}
                \put(126,-3){\line(0,1){6}}
                \put(125,37){\line(0,1){6}}

                \put(97,-22){\footnotesize $\Gamma_{8}=\langle (123) \rangle$}
        \put(120,-40){\footnotesize $\mathcal{G}_{14}$}

                \put(176,0){1}
        \put(176,40){4}
        \put(236,0){2}
        \put(236,40){3}
        \put(190,0){\color{darkgreen}{\line(0,1){40}}}
        \put(230,40){\color{purple}{\line(0,-1){40}}}
        \put(190,0){\color{darkgreen}{\line(1,0){40}}}
        \put(230,40){\color{purple}{\line(-1,0){40}}}
        \put(190,0){\color{purple}{\line(1,1){40}}}
        \put(230,0){\color{darkgreen}{\line(-1,1){40}}}
        \put(190,0){\color{blue}{\circle*{6}}}
        \put(190,40){\color{blue}{\circle*{6}}}
        \put(230,0){\color{blue}{\circle*{6}}}
        \put(230,40){{\circle*{6}}}
        \put(187,-11){$*$}
        \put(227,-11){$*$}
        \put(187,46){$*$}
                \put(187,20){\line(1,0){6}}
                \put(227,21){\line(1,0){6}}
                \put(227,19){\line(1,0){6}}
                \put(199,7){\line(0,1){6}}
                \put(201,7){\line(0,1){6}}
                \put(220,7){\line(0,1){6}}
            \put(211,37){\line(0,1){6}}
            \put(209,37){\line(0,1){6}}
            \put(210,-3){\line(0,1){6}}

                \put(182,-22){\footnotesize $\Gamma_{9}=\langle (124) \rangle$}
        \put(205,-40){\footnotesize $\mathcal{G}_{15}$}

                \put(261,0){1}
        \put(261,40){4}
        \put(321,0){2}
        \put(321,40){3}
        \put(275,0){\color{darkgreen}{\line(0,1){40}}}
        \put(315,40){\color{purple}{\line(0,-1){40}}}
        \put(275,0){\color{purple}{\line(1,0){40}}}
        \put(315,40){\color{darkgreen}{\line(-1,0){40}}}
        \put(275,0){\color{darkgreen}{\line(1,1){40}}}
        \put(315,0){\color{purple}{\line(-1,1){40}}}
        \put(275,0){\color{blue}{\circle*{6}}}
        \put(275,40){\color{blue}{\circle*{6}}}
        \put(315,0){{\circle*{6}}}
        \put(315,40){\color{blue}{\circle*{6}}}
        \put(272,-11){$*$}
        \put(272,46){$*$}
        \put(312,46){$*$}
                \put(272,20){\line(1,0){6}}
                \put(312,21){\line(1,0){6}}
                \put(312,19){\line(1,0){6}}
                \put(285,7){\line(0,1){6}}
                \put(306,7){\line(0,1){6}}
                \put(304,7){\line(0,1){6}}
            \put(295,37){\line(0,1){6}}
            \put(296,-3){\line(0,1){6}}
            \put(294,-3){\line(0,1){6}}

                \put(265,-22){\footnotesize $\Gamma_{10}=\langle (134) \rangle$}
        \put(290,-40){\footnotesize $\mathcal{G}_{16}$}

                \put(346,0){1}
        \put(346,40){4}
        \put(406,0){2}
        \put(406,40){3}
        \put(360,0){\color{purple}{\line(0,1){40}}}
        \put(400,40){\color{darkgreen}{\line(0,-1){40}}}
        \put(360,0){\color{purple}{\line(1,0){40}}}
        \put(400,40){\color{darkgreen}{\line(-1,0){40}}}
        \put(360,0){\color{purple}{\line(1,1){40}}}
        \put(400,0){\color{darkgreen}{\line(-1,1){40}}}
        \put(360,0){{\circle*{6}}}
        \put(360,40){\color{blue}{\circle*{6}}}
        \put(400,0){\color{blue}{\circle*{6}}}
        \put(400,40){\color{blue}{\circle*{6}}}
        \put(397,-11){$*$}
        \put(357,46){$*$}
        \put(397,46){$*$}
                \put(357,19){\line(1,0){6}}
                \put(357,21){\line(1,0){6}}
                \put(397,20){\line(1,0){6}}
                \put(369,7){\line(0,1){6}}
                \put(371,7){\line(0,1){6}}
                \put(390,7){\line(0,1){6}}
            \put(380,37){\line(0,1){6}}
            \put(379,-3){\line(0,1){6}}
            \put(381,-3){\line(0,1){6}}

                \put(350,-22){\footnotesize $\Gamma_{11}=\langle (234) \rangle$}
        \put(375,-40){\footnotesize $\mathcal{G}_{17}$}

                \put(431,0){1}
        \put(431,40){4}
        \put(491,0){2}
        \put(491,40){3}
        \put(445,0){\color{darkgreen}{\line(0,1){40}}}
        \put(485,40){\color{darkgreen}{\line(0,-1){40}}}
        \put(445,0){{\line(1,0){40}}}
        \put(485,40){{\line(-1,0){40}}}
        \put(445,0){\color{purple}{\line(1,1){40}}}
        \put(485,0){\color{purple}{\line(-1,1){40}}}
        \put(445,0){\color{blue}{\circle*{6}}}
        \put(445,40){\color{red}{\circle*{6}}}
        \put(485,0){\color{blue}{\circle*{6}}}
        \put(485,40){\color{red}{\circle*{6}}}
        \put(442,-11){$*$}
        \put(482,-11){$*$}
        \put(440,46){$**$}
        \put(480,46){$**$}
                \put(442,20){\line(1,0){6}}
                \put(482,20){\line(1,0){6}}
                \put(454,7){\line(0,1){6}}
                \put(456,7){\line(0,1){6}}
                \put(476,7){\line(0,1){6}}
                \put(474,7){\line(0,1){6}}

                \put(430,-22){\footnotesize $\Gamma_{12}=\langle (12)(34) \rangle$}

        \put(6,-90){1}
        \put(6,-50){4}
        \put(66,-90){2}
        \put(66,-50){3}
        \put(20,-90){\color{purple}{\line(0,1){40}}}
        \put(60,-50){\color{purple}{\line(0,-1){40}}}
        \put(20,-90){\color{darkgreen}{\line(1,0){40}}}
        \put(60,-50){\color{darkgreen}{\line(-1,0){40}}}
        \put(20,-90){{\line(1,1){40}}}
        \put(60,-90){{\line(-1,1){40}}}
        \put(20,-90){\color{blue}{\circle*{6}}}
        \put(20,-50){\color{red}{\circle*{6}}}
        \put(60,-90){\color{red}{\circle*{6}}}
        \put(60,-50){\color{blue}{\circle*{6}}}
        \put(17,-101){$*$}
        \put(55,-101){$**$}
        \put(15,-44){$**$}
        \put(57,-44){$*$}
                \put(17,-71){\line(1,0){6}}
                \put(17,-69){\line(1,0){6}}
                \put(57,-71){\line(1,0){6}}
                \put(57,-69){\line(1,0){6}}
                \put(40,-53){\line(0,1){6}}
                \put(40,-93){\line(0,1){6}}

                \put(3,-112){\footnotesize $\Gamma_{13}=\langle (13)(24) \rangle$}

                \put(91,-90){1}
        \put(91,-50){4}
        \put(151,-90){2}
        \put(151,-50){3}
        \put(105,-90){{\line(0,1){40}}}
        \put(145,-50){{\line(0,-1){40}}}
        \put(105,-90){\color{darkgreen}{\line(1,0){40}}}
        \put(145,-50){\color{darkgreen}{\line(-1,0){40}}}
        \put(105,-90){\color{purple}{\line(1,1){40}}}
        \put(145,-90){\color{purple}{\line(-1,1){40}}}
        \put(105,-90){\color{blue}{\circle*{6}}}
        \put(105,-50){\color{blue}{\circle*{6}}}
        \put(145,-90){\color{red}{\circle*{6}}}
        \put(145,-50){\color{red}{\circle*{6}}}
        \put(102,-101){$*$}
        \put(140,-101){$**$}
        \put(102,-44){$*$}
        \put(140,-44){$**$}
                \put(114,-83){\line(0,1){6}}
                \put(116,-83){\line(0,1){6}}
                \put(134,-83){\line(0,1){6}}
                \put(136,-83){\line(0,1){6}}
                \put(125,-93){\line(0,1){6}}
                \put(125,-53){\line(0,1){6}}

                \put(88,-112){\footnotesize $\Gamma_{14}=\langle (14)(23) \rangle$}

                \put(176,-90){1}
        \put(176,-50){4}
        \put(236,-90){2}
        \put(236,-50){3}
        \put(190,-90){\color{darkgreen}{\line(0,1){40}}}
        \put(230,-50){\color{darkgreen}{\line(0,-1){40}}}
        \put(190,-90){\color{darkgreen}{\line(1,0){40}}}
        \put(230,-50){\color{darkgreen}{\line(-1,0){40}}}
        \put(190,-90){\color{purple}{\line(1,1){40}}}
        \put(230,-90){\color{purple}{\line(-1,1){40}}}
        \put(190,-90){\color{blue}{\circle*{6}}}
        \put(190,-50){\color{blue}{\circle*{6}}}
        \put(230,-90){\color{blue}{\circle*{6}}}
        \put(230,-50){\color{blue}{\circle*{6}}}
        \put(187,-101){$*$}
        \put(227,-101){$*$}
        \put(187,-44){$*$}
        \put(227,-44){$*$}
                \put(187,-70){\line(1,0){6}}
                \put(227,-70){\line(1,0){6}}
                \put(199,-83){\line(0,1){6}}
                \put(201,-83){\line(0,1){6}}
                \put(219,-83){\line(0,1){6}}
                \put(221,-83){\line(0,1){6}}
            \put(210,-53){\line(0,1){6}}
            \put(210,-93){\line(0,1){6}}

                \put(178,-112){\footnotesize $\Gamma_{15}=\langle (1234) \rangle$}

                \put(261,-90){1}
        \put(261,-50){4}
        \put(321,-90){2}
        \put(321,-50){3}
        \put(275,-90){\color{purple}{\line(0,1){40}}}
        \put(315,-50){\color{purple}{\line(0,-1){40}}}
        \put(275,-90){\color{darkgreen}{\line(1,0){40}}}
        \put(315,-50){\color{darkgreen}{\line(-1,0){40}}}
        \put(275,-90){\color{darkgreen}{\line(1,1){40}}}
        \put(315,-90){\color{darkgreen}{\line(-1,1){40}}}
        \put(275,-90){\color{blue}{\circle*{6}}}
        \put(275,-50){\color{blue}{\circle*{6}}}
        \put(315,-90){\color{blue}{\circle*{6}}}
        \put(315,-50){\color{blue}{\circle*{6}}}
        \put(272,-101){$*$}
        \put(312,-101){$*$}
        \put(272,-44){$*$}
        \put(312,-44){$*$}
                \put(272,-69){\line(1,0){6}}
                \put(272,-71){\line(1,0){6}}
                \put(312,-69){\line(1,0){6}}
                \put(312,-71){\line(1,0){6}}
                \put(285,-83){\line(0,1){6}}
                \put(305,-83){\line(0,1){6}}
            \put(295,-53){\line(0,1){6}}
            \put(295,-93){\line(0,1){6}}

                \put(263,-112){\footnotesize $\Gamma_{16}=\langle (1243) \rangle$}

                \put(346,-90){1}
        \put(346,-50){4}
        \put(406,-90){2}
        \put(406,-50){3}
        \put(360,-90){\color{darkgreen}{\line(0,1){40}}}
        \put(400,-50){\color{darkgreen}{\line(0,-1){40}}}
        \put(360,-90){\color{purple}{\line(1,0){40}}}
        \put(400,-50){\color{purple}{\line(-1,0){40}}}
        \put(360,-90){\color{darkgreen}{\line(1,1){40}}}
        \put(400,-90){\color{darkgreen}{\line(-1,1){40}}}
        \put(360,-90){\color{blue}{\circle*{6}}}
        \put(360,-50){\color{blue}{\circle*{6}}}
        \put(400,-90){\color{blue}{\circle*{6}}}
        \put(400,-50){\color{blue}{\circle*{6}}}
        \put(357,-101){$*$}
        \put(397,-101){$*$}
        \put(357,-44){$*$}
        \put(397,-44){$*$}
                \put(357,-70){\line(1,0){6}}
                \put(397,-70){\line(1,0){6}}
                \put(370,-83){\line(0,1){6}}
                \put(390,-83){\line(0,1){6}}
            \put(379,-53){\line(0,1){6}}
            \put(381,-53){\line(0,1){6}}
            \put(379,-93){\line(0,1){6}}
            \put(381,-93){\line(0,1){6}}

                \put(348,-112){\footnotesize $\Gamma_{17}=\langle (1324) \rangle$}

        \end{normalsize}
\end{picture}\\
\vspace{90pt}
  \caption{Colourings in $\mathcal{K}_{[4]}$ which
are generated by $\Gamma=\langle \sigma \rangle$ for some $\sigma
\in S(V)$.}
  \label{fig_s4_orbits}
\end{figure}
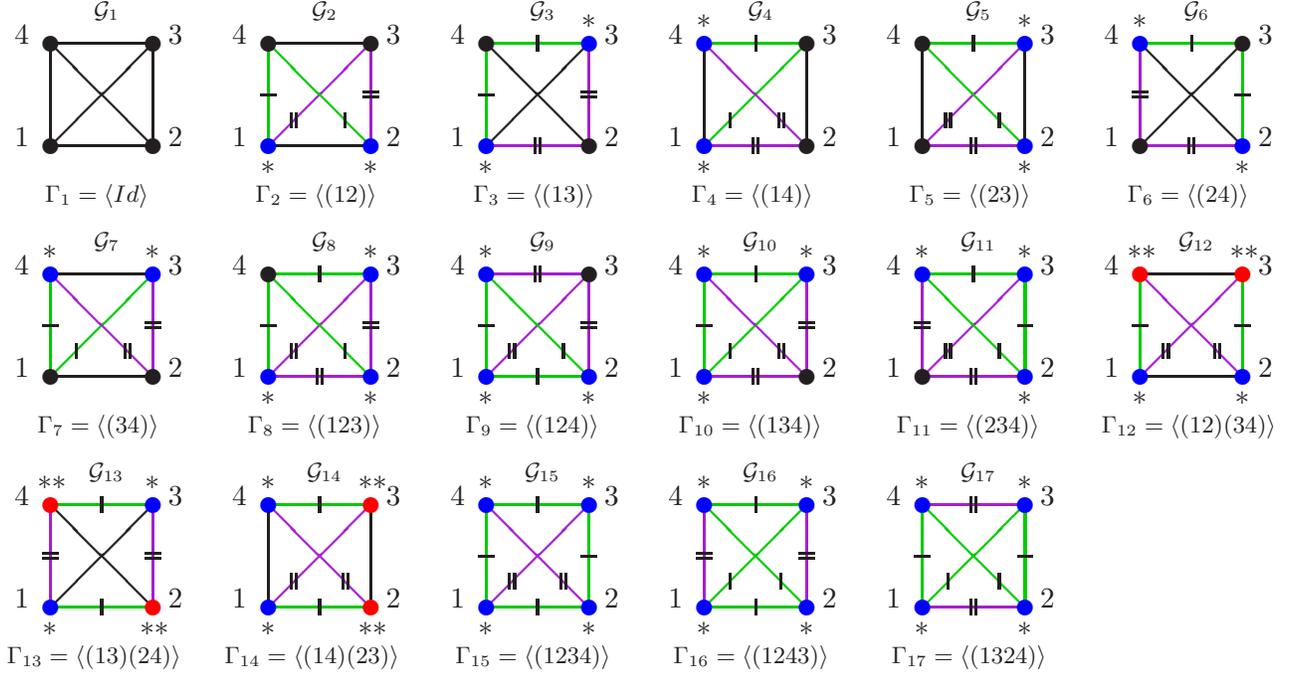

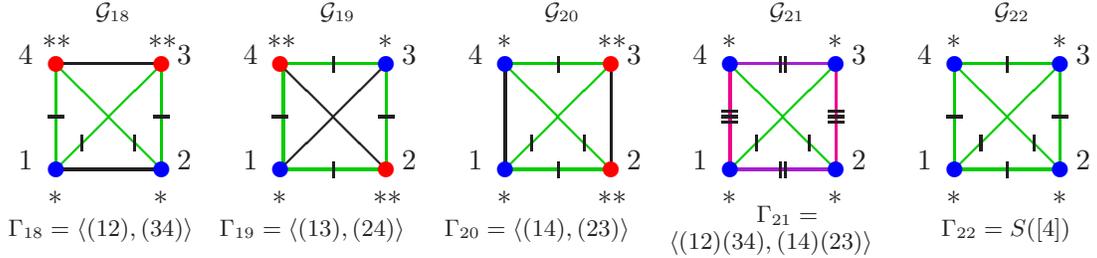
\begin{figure}[h!]
  \centering
  \vspace{90pt}
  \begin{picture}(400,90)(10,-30)
\begin{normalsize}
\thicklines
         \put(6,90){1}
        \put(6,130){4}
        \put(66,90){2}
        \put(66,130){3}
        \put(20,90){\color{darkgreen}{\line(0,1){40}}}
        \put(60,130){\color{darkgreen}{\line(0,-1){40}}}
        \put(20,90){{\line(1,0){40}}}
        \put(60,130){{\line(-1,0){40}}}
        \put(20,90){\color{darkgreen}{\line(1,1){40}}}
        \put(60,90){\color{darkgreen}{\line(-1,1){40}}}
        \put(20,90){\color{blue}{\circle*{6}}}
        \put(20,130){\color{red}{\circle*{6}}}
        \put(60,90){\color{blue}{\circle*{6}}}
        \put(60,130){\color{red}{\circle*{6}}}

        \put(17,77){$*$}
        \put(57,77){$*$}
        \put(15,136){$**$}
        \put(55,136){$**$}

        \put(17,110){\line(1,0){6}}
                \put(57,110){\line(1,0){6}}
                \put(30,97){\line(0,1){6}}
                \put(50,97){\line(0,1){6}}

        \put(35,147){\footnotesize $\mathcal{G}_{18}$}
        \put(2,65){\footnotesize $\Gamma_{18} = \langle(12), (34)\rangle$}

                \put(91,90){1}
        \put(91,130){4}
        \put(151,90){2}
        \put(151,130){3}
        \put(106,90){\color{darkgreen}{\line(0,1){40}}}
        \put(145,130){\color{darkgreen}{\line(0,-1){40}}}
        \put(105,90){\color{darkgreen}{\line(1,0){40}}}
        \put(145,130){\color{darkgreen}{\line(-1,0){40}}}
        \put(105,90){{\line(1,1){40}}}
        \put(145,90){{\line(-1,1){40}}}
        \put(105,90){\color{blue}{\circle*{6}}}
        \put(105,130){\color{red}{\circle*{6}}}
        \put(145,90){\color{red}{\circle*{6}}}
        \put(145,130){\color{blue}{\circle*{6}}}
        \put(102,77){$*$}
        \put(140,77){$**$}
        \put(100,136){$**$}
        \put(142,136){$*$}
                \put(102,110){\line(1,0){6}}
                \put(142,110){\line(1,0){6}}
                \put(125,87){\line(0,1){6}}
                \put(125,127){\line(0,1){6}}

        \put(120,147){\footnotesize $\mathcal{G}_{19}$}
        \put(82,65){\footnotesize $\Gamma_{19} = \langle(13), (24)\rangle$}

                \put(176,90){1}
        \put(176,130){4}
        \put(236,90){2}
        \put(236,130){3}
        \put(190,90){{\line(0,1){40}}}
        \put(230,130){{\line(0,-1){40}}}
        \put(190,90){\color{darkgreen}{\line(1,0){40}}}
        \put(230,130){\color{darkgreen}{\line(-1,0){40}}}
        \put(190,90){\color{darkgreen}{\line(1,1){40}}}
        \put(230,90){\color{darkgreen}{\line(-1,1){40}}}
        \put(190,90){\color{blue}{\circle*{6}}}
        \put(190,130){\color{blue}{\circle*{6}}}
        \put(230,90){\color{red}{\circle*{6}}}
        \put(230,130){\color{red}{\circle*{6}}}
        \put(187,77){$*$}
        \put(225,77){$**$}
        \put(187,136){$*$}
        \put(225,136){$**$}
                \put(200,97){\line(0,1){6}}
                \put(220,97){\line(0,1){6}}
            \put(210,127){\line(0,1){6}}
            \put(210,87){\line(0,1){6}}

        \put(205,147){\footnotesize $\mathcal{G}_{20}$}
        \put(167,65){\footnotesize $\Gamma_{20} = \langle(14), (23)\rangle$}

                \put(261,90){1}
        \put(261,130){4}
        \put(321,90){2}
        \put(321,130){3}
        \put(275,90){\color{magenta}{\line(0,1){40}}}
        \put(315,130){\color{magenta}{\line(0,-1){40}}}
        \put(275,90){\color{purple}{\line(1,0){40}}}
        \put(315,130){\color{purple}{\line(-1,0){40}}}
        \put(275,90){\color{darkgreen}{\line(1,1){40}}}
        \put(315,90){\color{darkgreen}{\line(-1,1){40}}}
        \put(275,90){\color{blue}{\circle*{6}}}
        \put(275,130){\color{blue}{\circle*{6}}}
        \put(315,90){\color{blue}{\circle*{6}}}
        \put(315,130){\color{blue}{\circle*{6}}}
        \put(272,77){$*$}
        \put(312,77){$*$}
        \put(272,136){$*$}
        \put(312,136){$*$}
                \put(285,97){\line(0,1){6}}
                \put(305,97){\line(0,1){6}}
            \put(294,127){\line(0,1){6}}
            \put(296,127){\line(0,1){6}}
            \put(294,87){\line(0,1){6}}
            \put(296,87){\line(0,1){6}}
            \put(272,110){\line(1,0){6}}
            \put(272,108){\line(1,0){6}}
            \put(272,112){\line(1,0){6}}
            \put(312,110){\line(1,0){6}}
            \put(312,108){\line(1,0){6}}
            \put(312,112){\line(1,0){6}}

        \put(290,147){\footnotesize $\mathcal{G}_{21}$}
        \put(285,70){\footnotesize $\Gamma_{21} = $}
        \put(252,60){\footnotesize $\langle(12)(34), (14)(23)\rangle$}

                \put(346,90){1}
        \put(346,130){4}
        \put(406,90){2}
        \put(406,130){3}
        \put(360,90){\color{darkgreen}{\line(0,1){40}}}
        \put(400,130){\color{darkgreen}{\line(0,-1){40}}}
        \put(360,90){\color{darkgreen}{\line(1,0){40}}}
        \put(400,130){\color{darkgreen}{\line(-1,0){40}}}
        \put(360,90){\color{darkgreen}{\line(1,1){40}}}
        \put(400,90){\color{darkgreen}{\line(-1,1){40}}}
        \put(360,90){\color{blue}{\circle*{6}}}
        \put(360,130){\color{blue}{\circle*{6}}}
        \put(400,90){\color{blue}{\circle*{6}}}
        \put(400,130){\color{blue}{\circle*{6}}}
        \put(357,77){$*$}
        \put(397,77){$*$}
        \put(357,136){$*$}
        \put(397,136){$*$}
                \put(357,110){\line(1,0){6}}
                \put(397,110){\line(1,0){6}}
                \put(370,97){\line(0,1){6}}
                \put(390,97){\line(0,1){6}}
            \put(380,127){\line(0,1){6}}
            \put(380,87){\line(0,1){6}}

        \put(375,147){\footnotesize $\mathcal{G}_{22}$}
        \put(355,65){\footnotesize $\Gamma_{22} = S([4])$}

        \end{normalsize}
\end{picture}\\
\vspace{-90pt}
  \caption{Remaining colourings in $\mathcal{K}_{[4]}$.}
  \label{fig_s4_orbits2}
\end{figure}

For $V=\{B_{1}, B_{2},L_{1}, L_{2}\}$ the symmetric group $S(V)$
contains $4!=24$ permutations and has 30 subgroups, 17 of which are
generated by a single permutation. Let $\mathcal{K}_{V}$ denote the
set of colourings of the complete graph these groups generate. The
graph colourings in $\mathcal{K}_{V}$ which are generated by a
single permutation, i.e.\ by $\Gamma = \langle \sigma
\rangle$ for $\sigma \in S(V)$, are displayed in Figure
\ref{fig_s4_orbits}, where, for the sake of legibility we label the
vertices by $V=\{1,2,3,4\}$. The remaining 13 subgroups generate
only 5 distinct colourings in $\mathcal{K}_{[4]}$, all of which are
shown in Figure \ref{fig_s4_orbits2} with one of their generating
groups.

The search space $\Pi_{[4]}$ consists of all models which are
represented by one of the 22 graphs in $\mathcal{K}_{[4]}
= \{\mathcal{G}_{1}, \ldots, \mathcal{G}_{22}\}$ displayed in
Figures \ref{fig_s4_orbits} and \ref{fig_s4_orbits2}, together with
those represented by graphs which can be obtained from the above by
dropping edge colour classes. Thus the size of the total search
space~is
\begin{eqnarray*}
|\Pi_{[4]}|&=& N_{1} + 6 N_{2} + 4 N_{8} + 3 N_{12} + 3(N_{15}-1) + 3(N_{18}-4) + (N_{21}-4) + N_{22}\\
&=& 2^{6} + 6 \cdot 2^{4} + 4 \cdot 2^{2} + 3 \cdot 2^{4} + 3
(2^{2}-1) + 3 (2^{3} -4) + (2^{3}-4) + 2 = 251
\end{eqnarray*}
where $N_{i}$ denotes the number of graphs one can obtain from graph
$\mathcal{G}_{i}$. The subtracted correction terms prevent some
graphs to be counted more than once.

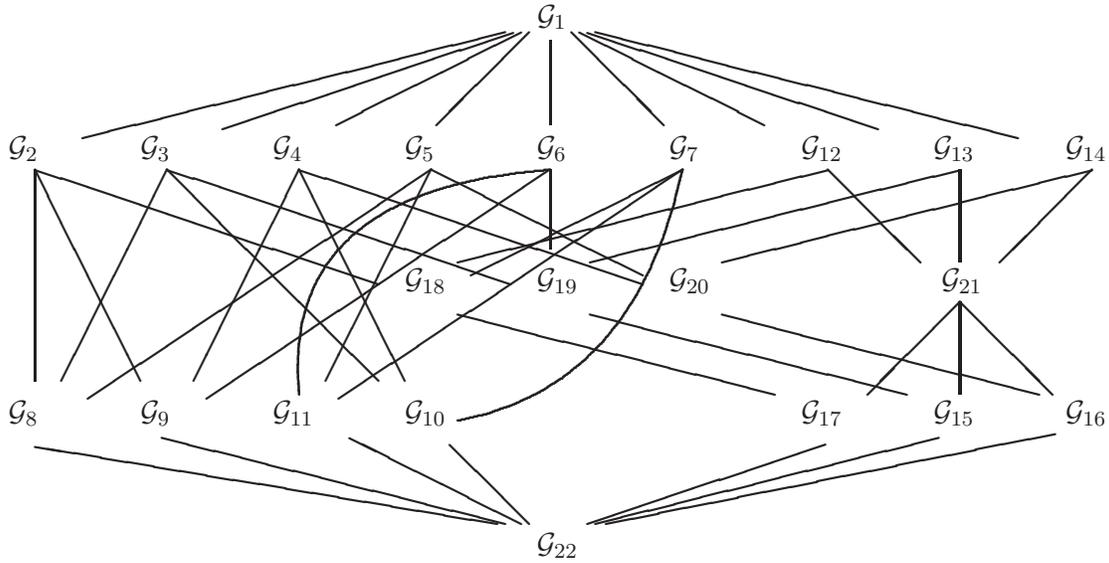
\begin{figure}[h!]
  \centering
\begin{picture}(450,185)(10,10)
\begin{normalsize}
\thicklines

\put(225,180){$\mathcal{G}_{1}$}

\put(213,177){\line(-4,-1){160}} \put(216,177){\line(-3,-1){110}}
\put(219,177){\line(-2,-1){70}} \put(222,177){\line(-1,-1){35}}
\put(230,174){\line(0,-1){32}}

\put(247,177){\line(4,-1){160}} \put(244,177){\line(3,-1){110}}
\put(241,177){\line(2,-1){70}} \put(238,177){\line(1,-1){35}}

\put(25,130){$\mathcal{G}_{2}$} \put(75,130){$\mathcal{G}_{3}$}
\put(125,130){$\mathcal{G}_{4}$} \put(175,130){$\mathcal{G}_{5}$}

\put(35,125){\line(3,-1){130}} \put(85,125){\line(3,-1){130}}
\put(135,125){\line(3,-1){130}} \put(185,125){\line(2,-1){80}}
\put(230,125){\line(0,-1){30}} \put(280,125){\line(-2,-1){80}}

\put(35,125){\line(0,-1){80}} \put(35,125){\line(1,-2){40}}
\put(85,125){\line(-1,-2){40}} \put(85,125){\line(1,-1){80}}
\put(135,125){\line(1,-2){40}} \put(135,125){\line(-1,-2){40}}
\put(185,125){\line(-1,-2){40}} \put(185,125){\line(-3,-2){130}}
\put(230,125){\line(-3,-2){130}} \curve(230,125, 155,100, 135,40)
\put(280,125){\line(-3,-2){130}} \curve(280,125, 250,60, 195,30)

\put(225,130){$\mathcal{G}_{6}$} \put(275,130){$\mathcal{G}_{7}$}
\put(325,130){$\mathcal{G}_{12}$} \put(375,130){$\mathcal{G}_{13}$}
\put(425,130){$\mathcal{G}_{14}$} \put(378,80){$\mathcal{G}_{21}$}

\put(335,125){\line(1,-1){35}} \put(385,125){\line(0,-1){35}}
\put(435,125){\line(-1,-1){35}} \put(385,75){\line(-1,-1){35}}
\put(385,75){\line(0,-1){35}} \put(385,75){\line(1,-1){35}}

\put(335,125){\line(-4,-1){140}} \put(385,125){\line(-4,-1){140}}
\put(435,125){\line(-4,-1){140}}

\put(195,70){\line(4,-1){120}} \put(245,70){\line(4,-1){120}}
\put(295,70){\line(4,-1){120}}

\put(25,30){$\mathcal{G}_{8}$} \put(75,30){$\mathcal{G}_{9}$}
\put(125,30){$\mathcal{G}_{11}$} \put(175,30){$\mathcal{G}_{10}$}

\put(325,30){$\mathcal{G}_{17}$} \put(375,30){$\mathcal{G}_{15}$}
\put(425,30){$\mathcal{G}_{16}$}

\put(175,80){$\mathcal{G}_{18}$} \put(225,80){$\mathcal{G}_{19}$}
\put(275,80){$\mathcal{G}_{20}$}

\put(225,-20){$\mathcal{G}_{22}$}

\put(213,-9){\line(-4,1){130}} \put(210,-9){\line(-6,1){175}}
\put(219,-9){\line(-2,1){65}} \put(222,-9){\line(-1,1){30}}

\put(251,-9){\line(5,1){170}} \put(247,-9){\line(4,1){130}}
\put(244,-9){\line(3,1){90}}

        \end{normalsize}
\end{picture}\\
\vspace{40pt}
  \caption{Hasse diagram of $\mathcal{K}_{[4]}$.}
  \label{fig_s4_lattice}
\end{figure}

Figure \ref{fig_s4_lattice} displays the Hasse diagram of
$\mathcal{K}_{[4]}$. By construction, it contains only graphs which
represent the saturated model. An Edwards--Havr{\'a}nek model
selection procedure searches along the full Hasse diagram of all
models in $\Pi_{[4]}$, which contains all 251 models and has the
diagram in Figure \ref{fig_s4_lattice} as a subgraph. At each stage,
the search moves along the edges in the diagram, passing each model
at most once. Once a model has been rejected, all models below it
are excluded from the future search; once a model has been accepted
all models above it are excluded.

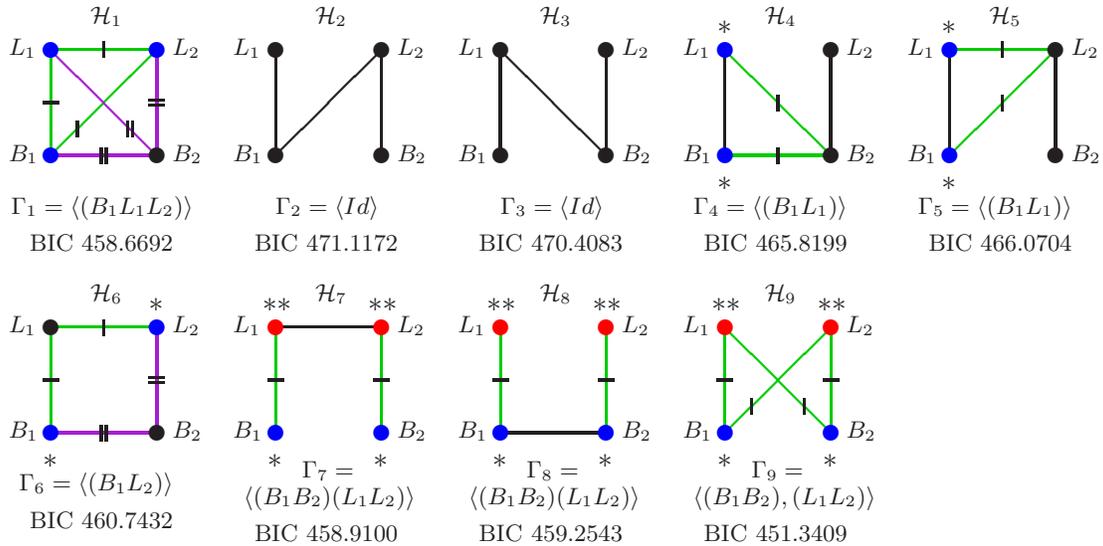
\begin{figure}[ht!]
  \centering
\begin{picture}(400,170)(10,-30)
\begin{normalsize}
\thicklines

                \put(4,88){\footnotesize $B_{1}$}
        \put(4,128){\footnotesize $L_{1}$}
        \put(66,88){\footnotesize $B_{2}$}
        \put(66,128){\footnotesize $L_{2}$}
        \put(20,90){\color{darkgreen}{\line(0,1){40}}}
        \put(60,130){\color{purple}{\line(0,-1){40}}}
        \put(20,90){\color{purple}{\line(1,0){40}}}
        \put(60,130){\color{darkgreen}{\line(-1,0){40}}}
        \put(20,90){\color{darkgreen}{\line(1,1){40}}}
        \put(60,90){\color{purple}{\line(-1,1){40}}}
        \put(20,90){\color{blue}{\circle*{6}}}
        \put(20,130){\color{blue}{\circle*{6}}}
        \put(60,90){{\circle*{6}}}
        \put(60,130){\color{blue}{\circle*{6}}}

        \put(17,110){{\line(1,0){6}}}
        \put(57,109){{\line(1,0){6}}}
        \put(57,111){{\line(1,0){6}}}
        \put(40,127){{\line(0,1){6}}}
        \put(41,87){{\line(0,1){6}}}
        \put(39,87){{\line(0,1){6}}}
        \put(30,97){{\line(0,1){6}}}
        \put(51,97){{\line(0,1){6}}}
        \put(49,97){{\line(0,1){6}}}

        \put(35,140){\footnotesize $\mathcal{H}_{1}$}
        \put(5,68){\footnotesize $\Gamma_{1}=\langle (B_{1}L_{1}L_{2}) \rangle$}
        \put(12,54){\footnotesize BIC 458.6692}

                \put(89,88){\footnotesize {$B_{1}$}}
        \put(89,128){\footnotesize {$L_{1}$}}
        \put(151,88){\footnotesize {$B_{2}$}}
        \put(151,128){\footnotesize {$L_{2}$}}
        \put(105,90){{\line(0,1){40}}}
        \put(145,130){{\line(0,-1){40}}}
        \put(105,90){{\line(1,1){40}}}
        \put(105,90){{\circle*{6}}}
        \put(105,130){{\circle*{6}}}
        \put(145,90){{\circle*{6}}}
        \put(145,130){{\circle*{6}}}

                \put(120,140){\footnotesize {$\mathcal{H}_{2}$}}
        \put(105,68){\footnotesize {$\Gamma_{2}=\langle Id\rangle$}}
        \put(97,54){\footnotesize BIC 471.1172}

                \put(174,88){\footnotesize $B_{1}$}
        \put(174,128){\footnotesize $L_{1}$}
        \put(236,88){\footnotesize $B_{2}$}
        \put(236,128){\footnotesize $L_{2}$}
        \put(190,90){{\line(0,1){40}}}
        \put(230,130){{\line(0,-1){40}}}
        \put(230,90){{\line(-1,1){40}}}
        \put(190,90){{\circle*{6}}}
        \put(190,130){{\circle*{6}}}
        \put(230,90){{\circle*{6}}}
        \put(230,130){\circle*{6}}

            \put(205,140){\footnotesize $\mathcal{H}_{3}$}
        \put(190,68){\footnotesize {$\Gamma_{3}=\langle Id\rangle$}}
        \put(182,54){\footnotesize BIC 470.4083}

                \put(259,88){\footnotesize {$B_{1}$}}
        \put(259,128){\footnotesize {$L_{1}$}}
        \put(321,88){\footnotesize {$B_{2}$}}
        \put(321,128){\footnotesize {$L_{2}$}}
        \put(275,90){{\line(0,1){40}}}
        \put(315,130){{\line(0,-1){40}}}
        \put(275,90){\color{darkgreen}{\line(1,0){40}}}
        \put(315,90){\color{darkgreen}{\line(-1,1){40}}}
        \put(275,90){\color{blue}{\circle*{6}}}
        \put(275,130){\color{blue}{\circle*{6}}}
        \put(315,90){{\circle*{6}}}
        \put(315,130){{\circle*{6}}}

                \put(295,87){{\line(0,1){6}}}
                \put(295,107){{\line(0,1){6}}}
                \put(272,136){$*$}
                \put(272,77){$*$}

            \put(290,140){\footnotesize $\mathcal{H}_{4}$}
        \put(263,68){\footnotesize {$\Gamma_{4}=\langle (B_{1}L_{1})\rangle$}}
        \put(267,54){\footnotesize BIC 465.8199}

                \put(344,88){\footnotesize {$B_{1}$}}
        \put(344,128){\footnotesize {$L_{1}$}}
        \put(406,88){\footnotesize {$B_{2}$}}
        \put(406,128){\footnotesize {$L_{2}$}}
        \put(360,90){{\line(0,1){40}}}
        \put(400,130){{\line(0,-1){40}}}
        \put(400,130){\color{darkgreen}{\line(-1,0){40}}}
        \put(360,90){\color{darkgreen}{\line(1,1){40}}}
        \put(360,90){\color{blue}{\circle*{6}}}
        \put(360,130){\color{blue}{\circle*{6}}}
        \put(400,90){{\circle*{6}}}
        \put(400,130){{\circle*{6}}}
        \put(357,136){$*$}
        \put(357,77){{$*$}}
            \put(380,127){{\line(0,1){6}}}
            \put(380,107){{\line(0,1){6}}}

            \put(375,140){\footnotesize {$\mathcal{H}_{5}$}}
        \put(348,68){\footnotesize {$\Gamma_{5}=\langle (B_{1}L_{1})\rangle$}}
        \put(352,54){\footnotesize BIC 466.0704}

                \put(4,-17){{\footnotesize $B_{1}$}}
        \put(4,23){{\footnotesize $L_{1}$}}
        \put(66,-17){{\footnotesize $B_{2}$}}
        \put(66,23){{\footnotesize $L_{2}$}}
        \put(20,-15){\color{darkgreen}{\line(0,1){40}}}
        \put(60,25){\color{purple}{\line(0,-1){40}}}
        \put(20,-15){\color{purple}{\line(1,0){40}}}
        \put(60,25){\color{darkgreen}{\line(-1,0){40}}}
        \put(20,-15){\color{blue}{\circle*{6}}}
        \put(20,25){{\circle*{6}}}
        \put(60,-15){{\circle*{6}}}
        \put(60,25){\color{blue}{\circle*{6}}}

        \put(17,-28){{$*$}}
        \put(57,31){$*$}

            \put(41,-18){{\line(0,1){6}}}
            \put(39,-18){{\line(0,1){6}}}
            \put(40,22){{\line(0,1){6}}}
            \put(17,5){{\line(1,0){6}}}
            \put(57,4){{\line(1,0){6}}}
            \put(57,6){{\line(1,0){6}}}

        \put(35,35){\footnotesize {$\mathcal{H}_{6}$}}
        \put(8,-37){\footnotesize {$\Gamma_{6}=\langle (B_{1}L_{2}) \rangle$}}
        \put(12,-51){\footnotesize BIC 460.7432}

                \put(89,-17){\footnotesize {$B_{1}$}}
        \put(89,23){\footnotesize {$L_{1}$}}
        \put(151,-17){\footnotesize {$B_{2}$}}
        \put(151,23){\footnotesize {$L_{2}$}}
        \put(105,-15){\color{darkgreen}{\line(0,1){40}}}
        \put(145,25){\color{darkgreen}{\line(0,-1){40}}}
        \put(145,25){{\line(-1,0){40}}}
        \put(105,-15){\color{blue}{\circle*{6}}}
        \put(105,25){\color{red}{\circle*{6}}}
        \put(145,-15){\color{blue}{\circle*{6}}}
        \put(102,-28){$*$}
        \put(142,-28){$*$}
                \put(100,31){$**$}
                \put(140,31){$**$}
        \put(145,25){\color{red}{\circle*{6}}}
                \put(142,5){{\line(1,0){6}}}
                \put(102,5){{\line(1,0){6}}}

                \put(120,35){\footnotesize {$\mathcal{H}_{7}$}}
        \put(115,-32){\footnotesize {$\Gamma_{7}=$}}
        \put(93,-42){\footnotesize {$\langle (B_{1}B_{2})(L_{1}L_{2})\rangle$}}
        \put(97,-56){\footnotesize BIC 458.9100}

                \put(174,-17){\footnotesize $B_{1}$}
        \put(174,23){\footnotesize $L_{1}$}
        \put(236,-17){\footnotesize $B_{2}$}
        \put(236,23){\footnotesize $L_{2}$}
        \put(190,-15){\color{darkgreen}{\line(0,1){40}}}
        \put(230,25){\color{darkgreen}{\line(0,-1){40}}}
        \put(190,-15){{\line(1,0){40}}}
        \put(190,-15){\color{blue}{\circle*{6}}}
        \put(190,25){\color{red}{\circle*{6}}}
        \put(230,-15){\color{blue}{\circle*{6}}}
        \put(230,25){\color{red}{\circle*{6}}}
        \put(187,-28){$*$}
        \put(227,-28){$*$}
                \put(185,31){$**$}
                \put(225,31){$**$}
                \put(227,5){{\line(1,0){6}}}
                \put(187,5){{\line(1,0){6}}}

                \put(205,35){\footnotesize {$\mathcal{H}_{8}$}}
        \put(200,-32){\footnotesize {$\Gamma_{8}=$}}
        \put(178,-42){\footnotesize {$\langle (B_{1}B_{2})(L_{1}L_{2})\rangle$}}
        \put(182,-56){\footnotesize BIC 459.2543}

                \put(259,-17){\footnotesize {$B_{1}$}}
        \put(259,23){\footnotesize {$L_{1}$}}
        \put(321,-17){\footnotesize {$B_{2}$}}
        \put(321,23){\footnotesize {$L_{2}$}}
        \put(275,-15){\color{darkgreen}{\line(0,1){40}}}
        \put(315,25){\color{darkgreen}{\line(0,-1){40}}}
        \put(275,-15){\color{darkgreen}{\line(1,1){40}}}
        \put(315,-15){\color{darkgreen}{\line(-1,1){40}}}
        \put(275,-15){\color{blue}{\circle*{6}}}
        \put(275,25){\color{red}{\circle*{6}}}
        \put(315,-15){\color{blue}{\circle*{6}}}
        \put(315,25){\color{red}{\circle*{6}}}
        \put(272,-28){$*$}
        \put(312,-28){$*$}
                \put(270,31){$**$}
                \put(310,31){$**$}
                \put(272,5){{\line(1,0){6}}}
                \put(312,5){{\line(1,0){6}}}
                \put(285,-8){{\line(0,1){6}}}
                \put(305,-8){{\line(0,1){6}}}

            \put(290,35){\footnotesize $\mathcal{H}_{9}$}
        \put(285,-32){\footnotesize $\Gamma_{9}=$}
        \put(263,-42){\footnotesize $\langle (B_{1}B_{2}), (L_{1}L_{2})\rangle$}
        \put(267,-56){\footnotesize BIC 451.3409}

                \end{normalsize}
        \end{picture}
        \vspace{25pt}
 \caption{Graphs of minimally accepted models in $\Pi_{[4]}$ for Frets' heads data.}
  \label{fig_frets_min}
\end{figure}

Exploiting the demonstrated lattice structure of $\Pi_{[4]}$, we
applied the algorithm to the Fret's heads data described in Example
\ref{rcop_example} with $\mathcal{A}$ initially consisting of the
saturated unrestricted model and $\mathcal{R}=\emptyset$. After
testing 48 models in 4 stages the algorithm arrived at 9 minimally
accepted models whose graphs and generating groups are given in
Figure \ref{fig_frets_min}. (The models are distributed between the
stages the following way. 1: 15 (9 accepted), 2: 16 (16 accepted),
3: 13 (13 accepted), 4: 4 (3 accepted).)

The minimally accepted model with the lowest BIC value is
represented by graph $\mathcal{H}_{9}$, which is considerably less
than the BIC value 471.2982 of the model fitted
in~\cite{lauritzen_sym} whose graph is displayed in Figure
\ref{fig:rcop_four_cycle}. Further, the model selected
in~\cite{lauritzen_sym} is a supermodel, in fact the supremum, of
the models represented by $\mathcal{H}_{7}$ and $\mathcal{H}_{8}$,
with a further edge to complete the four cycle. Interestingly,
$\mathcal{H}_{2}$ and $\mathcal{H}_{3}$ are two of the graphs found
in Section 8.3 in \citet{whittaker}, the other two being the
underlying uncoloured graphs of $\mathcal{H}_{7}$ and
$\mathcal{H}_{8}$. Note that the BIC value of the complete symmetry
model, whose graph is displayed on the right in Figure
\ref{fig:rcop_four_cycle}, lies between the smallest and the second
smallest BIC values of the minimally accepted models, the
corresponding graphs being $\mathcal{H}_{9}$ and $\mathcal{H}_{1}$
respectively. However it was (weakly) rejected by the procedure as
it is a submodel of a model rejected in stage 1.

\section{Discussion}

As we argued, graphical Gaussian models with equality constraints
are a promising model class as they combine parsimony in the number
of parameters with the concise and efficient graphical models
framework. We studied two model types introduced
by \citet{lauritzen_sym}: RCON models which place equality 
restrictions on the model covariance matrix $\Sigma^{-1}$ and RCOR
models which restrict the diagonal of $\Sigma^{-1}$ and the partial
correlations, which can both be represented by vertex and edge
coloured graphs.

We showed four model classes within the sets of RCON and RCOR
models, each possessing desirable statistical properties and being
more readily interpretable than RCON and RCOR models in general, to
form complete non-distributive lattices. This qualifies each of them
for an Edwards--Havr{\'a}nek model selection procedure. Two model
classes, those represented by edge regular and permutation-generated
colourings respectively, are most readily interpretable and possess
the most tractable structure out of the four and are thus most
suitable for a model search.

For the former model class we have developed an
Edwards--Havr{\'a}nek model selection algorithm, and demonstrated an
encouraging performance for the data set previously described in
Example \ref{rcon_example}. We further illustrated the principal
functionality of the Edwards--Havr{\'a}nek procedure on the lattice
of models represented by permutation-generated colourings with the
example of Fret's heads data from Example \ref{rcop_example}. Here
as well, the algorithm performed in a satisfactory fashion. In order
to fully generalise it to work for any number of variables $|V|$,
further investigation into the relationship between permutation
groups acting on $V$ and their orbits in $V$ and $V \times V$ is
necessary.

Some potential concerns need to be mentioned: firstly, while the
algorithm's performance in the above examples was encouraging, it
has to be taken into account that the number of variables is rather
small in both cases. Further, it is at this stage unknown how much
this behaviour relied on strong/weak conditional independence and
symmetry relations in the data sets considered. It may be that the
number of models to be tested can still grow in an unmanageable
fashion. Secondly, a general concern with the Edwards--Havr{\'a}nek
model selection procedure is that its sampling properties are
intractable. In particular the procedure does not control the
overall error rate.


Especially in view of the above, it may be
worthwhile to explore alternative model selection approaches. We
argued that neighbourhood selection with the
lasso \citep{meinshausen2006}, stability
selection \citep{meinshausen2010}, and the SINful
approach \citep{drton2008} were not directly applicable to the
lattices of models studied in this article due to their
non-distributivity. Modified variants may still be feasible, which
could be investigated.

One further viable alternative may be a symmetry variant of the
graphical lasso \citep{glasso, wainwright2008}, which in its original
form seeks to maximise the penalised log-likelihood
\begin{equation} \label{lasso}
\log \det \Sigma^{-1} - \textup{tr}(S\Sigma^{-1}) - \rho \|
\Sigma^{-1} \|_{1}
\end{equation}
over non-negative definite matrices $\Sigma^{-1}$, with $S$ denoting
the empirical covariance matrix of the observations,
$\textup{tr}(\cdot)$ being the trace, $\| \cdot \|_{1}$ the $l_{1}$
norm giving the sum of the absolute values of the elements in the
argument matrix and $\rho$ being the penalisation parameter. Testing
for equality constraints on the entries of
$\Sigma^{-1}=(k_{\alpha\beta})_{\alpha, \beta \in V}$ can be
enforced by replacing equation (\ref{lasso}) by the following
function
\begin{equation*}
\log \det \Sigma^{-1} - \textup{tr}(S\Sigma^{-1}) - \rho_{1} \|
\Sigma^{-1} \|_{1} - \rho_{2} \sum_{\alpha, \beta \in
V}|k_{\alpha\alpha} - k_{\beta\beta}| - \rho_{3}
\sum_{\substack{\alpha,\beta, \gamma,\delta \in V,
\\\alpha\not=\beta, \gamma\not=\delta}}|k_{\alpha\beta} -
k_{\gamma\delta}|
\end{equation*}
This lies in direct analogy to the development of the fused
lasso \citep{fused_lasso} from the standard lasso for linear
regression. However how to maximise the function for a given model
class, for example over models with edge regular colourings to
ensure scale invariance, seems non-trivial.

\section*{Acknowledgments}

I would like to thank my PhD supervisor Steffen Lauritzen for many
helpful discussions on the topic and Richard Samworth for inspiring
conversations. I am very grateful to Michael Perlman and an
anonymous referee for their constructive criticism of an earlier
version of this article.

\bibliographystyle{chicago}
\bibliography{lattices_gehrmann}

\end{document}